\documentclass[12pt,a4paper,reqno]{amsart}
\usepackage{amsmath}
\usepackage{amsfonts}
\usepackage{amssymb}
\usepackage{amscd}
\numberwithin{equation}{section}

\theoremstyle{plain}

\newtheorem{theorem}[subsection]{Theorem}
\newtheorem{proposition}[subsection]{Proposition}
\newtheorem{lemma}[subsection]{Lemma}
\newtheorem{corollary}[subsection]{Corollary}

\newtheorem*{main-theorem-1}{Main Theorem}
\newtheorem*{nil-to-local-repeat}{Proposition \ref{nil-to-local}}
\newtheorem*{restate-tech-thm}{Theorem \ref{mainthm-tech}}

\theoremstyle{definition}

\newtheorem{definition}[subsection]{Definition}

\theoremstyle{remark}

\newtheorem{example}{Example}
\newtheorem{examples}{Examples}
\newtheorem*{remark}{Remark}
\newtheorem*{remarks}{Remarks}

\renewcommand{\leq}{\leqslant}
\renewcommand{\geq}{\geqslant}

\newsavebox{\proofbox}
\savebox{\proofbox}{\begin{picture}(7,7)%
  \put(0,0){\framebox(7,7){}}\end{picture}}
\def\boxeq{\tag*{\usebox{\proofbox}}}

\newcommand{\md}[1]{\ensuremath{(\operatorname{mod}\, #1)}} 

\def\endproof{\hfill{\usebox{\proofbox}}}
\def\E{\mathbb{E}}
\def\Z{\mathbb{Z}}
\def\R{\mathbb{R}}

\def\C{\mathbb{C}}
\def\N{\mathbb{N}}

\def\b{{\mathbf b}}

\def\sgn{\operatorname{sgn}}

\def\th{{\operatorname{th}}}
\def\Lip{\operatorname{Lip}}
\def\TV{\operatorname{TV}}

\def\eps{\varepsilon}

\def\half{\textstyle \frac{1}{2}\displaystyle}
\def\typeI{T_{\operatorname{I}}}
\def\typeII{T_{\operatorname{II}}}
\def\righti{\rangle_{\operatorname{I}}}
\def\rightii{\rangle_{\operatorname{II}}}
\def\Supp{\operatorname{Supp}}

\begin{document}

\title{Quadratic uniformity of the M\"obius function
}

\author{Ben Green}
\address{School of Mathematics, University of Bristol, University Walk, Bristol BS8 1TW.
}
\email{b.j.green@bristol.ac.uk}

\author{Terence Tao}
\address{UCLA Department of Mathematics, Los Angeles, CA 90095-1596.
}
\email{tao@math.ucla.edu}

\thanks{The first author is a Clay Research Fellow and gratefully acknowledges the support of the Clay Institute. He also spent time, while this work was being carried out, at Trinity College, Cambridge and at the Massachusetts Institute of Technology, and is very happy to acknowledge the kind hospitality of both institutions. The second author is supported by a grant from the Packard Foundation.}

\begin{abstract} 
This paper is a part of our programme to generalise the Hardy-Littlewood method to handle systems of linear questions in primes. This programme is laid out in our paper \emph{Linear equations in primes \cite{green-tao-linearprimes}.} In particular, the results of this paper may be used, together with the machinery of \cite{green-tao-linearprimes}, to establish an asymptotic for the number of four-term progressions $p_1 < p_2 < p_3 < p_4 \leq N$ of primes, and more generally any problem counting prime points inside a ``non-degenerate'' affine lattice of codimension at most $2$.

The main result of this paper is a proof of the \emph{M\"obius and Nilsequences Conjecture} for $1$ and $2$-step nilsequences. This conjecture is introduced in \cite{green-tao-linearprimes} and amounts to showing that if $G/\Gamma$ is an $s$-step nilmanifold, $s \leq 2$, if $F : G/\Gamma \rightarrow [-1,1]$ is a Lipschitz function, and if $T_g : G/\Gamma \rightarrow G/\Gamma$ is the action of $g \in G$ on $G/\Gamma$, then
\[ N^{-1}\sum_{n \leq N} \mu(n) F(T_g^n \cdot x) \ll_{A,G/\Gamma} \Vert F \Vert_{\Lip} \log^{-A} N\] uniformly in $g \in G$ and $x \in G/\Gamma$, for any $A > 0$.  This can be viewed as a ``quadratic'' generalisation of an exponential sum estimate of Davenport \cite{davenport-early}, and is proven by following the methods of Vinogradov and Vaughan.
\end{abstract}

\maketitle

\section{Introduction}

 The \emph{M\"obius function} $\mu: \N \to \{-1,0,+1\}$, defined by
\[ \mu(n) := \left\{ \begin{array}{ll} (-1)^{k} & \mbox{if $n = p_1p_2 \dots p_k$ for distinct primes $p_1,\dots,p_k$}\\ 0 & \mbox{if $n$ is not squarefree} \\1, & \mbox{if $n = 1$}\end{array}\right.\]
plays a fundamental role in analytic number theory, especially with regard to the distribution of primes.  A well-known metaprinciple holds that $\mu$ fluctuates so ``randomly'' that it is asymptotically orthogonal to any ``low complexity'' bounded sequence $f : \N \rightarrow \C$. We do not have a formal definition of ``low complexity'', but the examples of this section should convey the general flavour. Functions which arise from geometry and algebra, such as characters $n \mapsto e(n\alpha)$, are certainly of low complexity, whereas functions which depend on the prime factorisation of $n$, such as $\mu$ itself, the von Mangoldt function $\Lambda$, and certain divisor sums arising in sieve theory, are not.

In our first example, and throughout the paper, we will use the following notation. We write $[N] := \{1,\ldots,N\}$ to denote the integers from $1$ to $N$, and $\E_{n \in A} f(n) := \frac{1}{|A|} \sum_{n \in A} f(n)$ to denote the average of a function $f: A \to \C$ on a non-empty finite set $A$.  We also use $X \ll Y$ or $X = O(Y)$ to denote the claim that $|X| \leq CY$ for some absolute constant $C > 0$.

\begin{example}[$\mu$ is strongly orthogonal to the constant function]\label{ex1} We have \begin{equation}\label{mu-1}
\E_{n \in [N]} \mu(n) \ll e^{-c\sqrt{\log N}}
\end{equation}
for all $N > 1$ and some absolute constant $c > 0$. 
\end{example}

\begin{remark} This is essentially equivalent to the prime number theorem with the classical error term of Hadamard and de la Vall\'ee Poussin.\end{remark}

 In the next example, and throughout the paper, we use $X \ll_A Y$ or $X = O_A(Y)$ to denote the claim that $|X| \leq C_A Y$ for some constant $C_A > 0$ depending on $A$.  

\begin{example}[$\mu$ is strongly orthogonal to Dirichlet characters]\label{ex2} For any $A > 0$ we have
\begin{equation}\label{mu-chi} \E_{n \in [N]} \mu(n) \overline{\chi(n)} \ll_A q^{1/2} \log^{-A} N\end{equation} for all $N$ and all Dirichlet characters $\chi$ to modulus $q$. 
\end{example}

\begin{remark}
See for instance \cite[Corollary 5.29]{iwaniec-kowalski}.
This may be used to prove the Siegel-Walfisz theorem concerning the distribution of primes in arithmetic progressions. 
\end{remark}

 The form of the bound in \eqref{mu-chi} may appear strange at first sight. A key point to appreciate is that the implied constant $C = C_A$ is \textit{ineffective}, due to the possible existence of Landau-Siegel zeros. The book \cite{davenport} may be consulted for further information. It is useful to have a name for bounds of this kind.

\begin{definition}[Strong asymptotic orthogonality]
If $f: \N \to \C$ and $g: \N \to \C$ are two sequences on the natural numbers $\N = \{ 1, 2, 3, \ldots \}$, we say that $f$ and $g$ are \emph{strongly asymptotically orthogonal} if we have the estimate
\[ \E_{n \in [N]} f(n) \overline{g(n)} \ll_A \log^{-A} N \]
for all $N > 1$ and all $A > 0$. We allow the implied constant $C_A$ to be \emph{ineffective}, in that we may have no explicit bounds on $C_A$ other than that it is finite.
\end{definition}

Thus Example \ref{ex2} shows that $\mu$ is strongly asymptotically orthogonal to all Dirichlet characters, and some Fourier analysis then shows that it is in fact strongly asymptotically orthogonal to any periodic sequence.  In fact, more is true, as we shall see in the next example. Here, and throughout the paper, we use $e()$ to denote the standard character $e(x) := \exp(2\pi i x)$.

\begin{example}[$\mu$ is strongly orthogonal to linear phases]\label{ex3}
For any $\alpha \in \R/\Z$ and for any $A > 0$, we have
\begin{equation}\label{mu-alpha} \E_{n \in [N]} \mu(n) e(-\alpha n) \ll_A \log^{-A} N,\end{equation} uniformly in $\alpha \in \R/\Z$.\end{example}

This bound is due to Davenport \cite{davenport-early} and can be deduced from \eqref{mu-chi} by an application of Vinogradov's version of the Hardy-Littlewood major/minor arc decomposition of $\R/\Z$. See, for example, \cite[Theorem 13.10]{iwaniec-kowalski}. 
For pedagogical reasons, and because we need this result for later sections, we give the derivation in \S \ref{linear-sec}. Davenport's result may be used on its own to obtain a number of self-correlation estimates on $\mu$.  For instance, by combining \eqref{mu-alpha} with elementary Fourier analysis (the circle method)
we easily obtain the estimates
\begin{equation}\label{mu1} \E_{x,d \in [N]} \mu(x) \mu(x+d) \mu(x+2d) \ll_A \log^{-A} N
\end{equation}
and
\begin{equation}\label{mu2}
 \E_{x,h_1,h_2 \in [N]} \mu(x) \mu(x+h_1) \mu(x+h_2) \mu(x+h_1+h_2) \ll_A \log^{-A} N. 
\end{equation}
Similar expressions in which $\mu$ is replaced by $\Lambda$, the von Mangoldt function, may be analysed using \eqref{mu-alpha} as a key ingredient. The answers have a more complicated form involving a main term which is a product of local factors or \emph{singular series}. See \cite[\S 13]{iwaniec-kowalski} and \cite{green-tao-linearprimes} for different approaches to this\footnote{While the von Mangoldt function $\Lambda$ is more directly related to the primes, the M\"obius function $\mu$ is somewhat easier to handle analytically, being bounded by $1$ and not encountering the ``local'' irregularities in small residue classes that $\Lambda$ faces; in particular, the ``major arc'' terms will have a significantly simpler form.  Also, the Vaughan identity for $\mu$ is slightly cleaner than that for $\Lambda$ (see Lemma \ref{Vaughan}).  Thus in this series of papers we have adopted a ``M\"obius first'' philosophy, in which we obtain estimates on the M\"obius function $\mu$ using ``hard'' analytic tools, and then use ``softer'' techniques to transfer the bounds on $\mu$ to the bounds on $\Lambda$.}.

A full discussion of results such as \eqref{mu1}, \eqref{mu2} and the corresponding results for $\Lambda$ is given in \cite{green-tao-linearprimes}. For comparison with that paper, we remark that the two systems of linear forms in \eqref{mu1} and \eqref{mu2}, namely $(x,x+d,x+2d)$ and $(x, x+ h_1, x+ h_2, x+ h_1 + h_2)$, both have \emph{complexity} equal to one.  This notion of complexity 1 essentially marks the limit of the classical Hardy-Littlewood circle method. The main goal of this paper is to provide some of the technical machinery needed to address the case of complexity 2.

We can reformulate \eqref{mu-alpha} in a manner which may appear strange at first, but is well suited to generalisations, as we shall soon see.  If $X$ is any metric space, define a \emph{Lipschitz function}\footnote{The Lipschitz class is a convenient regularity class for us to use; it is smooth enough that one approximate uniformly and quantitatively by trigonometric series (see Lemma \ref{fourier-lip}), yet rough enough that one can easily extend a function in this class from a small domain to a larger domain (see Lemma \ref{lip-extend}).  Also, the Lipschitz class is meaningful in both discrete and continuous settings.  Of course, the results of this paper also hold in smoother classes such as $C^\infty$, and qualitative versions of these results (with decay factors such as $\log^{-A} N$ replaced by $o(1)$) hold for rougher classes such as the continuous class $C^0$, or even piecewise continuous classes, by standard limiting arguments.}
 on $X$ to be any function $f: X \to \C$ whose (inhomogeneous) Lipschitz norm
$$ \Vert f \Vert_{\Lip} := \sup_{x \in X} |f(x)| + \sup_{x,y \in X: x \neq y} \frac{|f(x)-f(y)|}{d(x,y)}$$
is finite.

\begin{example}[$\mu$ is strongly orthogonal to $1$-step nilsequences]\label{ex4}
Suppose that $G$ is a connected, simply-connected abelian Lie group (written multiplicatively) with a smooth metric $d$, and that $\Gamma$ is a closed subgroup of $G$ which is cocompact. Then $G/\Gamma$ is called a \emph{$1$-step nilmanifold}; it is a torus. Let $F : G/\Gamma \rightarrow \C$ be a Lipschitz function, and let $T_g : G/\Gamma \rightarrow G/\Gamma$ denote the action of $g$ on $G/\Gamma$. Then we have the estimate
\begin{equation}\label{mu-1-step}
\E_{n \in [N]} \mu(n) \overline{F(T_g^n x)} \ll_{A,G/\Gamma} \Vert F \Vert_{\Lip} \log^{-A} N
\end{equation}
for all $N > 1$, uniformly in $g \in G$ and $x \in G/\Gamma$.
\hfill 
\end{example}

The sequence $n \mapsto F(T_g^n x)$ is called a $1$-\textit{step nilsequence}.  If we specialize to the \emph{circle nilflow} case
$$ G := \begin{pmatrix}
1 & \R \\
0 & 1
\end{pmatrix} := \left\{ 
\begin{pmatrix}
1 & x \\
0 & 1
\end{pmatrix}: x \in \R \right \}; \quad \Gamma := 
\begin{pmatrix}
1 & \Z \\
0 & 1
\end{pmatrix} :=
\left\{ 
\begin{pmatrix}
1 & n \\
0 & 1
\end{pmatrix}: n \in \Z \right \}$$
then $G/\Gamma$ is isomorphic to the unit circle $\R/\Z$, and if we identify a real number $\alpha$ with the group element $\left(\begin{smallmatrix}
1 & \alpha \\
0 & 1
\end{smallmatrix}\right)$, then $T_\alpha: \R/\Z \to \R/\Z$ is just the shift $x \mapsto x+\alpha \md{1}$.  Using the standard character $e: \R/\Z \to \C$ as the Lipschitz function $F$, one then sees that \eqref{mu-alpha} is a special case of \eqref{mu-1-step}.  In fact, the two examples are more-or-less equivalent, as we shall see in \S \ref{1-nil-sec} where \eqref{mu-1-step} will be established.

The main aim of this paper is to generalise \eqref{mu-1-step} to cover $2$-step nilsequences.  In the companion paper \cite{green-tao-linearprimes} to this paper, we shall show how such estimates can be used to prove various ``complexity $2$'' estimates for the M\"obius and von Mangoldt functions.

Before stating our main result, we give the definition of $s$-step nilsequences in general, followed by some examples.

\begin{definition}[Nilmanifolds and nilsequences]\label{nil-def}
Let $G$ be a connected, simply connected, Lie group. We define the \emph{central series} $G_0 \supseteq G_1 \supseteq G_2 \supseteq \dots$ by defining $G_0 = G_1 = G$, and $G_{i+1} = [G, G_i]$ for $i \geq 2$, where the commutator group $[G,G_i]$ is the group generated by $\{ g h g^{-1}h^{-1}: g \in G, h \in G_i \}$. We say that $G$ is \emph{$s$-step nilpotent} if $G_{s+1} = {1}$. Let $\Gamma \subseteq G$ be a discrete, cocompact subgroup. Then the quotient $G/\Gamma$ is called an $s$-step nilmanifold. If $g \in G$ then $g$ acts on $G/\Gamma$ by left multiplication, $x \mapsto gx$. By a \emph{\textup{(}basic\textup{)} $s$-step nilsequence}, we mean a sequence of the form $(F(T_g^n \cdot  x))_{n \in \N}$, where $x \in G/\Gamma$ is a point, $F : G/\Gamma \rightarrow \C$ is a continuous function and $T_g : G/\Gamma \rightarrow G/\Gamma$ is left multiplication by $g$.  We say that the nilsequence is \emph{bounded} if $|F|$ takes values in $[-1,1]$. We may (arbitrarily) endow $G/\Gamma$ with a smooth Riemannian metric $d_{G/\Gamma}$. If the function $F$ is Lipschitz with respect to this metric, we shall refer to the nilsequence $(F(T_g^n \cdot x))_{n \in \N}$ as Lipschitz.
\end{definition}

\begin{remark}
In this paper we will usually suppress explicit mention of the metric $d_{G/\Gamma}$. Whenever an estimate is said to depend on a nilmanifold $G/\Gamma$, it should be assumed that it also depends on the choice of metric. See \cite{green-tao-linearprimes} for a more detailed discussion.
\end{remark}

Clearly every $1$-step nilsequence is a $2$-step nilsequence.  The next simplest example of nilsequences are quadratic phases.

\begin{example}[The Heisenberg nilflow, I]\label{nil2} Consider the example\footnote{For more detail on the Heisenberg nilflow, Appendix \ref{appendixA} may be consulted.  One can also generate quadratic phase sequences such as $e(n^2 \theta)$ using the slightly simpler skew shift nilflow (see e.g. \cite[Example 12.3]{green-tao-inverseu3}), but we shall refrain from doing so here as the underlying Lie group is disconnected and thus does not quite fall within the framework of Definition \ref{nil-def}.
}
\[ G := \left(\begin{smallmatrix}
1 & \R & \R\\
0 & 1  & \R\\
0 & 0  & 1
\end{smallmatrix}\right); \quad
\Gamma := \left(\begin{smallmatrix}
1 & \Z & \Z\\
0 & 1  & \Z\\
0 & 0  & 1
\end{smallmatrix}\right). \]
Then $G/\Gamma$ is a 2-step nilmanifold. Apart from a set of zero measure, $G/\Gamma$ may be identified with the fundamental domain
\[ \mathcal{F} := \left\{ \left(\begin{smallmatrix} 1 & x & y\\
0 & 1  & z\\
0 & 0  & 1
\end{smallmatrix}\right) : -1/2 < x,y,z \leq 1/2 \right\}\]
using the easily-verified fact that
\[ \left(\begin{smallmatrix} 1 & x & y\\
0 & 1  & z\\
0 & 0  & 1
\end{smallmatrix}\right) \equiv \left(\begin{smallmatrix} 1 & \{x\} & \{y - x[z]\}\\
0 & 1  & \{z\}\\
0 & 0  & 1
\end{smallmatrix}\right) \md{\Gamma}.\]
Here, $\{x\}$ refers to the fractional part of $x$ lying in the interval $(-1/2,1/2]$ and $[x] := x - \{x\}$. Writing
\[ g := \left(\begin{smallmatrix} 1 & -\theta & -\theta\\
0 & 1  & 2\\
0 & 0  & 1
\end{smallmatrix}\right),\] where $\theta \in \R$, one may check that 
\[ g^n \equiv \left(\begin{smallmatrix} 1 & \{-n\theta\} & \{n^2 \theta\}\\
0 & 1  & 0\\
0 & 0  & 1
\end{smallmatrix}\right) \md{\Gamma}.\]
Thus we see how functions with ``quadratic'' behaviour arise from 2-step nilsequences. The rather natural function $e(n^2\theta)$ does not quite arise as a Lipschitz nilsequence on the $3 \times 3$ Heisenberg group, since the function 
\[ \left(\begin{smallmatrix} 1 & x & y\\
0 & 1  & z\\
0 & 0  & 1
\end{smallmatrix}\right) \mapsto e(y)\]
on $\mathcal{F}$ does not extend to a continuous function on $G/\Gamma$.
The situation may be remedied by splitting $e(n^2 \theta)$ as the sum of (say) 10 functions $\chi(\{n\theta\})e(n^2\theta)$ where $\chi$ is a Lipschitz cutoff supported on an interval of width $1/5$. Each of the 100 functions 
\[ \left(\begin{smallmatrix} 1 & x & y\\
0 & 1  & z\\
0 & 0  & 1
\end{smallmatrix}\right) \mapsto \chi(x)\chi'(z)e(y)\]
\emph{does} extend to a Lipschitz function on $G/\Gamma$. By taking products one may realise $e(n^2 \theta)$ as a Lipschitz nilsequence on the $2$-step nilmanifold $(G/\Gamma)^{100}$.
\end{example}

In view of the previous example and our general intent in this paper, it is natural to ask for the estimate
\begin{equation}\label{mu-alphabeta} \E_{n \in [N]} \mu(n) e(-\alpha n^2 - \beta n - \gamma) \ll_A \log^{-A} N,\end{equation} with an implied constant independent of $\alpha,\beta$ and $\gamma$. We will prove such an estimate in \S \ref{quadsec}. Like \eqref{mu-alpha}, this bound is a fairly standard application of Vinogradov's version of the Hardy-Littlewood method, though somewhat more complicated due to the need to estimate quadratic exponential sums rather than just linear exponential sums. The proof of it has much in common with techniques pioneered by Hua \cite{hua-1} and Vinogradov \cite{vinogradov-1} in connection with the Goldbach-Waring problem. It should be thought of as a warm up for the main business of the paper.

As we have already mentioned, in \S \ref{1-nil-sec} we shall see that orthogonality to linear phases is more-or-less equivalent to orthogonality to $1$-step nilsequences. However, orthogonality to quadratic phases is significantly weaker than orthogonality to $2$-step nilsequences. This is because there are examples of $2$-step nilsequences which do not look much like quadratic phases.  

\begin{example}[The Heisenberg flow, II]\label{nil3} We repeat the analysis of the previous example, but with a less restrictive choice of $g$.
Take 
\[ g := \left(\begin{smallmatrix}
1 & \alpha  & \beta \\
0 & 1  & \gamma \\
0 & 0  & 1
\end{smallmatrix}\right). \]
A simple induction confirms that 
\[ g^n \cdot \left(\begin{smallmatrix}
1 & x  & y\\
 0 & 1 & z\\
0 & 0 & 1 \end{smallmatrix}\right)    = \left(\begin{smallmatrix}
1 & x + n\alpha  & y + n\beta + \frac{1}{2} n(n+1)\alpha \\
0 & 1  & z + n\gamma \\
0 & 0  & 1
\end{smallmatrix}\right)\]
When reduced to lie in the fundamental domain $\mathcal{F}$, one can end up with functions taking the form $[n\alpha]n\gamma$ (and related forms). These functions are known as \emph{generalised quadratics}, and they capture the spirit of 2-step nilsequences much more completely than genuine quadratic functions do. By repeating the tricks mentioned in the previous example one may actually approximate $e(-[n\sqrt{2}]n\sqrt{3})$ (say) outside of sets of arbitrarily small density as a Lipschitz nilsequence on some product of several copies of the Heisenberg example.
\end{example}

The previous two examples give some idea of what a 2-step nilsequence looks like. Our main result in this paper is that the M\"obius function is strongly asymptotically orthogonal to all such functions. This estimate is the case $s = 2$ of the M\"obius and Nilsequences Conjecture $\mbox{MN}(s)$: see \cite[\S 6]{green-tao-linearprimes} for further discussion.

\begin{main-theorem-1}[\mbox{MN}(2) conjecture] Suppose that $G/\Gamma$ is a $2$-step nilmanifold, and that $F : G/\Gamma \rightarrow \C$ is a Lipschitz function. Then for every $A > 0$ we have the estimate
\begin{equation}\label{mu-2-step}
\E_{n \in [N]} \mu(n) \overline{F(T_g^n x)} \ll_{A,G/\Gamma} \Vert F \Vert_{\Lip} \log^{-A} N
\end{equation}
uniformly in $g \in G$ and $x \in G/\Gamma$.
\end{main-theorem-1}

\begin{remark}
We conjecture that $\mbox{MN}(s)$ holds for arbitrary $s$, that is to say there is an analogue of the Main Theorem for $s$-step nilmanifolds for any $s \geq 1$.  The fact that the bound \eqref{mu-2-step} is uniform in $x$ is unsurprising (since $G/\Gamma$ is compact), as is the uniformity among all $F$ with fixed Lipschitz norm (thanks to the Arzel\`a-Ascoli theorem). The uniformity in $g$ is less trivial, and is quite important for applications.
\end{remark}

We shall prove the Main Theorem as a consequence of a similar result, Theorem \ref{mainthm-tech} below, in which the notion of a $2$-step nilsequence is replaced by a more technical type of sequence (a $1$-step nilsequence twisted by a locally quadratic phase) 
that is more tractable for analysis.  The proof of
Theorem \ref{mainthm-tech} is by far the most difficult portion of the paper and will occupy \S \ref{periodic-sec}-- \S\ref{major-sec}.  In comparison, the deduction of the Main Theorem from Theorem \ref{mainthm-tech} is more standard and is performed in \S \ref{technical-sec} and Appendix \ref{appendixA}.

The estimate \eqref{mu-alphabeta}, as well as estimates for generalised quadratic phases such as 
\[ \E_{n \in [N]} \mu(n) e(- [n\sqrt{2}] n\sqrt{3} ) = o(1),\]
are consequences of our main theorem. 

\begin{remark} The main result of this paper can then be combined with the \emph{Gowers Inverse Theorem} from \cite{green-tao-inverseu3} to obtain a number of new correlation estimates for the M\"obius function, such as
$$ \E_{x,d \in [N]} \mu(x) \mu(x+d) \mu(x+2d) \mu(x+3d) = o_{N \to \infty}(1)$$
and
\begin{align*} &\E_{x,h_1,h_2,h_3 \in [N]} \mu(x) \mu(x+h_1) \mu(x+h_2) \mu(x+h_3) \\
&\quad \mu(x+h_1+h_2) \mu(x+h_1+h_3) \mu(x+h_2+h_3) \mu(x+h_1+h_2+h_3) = o_{N \to \infty}(1)
\end{align*}
(compare with \eqref{mu1}, \eqref{mu2}).
It can also be used (with some additional effort) to establish an asymptotic for expressions such as
$$ \E_{x,d \in [N]} \Lambda(x) \Lambda(x+d) \Lambda(x+2d) \Lambda(x+3d)$$
as $N \to \infty$, thus enabling one to count the quadruples of number of primes $p_1 < p_2 < p_3 < p_4 \leq N$ in arithmetic progression up to a fixed level $N$.  We defer all of these applications to the companion paper \cite{green-tao-linearprimes}.
\end{remark}

\section{A technical reduction}\label{technical-sec}

In this section we present a technical counterpart of the Main Theorem, namely Theorem \ref{mainthm-tech} below, in which the $2$-step nilsequence is replaced by a more analytically tractable object, namely a $1$-step nilsequence twisted by a locally quadratic phase.  We then 
discuss how this result implies the Main Theorem.  The proof of Theorem \ref{mainthm-tech} will then occupy the rest of the paper (except for the Appendices).  We first need some notation.

\begin{definition}[Locally polynomial phases]\label{polynomial-def}  Let $S \subset \Z$ be a set of integers, and let $d \geq 0$.  A phase function $\phi : S \rightarrow \R/\Z$ is said to be \emph{locally degree $d$} on $S$ if whenever $n,h_1,\ldots,h_{d+1}$ are such that the $2^{d+1}$ quantities $n + \epsilon_1 h_1 + \ldots + \epsilon_{d+1} h_{d+1}$, $\epsilon_i \in \{0,1\}$ lie in the set $S$, we have
\begin{equation}\label{quadratic-def}
 \sum_{\epsilon \in \{0,1\}^{d+1}} (-1)^{\epsilon_1 + \ldots + \epsilon_{d+1}} \phi(n + \epsilon_1 h_1 + \ldots + \epsilon_{d+1} h_{d+1}) = 0.
 \end{equation}
We refer to phases of local degree $1$ as \emph{locally linear}, phases of local degree $2$ as \emph{locally quadratic}, and so forth.
\end{definition}

\begin{examples}  Constant phases have local degree $0$, while linear phases $\phi(n) := \alpha n$ for $\alpha \in \R$ have local degree $1$.
If $\alpha, \beta, \gamma$ are real numbers, then the phase $\phi(n) := \alpha n^2 + \beta n + \gamma \md{1}$ is globally quadratic (i.e. quadratic on all of $\Z$).  The phase $\phi(n) := \{ \alpha n \} \{ \beta n \} \gamma \md{1}$ is not globally quadratic, but it is locally quadratic on the \emph{Bohr set} $S := \{ n \in \Z: |\{ \alpha n\}|, |\{ \beta n \}| \leq 0.1 \}$, which is a set of positive density in $\Z$.  The phase $\phi(n) := \{ \alpha n\} \gamma \md{1}$ is locally linear on the same set. 
\hfill
\end{examples}

\begin{theorem}[$\mu$ is strongly orthogonal to local quadratics]\label{mainthm-tech} Let $G/\Gamma$ be a $1$-step nilmanifold, let
$F: G/\Gamma \to \C$ be a Lipschitz function, and let $g \in G$ and $x \in G/\Gamma$
be arbitrary.  Let $\phi : B_N \rightarrow \R/\Z$ be a phase which is locally quadratic on the \emph{Bohr set}\footnote{This definition of a Bohr set is not quite identical to other Bohr sets in the literature, for instance in \cite{green-tao-inverseu3}, but it is very closely related; see the proof of Lemma \ref{lem12.4}.} 
$B_N := \{ n \in [N] : F(T_g^n x) \neq 0\}$. Then we have
\[ \E_{n \in [N]} \mu(n) \overline{F(T_g^n x)} e(-\phi(n)) \ll_{G/\Gamma, A} \Vert F \Vert_{\Lip} \log^{-A} N.\]
\end{theorem}

The proof of Theorem \ref{mainthm-tech} is rather lengthy.
Let us assume it for now and deduce the Main Theorem.  The main proposition in achieving this deduction is

\begin{proposition}[$2$-step nilsequences as averages of twisted $1$-step nilsequences]\label{nil-to-local} Let\\
$G/\Gamma$ be a $2$-step nilmanifold and let $0 < \eps < 1/2$.
Let $F: G/\Gamma \to \C$ be a Lipschitz function with $\Vert F \Vert_{\Lip} \leq 1$, and let $g \in G$ and $x \in G/\Gamma$ be arbitrary.  Then there exists a $1$-step nilmanifold $\widetilde G/\widetilde \Gamma$ depending only on $G/\Gamma$ and a decomposition
\begin{equation}\label{fgn}
F(T_g^n x) = \E_{i \in I} w_i F_i(T_{g_i}^n x_i) e(-\phi_i(n)) + O(\eps)
\end{equation}
where
\begin{itemize}
 \item $I$ is a finite index set;
\item For each $i \in I$ the $w_i$ are complex numbers with $\E_{i \in I} |w_i| \ll \eps^{-O_{G/\Gamma}(1)}$;
\item $F_i: \tilde G/\tilde \Gamma \to \C$ is Lipschitz with norm $O_{G/\Gamma}(1)$;
\item $g_i \in \tilde G$;
\item $x_i \in \tilde G/\tilde \Gamma$;
\item 
$\phi_i: B_i \to \R/\Z$ is a phase function which is locally quadratic on the \emph{generalized Bohr set} 
$B_i := \{ n \in [N] : F_i(T_{g_i}^n x_i) \neq 0\}$.
\end{itemize}
\end{proposition}

We have a proof of a generalisation of this proposition to $k$-step nilsequences (they are averages of twisted $(k-1)$-step nilsequences). This proceeds using 
some rather algebraic considerations involving ``Hall-Petresco parallelepiped groups'' associated to the nilmanifold $G/\Gamma$. These considerations are very similar to, but more complicated than, the material in \cite[Appendix E]{green-tao-linearprimes}. We anticipate presenting the proof of this result in a future paper concerned with the generalisation of the Main Theorem to nilmanifolds of arbitrary step. 

In this paper we present a more computational approach involving so called Mal'cev bases \cite{corwin-greenleaf,malcev}. This approach is completely explicit when the group $G$ is a product of Heisenberg groups $\left(\begin{smallmatrix} 1 & \R & \R \\ 0 & 1 & \R \\ 0 & 0 & 1 \end{smallmatrix} \right)$. The reader will find remarks in \cite{green-tao-linearprimes} explaining that, in the theory of linear systems of complexity 2 (such as four-term APs) only examples of this type need be considered. 

The use of bases may seem overly explicit to some, but it should be noted that Mal'cev bases are in fact required to prove certain foundational topological properties of nilmanifolds. Those results are needed for the approach, just alluded to, that is taken in \cite[Appendix E]{green-tao-linearprimes}.

The proof of Proposition \ref{nil-to-local} may be found in Appendix \ref{appendixA}.  
Assuming it and Theorem \ref{mainthm-tech}, we can now derive the Main Theorem as follows.

\begin{proof}[Proof of the Main Theorem assuming Theorem \ref{mainthm-tech} and Proposition \ref{nil-to-local}]
Let $G/\Gamma$, $F$, $A$ be as in the Main Theorem. By renormalising we may assume that $\Vert f \Vert_{\Lip} \leq 1$. We apply Proposition \ref{nil-to-local} with $\eps := \log^{-A} N$ and obtain
a decomposition \eqref{fgn}.  Taking inner products with $\mu$, we obtain
$$ \E_{n \in [N]} \mu(n) \overline{F(T_g^n x)}
\ll \E_{i \in I} |w_i| \E_{n \in [N]} \mu(n) \overline{F_i(T_{g_i}^n x_i)} e(-\phi_i(n)) + \log^{-A} N.$$
Applying Theorem \ref{mainthm-tech}, we conclude that
$$ \E_{n \in [N]} \mu(n) \overline{F(T_g^n x)}
\ll_{A',\tilde G/\tilde \Gamma} \E_{i \in I} |w_i| \log^{-A'} N + \log^{-A} N$$
for any $A'$.  But $\E_{i \in I} |w_i| \ll (\log^A N)^{O_{G/\Gamma}(1)}$, so the claim follows by taking $A'$ suitably large.
\end{proof}

\begin{remark} Conversely it is also possible to deduce Theorem \ref{mainthm-tech} from the Main Theorem by obtaining a suitable converse to Proposition \ref{nil-to-local} (cf. the proof of \cite[Theorem 12.8]{green-tao-inverseu3}), but we will not do so here.
\end{remark}

\section{Orthogonality to periodic functions}\label{periodic-sec}

We now begin the proof of Theorem \ref{mainthm-tech}, which is the heart of this paper.  (The other major component of the paper is 
the proof of Proposition \ref{nil-to-local} in Appendix \ref{appendixA}. This can mostly be read independently of the part of the paper concerned with Theorem \ref{mainthm-tech}, though it will utilize the harmonic analysis tools collected in Appendix \ref{tools-sec}.)

Our strategy in proving Theorem \ref{mainthm-tech} shall be to establish the strong asymptotic orthogonality of the M\"obius function to increasingly large classes of sequences, starting with very simple ones and then moving on to ``higher degree'' sequences.  
 Let us begin with some generalities on how one can go about proving that $\mu$ is orthogonal to some function $F$. There are essentially two complementary methods for doing this. The first, which will feature prominently in this section, is appropriate when $F$ is multiplicative, for example $F = 1$ or $F = \chi$, where $\chi$ is some Dirichlet character to the modulus $q$. Then one may relate $\E_{n \in [N]} \mu(n) F(n)$ via Perron's Formula to zeros of $L$-functions such as $\zeta(s)$ and $L(s,\chi)$ in the critical strip, the orthogonality coming from the non-existence of zeros close to $\Re s = 1$. Siegel's theorem, concerning a possible zero near $s = 1$ when $\chi$ is real, is of particular importance. It implies the bound \eqref{mu-chi}, which we recall now:

\begin{proposition}\label{mu-chi-prop}
For any $A > 0$ we have
\begin{equation}\label{mu-chi-new} \E_{n \in [N]} \mu(n) \overline{\chi(n)} \ll_A q^{1/2} \log^{-A} N\end{equation} for all Dirichlet characters $\chi$ to modulus $q$.
\end{proposition}

\begin{remark} For the proof, see \cite[Prop. 5.29]{iwaniec-kowalski}. As noted in \cite[p. 124]{iwaniec-kowalski} there are difficulties involved in applying the standard Perron's formula approach to $\E_{n \in [N]} \mu(n) \chi(n)$ directly, and it is rather easier to first obtain bounds on $\E_{n \in [N]} \Lambda(n) \chi(n)$.  
Note that the bound is only non-trivial when the period $q$ is very small (e.g. $O(\log^A N)$) compared to $N$.  If one assumed GRH then one could improve the logarithmic decay here to a polynomial decay, which would of course lead to improvements in the other bounds in this paper.
\end{remark}

 As we will see later in this section, the need to consider zeros of $L$-functions also appears when dealing with functions $F$ which are not quite multiplicative. For example, they must play a role in the case $F(n) = e(an/q)$, since any Dirichlet character to modulus $q$ is a linear combination of a few such functions $F$.

 At the other end of the spectrum one has functions $F$ which are far from multiplicative, such as $F(n) = e(n\sqrt{2})$. For these functions a completely different method, due originally to Vinogradov, may be brought to bear. The sum $\E_{n \in [N]} \mu(n) F(n)$ is decomposed into so-called \emph{Type I} and \emph{Type II} sums, which can be estimated without having to understand the oscillation of $\mu$. Provided $F$ is not close to being multiplicative, those sums can often be shown to be small by (effective) harmonic analysis methods. We will discuss this method, in a modern and very neat incarnation due to Vaughan, in \S \ref{vaughan-sec}.

We now begin the proof of Theorem \ref{mainthm-tech}
by establishing the asymptotic orthogonality of the M\"obius function to periodic sequences, which can 
be viewed in some sense as ``$0$-step nilsequences'', and which will be needed to handle the ``major arc'' case when moving on to linear phases.  
More precisely, we show

\begin{proposition}[M\"obius is orthogonal to periodic sequences]\label{mob-period}
Let $f: \N \to \C$ be a sequence bounded in magnitude by $1$ which is periodic of some period $q \geq 1$. 
Then we have
$$ \E_{n \in [N]} \mu(n) \overline{f(n)} \ll_A q \log^{-A} N$$
for all $A > 0$, where the implied constant is ineffective.
\end{proposition}

\begin{proof}
We first establish the estimate under the additional assumption that $f(n)$ vanishes whenever $(n,q) \neq 1$.  Then $f$
can be viewed as a function on the multiplicative group $(\Z/q\Z)^\times$, and thus has a Fourier expansion
$$ f(n) = \sum_\chi \hat f(\chi) \chi(n), \hbox{ where } \hat f(\chi) := \E_{n \in (\Z/q\Z)^\times} f(n) \overline{\chi(n)},$$
with $\chi$ ranging over all the characters on $(\Z/q\Z)^\times$.  Applying Proposition \ref{mu-chi-prop} and the triangle inequality, we conclude
$$ \E_{n \in [N]} \mu(n) \overline{f(n)} \ll_A q^{1/2} \log^{-A} N \big(\sum_\chi |\hat f(\chi)|\big).$$
But from Cauchy-Schwarz and Plancherel we have
$$ \sum_\chi |\hat f(\chi)| \leq \phi(q)^{1/2} (\sum_\chi |\hat f(\chi)|^2)^{1/2}
= \phi(q)^{1/2} (\E_{n \in (\Z/q\Z)^\times} |f(n)|^2)^{1/2} = O( \phi(q)^{1/2} ),$$
where $\phi(q) := |(\Z/q\Z)^\times|$ is the Euler totient function.  Since $\phi(q) \leq q$, the claim follows.

 Now we consider the general case, in which $(n,q)$ is not necessarily equal to $1$ on the support of $f$.
Observe that if $\mu(n)$ is non-zero, then $n$ is square-free, and we can split $n = dm$, where $d = (n,q)$ is
square-free (so $\mu^2(d)=1$) and $m$ is coprime to $q$.  Furthermore we have $\mu(n) = \mu(d) \mu(m)$.  We thus obtain the decomposition
\begin{equation}\label{2-star} \E_{n \in [N]} \mu(n) \overline{f(n)} = \frac{1}{N} \sum_{d|q; \mu^2(d) = 1} \mu(d) \sum_{1 \leq m \leq N/d} \mu(m) \overline{f(dm)} 1_{(m,q)=1}.\end{equation}
The sequence $m \mapsto f(dm) 1_{(m,q)=1}$ is periodic of period $q/d$ and vanishes whenever $(m,q/d) \neq 1$, hence by the preceding arguments
$$ \sum_{1 \leq m \leq N/d} \mu(m) \overline{f(dm)} 1_{(m,q)=1} \ll_A \frac{Nq}{d^2}  \log^{-A} N.$$
Thus from \eqref{2-star} we have
\[ \E_{n \in [N]} \mu(n) \overline{f(n)} \ll_A q\log^{-A} N \sum_{d|q} \frac{1}{d^2}\ll q \log^{-A} N,\]
concluding the proof of Proposition \ref{mob-period}.
\end{proof}

\section{Vaughan's identity}\label{vaughan-sec}

 In this section we discuss Vinogradov's method for proving that the M\"obius function $\mu$ is orthogonal to a function $F : \mathbb{N} \rightarrow \C$. As we remarked in \S \ref{periodic-sec}, this involves a decomposition of $\E_{n \in [N]} \mu(n) F(n)$ into Type I and Type II sums. The first argument of this type was due to Vinogradov (who worked with the von Mangoldt function $\Lambda$ instead of $\mu$). We will use a particularly simple identity due to Vaughan \cite{vaughan} to effect our decomposition into Type I and II sums. See \cite[Chapter 13]{iwaniec-kowalski} for a nice discussion of this and related identities.

 Let us begin with a few preliminary remarks on our strategy for dealing with Vinogradov's method throughout the paper. The normal method for proving Davenport's estimate \eqref{mu-alpha} would be to divide all $\alpha \in \R/\Z$ into two classes: the \textit{major arcs}, where $\alpha \approx a/q$ for some reasonably small $q$, and the \textit{minor arcs}, consisting of all other $\alpha$. If $\alpha$ lies in a major arc then one would use Proposition \ref{mob-period} to estimate $\E_{n \in [N]} \mu(n) e(\alpha n)$. If, by contrast, $\alpha$ lies in a minor arc then one would establish that Type I and II sums involving $f(n) = e(\alpha n)$ are small (see below for an explanation of what this means). Vaughan \cite[Chapter 3]{vaughan-hlm} may be consulted for details.

 We will adopt what we call an ``inverse'' strategy. In \S \ref{linear-sec} we will provide a proof of Davenport's estimate. There we will \textit{assume} that either a Type I or a Type II sum involving $f(n) = e(\alpha n)$ is large, and then \textit{deduce} that $\alpha$ lies in a major arc. The distinction between our argument and the standard one may seem rather unimportant, and indeed the two proofs are logically equivalent. However when it comes to dealing with more complicated functions $f$, such as locally quadratic phases which arise from the consideration of $2$-step nilsequences, the inverse strategy is very helpful. There it is much less obvious what one should mean by a ``major arc'', and even once the definition is made it is not obvious how to handle it in the context of Type I and II sums.

 In light of Lemma \ref{telescope}, it suffices to establish decay estimates for $\E_{N < n \leq 2N} \mu(n) f(n)$. The next lemma gives Vaughan's decomposition of sums of this kind.

\begin{lemma}[Vaughan's identity]\label{Vaughan}  Let $U, V, N$ be positive integers with $U V \leq N$, 
and $f: \N \to \C$ be a sequence.  
Then we have
\begin{equation}\label{vaughn-id} \E_{N < n \leq 2N} \mu(n) \overline{f(n)} = -\typeI + \typeII
\end{equation}
where $\typeI$ is the \emph{Type I} expression
\begin{equation}\label{TI-def}
 \typeI := \frac{1}{N}\sum_{1 \leq d \leq UV} a_d \sum_{N/d < w \leq 2N/d} \overline{f(dw)} 
 \end{equation}
 in which 
 \[ a_d := \sum_{bc=d: b \leq U, c \leq V} \mu(b) \mu(c), \]
and $\typeII$ is the \emph{Type II} expression
\begin{equation}\label{TII-def}
\typeII := \frac{1}{N} \sum_{V < d \leq 2N/U} \sum_{\max(U,N/d) < w \leq 2N/d} \mu(w)b_d  \overline{f(dw)}
\end{equation}
in which
\[ b_d := \sum_{c|d: c > V} \mu(c).\]
\end{lemma}

\begin{remark} One of the key points in the analysis of Type I sums is that the precise form of the coefficients $a_d$ is almost completely irrelevant: we will apply the Cauchy-Schwarz inequality, and so only the mean square size of these coefficients will concern us. The same is true in the analysis of Type II sums. In this case it is the coefficients $\mu(w)$ and $b_d$ which get removed by the Cauchy-Schwarz inequality.
\end{remark}

 There is considerable flexibility in the choice of the parameters $U$ and $V$. We will take $U = V = N^{1/3}$ in our applications.

\begin{proof} We follow \cite[\S 13.4 -- 5]{iwaniec-kowalski}. Observe that for any positive integer $n$ we have
\[ \mu(n) = \sum_{b,c : bc | n} \mu(b)\mu(c).\]
We split the range of the sum over $b,c$ into four ranges: (i) $b \leq U$, $c \leq V$; (ii) $b > U$, $c \leq V$; (iii) $b \leq U$, $c > V$ and (iv) $b > U$, $c > V$. Denoting the associated sums $\Sigma_1,\dots,\Sigma_4$, it is easy to check that $\Sigma_2 = \Sigma_3 = -\Sigma_1$. It follows that
\[ \mu(n) = -\Sigma_1 + \Sigma_4 = -\sum_{\substack{b \leq U \\ c \leq V \\ bc | n}} \mu(b)\mu(c) + \sum_{\substack{b > U \\ c > V \\ bc | n}} \mu(b)\mu(c).\]
Multiplying by $\overline{f(n)}$ and summing over $N < n \leq 2N$, we have Vaughan's identity:
\begin{eqnarray*} \E_{N < n \leq 2N} \mu(n) \overline{f(n)} & = & -\E_{N < n \leq 2N}\sum_{\substack{b \leq U \\ c \leq V \\ bc | n}} \mu(b)\mu(c)\overline{f(n)} + \E_{N < n \leq 2N}\sum_{\substack{b > U \\ c > V \\ bc | n}} \mu(b)\mu(c)\overline{f(n)} \\ & := &  -\typeI + \typeII.\end{eqnarray*}
It is an easy matter to confirm that $\typeI$ may be written in the form \eqref{TI-def}, after making the substitution $d = bc$ and $n = dw$. One may also check that $\typeII$ may be written in the form \eqref{TII-def} after making the substitution $w = b$ and $n = dw$.\end{proof}

 Vaughan's identity tells us that if $\E_{N < n \leq 2N} \mu(n) \overline{f(n)}$ is large then either $\typeI$ or $\typeII$ is large. The next proposition shows how this information is processed, by using the Cauchy-Schwarz inequality to eliminate the parameters $a_d$, $b_w$ and $\mu(w)$, leaving behind estimates which only involve the explicit function $f$.

\begin{proposition}[Inverse theorem for $\E_{N < n \leq 2N} \mu(n) \overline{f(n)}$]\label{inverse-prop}  
Let $U, V, N$ be positive integers with $UV \leq N$, and let
$f: \N \to \C$ be a function with $\Vert f \Vert_{\infty} = O(1)$ such that
\[ |\E_{N < n \leq 2N} \mu(n) \overline{f(n)}| \geq \delta\]
for some $\delta > 0$.  Then one of the following statements holds:
\begin{itemize}
\item \textup{(Type I sum is large)} There exists an integer $1 \leq D \leq UV$ such that
\begin{equation}\label{typeI}
|\E_{N/d < w \leq 2N/d} f(dw)| \gg \delta \log^{-5/2} N
\end{equation}
for $\gg \delta^2 D \log^{-5} N$ integers $d$ such that $D < d \leq 2D$.\vspace{11pt}

\item \textup{(Type II sum is large)} There exist integers $D,W$ with $V/2 \leq D \leq 4N/U$ and $N/4 \leq DW \leq 4N$, such that
\begin{equation}\label{typeII}
 |\E_{D < d,d' \leq 2D} \E_{W < w,w' \leq 2W} f(dw) \overline{f(d'w)} \overline{f(dw')} f(d'w')| \gg \delta^4 \log^{-14} N.
 \end{equation}
\end{itemize}
\end{proposition}

\begin{remark} The estimate \eqref{typeI} may be viewed as an assertion that $f$ behaves periodically, while \eqref{typeII} is an assertion that $f$ behaves multiplicatively.  The numerical exponents could probably be improved slightly here, but we will not need such refinements here (especially since our bounds will eventually become ineffective anyway).
\end{remark}

\begin{proof}  We may of course take $N$ to be large.  Applying Lemma \ref{Vaughan}, we see that either $|\typeI| \geq \delta/2$
or $|\typeII| \geq \delta/2$.

 Suppose first that the Type I expression is large, that is to say $|\typeI| \geq \delta/2$ where $\typeI$ is given by \eqref{TI-def}. Using the crude bound $|a_d| \leq \tau(d)$, where 
$\tau(d) := \sum_{b|d} 1$ is the divisor function, we have
\[ \sum_{1 \leq d \leq UV} \frac{\tau(d)}{d} \left| \E_{N/d < w \leq 2N/d} f(dw) \right| \gg \delta.\]
By Cauchy-Schwarz inequality this implies that
\[ \sum_{1 \leq d \leq UV} \frac{1}{d} \left| \E_{N/d < w \leq 2N/d} f(dw) \right|^2 \gg \delta^2 \big( \sum_{1 \leq d \leq UV} \frac{\tau^2(d)}{d}\big)^{-1}.\] Invoking the divisor moment estimate \eqref{taud}, it follows that
\[ \sum_{1 \leq d \leq UV} \frac{1}{d} \left| \E_{N/d < w \leq 2N/d} f(dw) \right|^2 \gg \delta^2 \log^{-4} N.\]
Dividing the region $1 \leq d \leq UV$ into dyadic blocks $D < d \leq 2D$ (allowing for some slight overlap) and applying the pigeonhole
principle we obtain 
\[ \sum_{D < d \leq 2D} \left| \E_{N/d < w \leq 2N/d} f(dw) \right|^2 \gg \delta^2 D\log^{-5} N\]
for some $D$, $1 \leq D \leq UV$.  Since the summand is bounded by $O(1)$, a simple averaging argument confirms that $\left| \E_{N/d < w \leq 2N/d} f(dw) \right| \gg \delta \log^{-5/2} N$ for at least $\gg \delta^{-2} D\log^{-5} N$ values of $d$, which is what we wanted to prove.

 Now suppose instead that the Type II expression is large, that is $|\typeII| \geq \delta/2$. Using the evident bound $|b_d| \leq \tau(d)$, we conclude
\[ \sum_{V < d \leq 2N/U} \tau(d) \big|\sum_{N/d < w \leq 2N/d} 1_{w>U} \mu(w) f(dw)\big| \gg N \delta.\]
Applying Cauchy-Schwarz and the divisor moment estimate \eqref{taud} once again, we conclude that
\[ \sum_{V < d \leq 2N/U} d \big|\sum_{N/d < w \leq 2N/d} 1_{w > U} \mu(w) f(dw)\big|^2 \gg N^2 \delta^2 \log^{-4} N.\]
By dyadic decomposition, we thus can find integers $D,W$ with $V/2 \leq D \leq 4N/U$ and $N/4 \leq DW \leq 4N$ such that
\[ \sum_{D < d \leq 2D} \big|\sum_{W < w \leq 2W} 1_{I_d}(w) \mu(w) f(dw)\big|^2 \gg \frac{N^2}{D} \delta^2 \log^{-5} N,\]
where $I_d$ is the discrete interval $\{ w > U: N/d < w \leq 2N/d \}$.  
Applying Lemma \ref{complete} to remove the cutoff $1_{I_d}(w)$, we obtain
\[ \sum_{D < d \leq 2D} \big|\sum_{W < w \leq 2W} \mu(w) f(dw) e(\alpha w)\big|^2 \gg N \delta^2 \log^{-7} N.\]
for some $\alpha \in \R/\Z$. Expanding the left-hand side as
\[ \sum_{W < w,w' \leq 2W} \sum_{D < d \leq 2D} \b(w,w') f(dw) \overline{f(dw')},\]
where we use $\b()$ to denote a bounded function whose exact form we do not care about (see Appendix \ref{tools-sec}),
the required inequality \eqref{typeII} follows from the Cauchy-Schwarz inequality in the form of Lemma \ref{cz}.
\end{proof}

\section{Orthogonality to linear phase functions}\label{linear-sec}

 As a first application of Proposition \ref{inverse-prop}, let us recall the standard proof of Davenport's estimate \eqref{mu-alpha}.  We do this partly for expository
reasons, to illustrate the ``inverse'' approach to dealing with Type I and II sums, and also because we will need \eqref{mu-alpha} to treat the ``major arc'' case of quadratic phases in later sections.  As we shall see, the linear case is particularly
easy, because the exponential sums can be easily computed (using \eqref{exponential-sum}).  Here and in the rest of the paper
we will be using some standard tools from harmonic analysis, together with the notations $\|x\|_{\R/\Z}$ and $\|x\|_{\R/\Z,Q}$,
which we summarize in Appendix \ref{tools-sec}.

 We begin with a partial result, which is weaker than \eqref{mu-alpha} in that it only resolves the theorem for the ``minor arc'' values of $\alpha$, but
has the advantage of being completely effective, as it does not require any information on Siegel zeroes.

\begin{proposition}[Correlation with a linear phase implies major arc]\label{corr-phase}  Let $\alpha \in \R$, let $A > 0$, and let $N$ be a large integer such that
\begin{equation}\label{mu-alpha-fail}
 |\E_{N < n \leq 2N} \mu(n) e(-\alpha n)| \geq \log^{-A} N.
\end{equation}
Then there exists $D$, $1 \leq D \ll N^{2/3}$, such that
\begin{equation}\label{d-spin}
 \#\{ 1 \leq d \leq 2D: \| \alpha d \|_{\R/\Z} \ll \frac{D}{N} \log^{4A + 14} N \} \gg D \log^{-4A - 14} N.
\end{equation}
\end{proposition}

\begin{proof}  We apply Proposition \ref{inverse-prop} with $U=V=N^{1/3}$ and conclude one of the following statements:
\begin{itemize}
\item (Type I sum is large) There exists $D$, $1 \leq D \leq N^{2/3}$, such that
$$ |\E_{N/d < w \leq 2N/d} e(\alpha dw)| \gg \log^{-A-5/2} N$$
for $\gg D \log^{-2A-5} N$ values of $D < d \leq 2D$.\vspace{11pt}

\item (Type II sum is large) There exist integers $D,W$ with $N^{1/3} \ll D \ll N^{2/3}$ and $N/8 \leq DW \leq 8N$ such that
$$
 |\E_{D < d,d' \leq 2D} \E_{W < w,w' \leq 2W} e(\alpha dw - \alpha d' w - \alpha dw' + \alpha d'w')| \gg \log^{-4A-14} N.
$$
\end{itemize}
Suppose first that the Type I sum is large.  Applying \eqref{exponential-sum} we conclude
that there are $\gg D \log^{-2A-5} N$ values of $d$, $D < d \leq 2D$, for which
\[ \| \alpha d \|_{\R/\Z} \ll \frac{D}{N}\log^{-A - 5/2} N.\]
This implies \eqref{d-spin} with some room to spare.

 Now suppose instead that the Type II sum is large.  By the pigeonhole principle we can find $d', w'$ such that
$$  |\E_{D < d \leq 2D} \E_{W < w \leq 2W} e(\alpha dw - \alpha d' w - \alpha dw' + \alpha d'w')| \gg \log^{-4A-14} N$$
and hence by the triangle inequality
$$  \E_{D < d \leq 2D} |\E_{W < w \leq 2W} e(\alpha (d - d') w)| \gg \log^{-4A-14} N.$$
Applying \eqref{exponential-sum} we obtain
\[ \E_{D < d \leq 2D} \min\big(1, \frac{D}{N \Vert \alpha (d - d') \Vert_{\R/\Z}}\big) \gg \log^{-4A - 14} N,\] and thus after a simple averaging argument we establish
\[ \#\{ D < d \leq 2D: \| \alpha d - \alpha d' \|_{\R/\Z} \ll \frac{D}{N} \log^{4A+14} N \} \gg D \log^{-4A-14} N.\]
Substituting $\tilde d := d-d'$, we conclude
$$ \#\{ -2D \leq \tilde d \leq 2D: \| \alpha \tilde d \|_{\R/\Z} \ll \frac{D}{N} \log^{4A+14} N \} \gg D \log^{-4A-14} N.$$
Since $D \geq N^{1/3}$, we can easily remove the degenerate contribution when $\tilde d=0$.  The claim \eqref{d-spin} then follows by
symmetry.
\end{proof}

 The next task is to understand exactly what the condition \eqref{d-spin} implies.  It is clear that it is some sort of ``major arc'' condition,
as it forces $\alpha$ to lie close to a rational number with reasonably small denominator.  A na\"{\i}ve inspection of \eqref{d-spin} would lead one
to guess that this denominator is of size $D$ or so; however it turns out that one can reduce the size of the denominator substantially, to be
a power of $\log N$.  Indeed, we have

\begin{corollary}[Correlation with a linear phase implies major arc, II]\label{corr-phase-2}  Let $\alpha \in \R$, let $A > 0$ and let $N$ be a large integer such that \eqref{mu-alpha-fail} holds.  Then 
$$\|\alpha\|_{\R/\Z, 16\log^{8(A+4)} N} \ll \frac{\log^{28(A+4)} N}{N}.$$
The implied constant is effective.
\end{corollary}

\begin{proof} We apply Proposition \ref{corr-phase} to obtain $D$, $1 \leq D \leq N^{2/3}$, obeying \eqref{d-spin}.  If $D \leq \log^{8(A+4)}N$
then the claim follows directly from \eqref{d-spin}. If instead
$D \geq \log^{8(A+4)} N$, we may apply Lemma \ref{lem3.1}(ii) with $I = \{1,\ldots,2D\}$, $\delta_1 \ll \frac{D}{N} \log^{4(A+4)} N$, and
$\delta_2 \gg \log^{-4(A+4)} N$ to obtain the claim.
\end{proof}

 When $\alpha$ \textit{is} major arc, i.e. when $\Vert \alpha \Vert_{\R/\Z,Q}$ is small, we may proceed using Proposition \ref{mu-chi-prop}.

\begin{proposition}[Major arc phases are orthogonal to M\"obius]\label{major-orthog}  
Let $N$ be a large integer, let $\alpha$ be a real number, and let $Q, K \geq 1$
be such that $\|\alpha\|_{\R/\Z,Q} \leq K/N$.  Then we have
$$ |\E_{N < n \leq 2N} \mu(n) e(-\alpha n)| \ll_A Q^{1/2} K^{1/2} \log^{-A} N
$$
for any $A > 0$ \textup{(}the implied constant is ineffective\textup{)}.
\end{proposition}

\begin{proof}  Let $1 \leq M < N$ be a parameter to be chosen later.  Then by partitioning the interval $\{ N < n \leq 2N \}$ into
intervals of length $M$, plus a remainder, we conclude that
$$  |\E_{N < n \leq 2N} \mu(n) e(-\alpha n)| \leq \sup_{|I| = M; I \subset [N,2N]} |\frac{1}{M} \sum_{n \in I} \mu(n) e(\alpha n)| + O( \frac{M}{N} ).$$
By hypothesis, we have integers $a$ and $1 \leq q \leq Q$ such that $|\alpha - \frac{a}{q}| \leq \frac{K}{N}$.  We thus have
$$ e(-\alpha n) = e(-a n/q) e( -(\alpha - a/q) n) = e(-a n/q) e( -(\alpha - a/q) n_I) + O( \frac{K M}{N} )$$
for any $n, n_I \in I$.  Discarding the constant phase $e( -(\alpha - a/q) n_I)$, we conclude
$$  |\E_{N < n \leq 2N} \mu(n) e(-\alpha n)| \leq \sup_{|I| = M; I \subset [N,2N]} |\frac{1}{M} \sum_{n \in I} \mu(n) e(-a n / q)| + 
 O( \frac{K M}{N} ).$$
Applying Proposition \ref{mob-period} (replacing $A$ by $2A$) we have
$$ |\frac{1}{M} \sum_{n \in I} \mu(n) e(-a n / q)| \ll_A  \frac{qN}{M} \log^{-2A} N.$$
Combining these estimates and making the optimal choice $M = q^{1/2} K^{-1/2} N \log^{-A} N$, we obtain the claim.
\end{proof}

 Combining Corollary \ref{corr-phase-2} with Proposition \ref{major-orthog} (and selecting the parameters $A$ appropriately) we conclude
the unconditional estimate
$$ |\E_{N < n \leq 2N} \mu(n) e(-\alpha n)| \ll_A \log^{-A} N,$$
uniformly in $\alpha \in \R/\Z$ and for any $A > 0$. Here the implied constant is ineffective.  Davenport's estimate \eqref{mu-alpha} then follows from Lemma \ref{telescope} (with $\varphi \equiv 1$),
observing that the additional linear phase created by that lemma can be easily absorbed.
\endproof

\section{Orthogonality to linear objects}\label{1-nil-sec}

 Our aim in this section is to prove that the M\"obius function $\mu$ is orthogonal to various functions $f : \Z \rightarrow \C$ of ``linear'' type. We begin by proving \eqref{mu-1-step}, which asserts that $\mu$ is orthogonal to $1$-step nilsequences. Then, in Proposition \ref{ortho-prog},we confirm that $\mu$ is orthogonal to a certain type of locally linear phase function. This proposition is needed for our later analysis of $2$-step nilsequences (indeed, it essentially forms the ``major arc'' part of that analysis; see \S \ref{major-sec}).

\begin{proof}[Proof of \eqref{mu-1-step}] Let us begin by recalling what it is we are trying to prove. We have an abelian Lie group $G$ and a cocompact discrete subgroup $\Gamma \leq G$. Let $F : G/\Gamma \rightarrow \C$ be any Lipschitz function. 
Then we wish to show that
\begin{equation}
\E_{n \in [N]} \mu(n) \overline{F(g^n x)} \ll_{A,G/\Gamma} \Vert F \Vert_{\Lip} \log^{-A} N
\end{equation}
uniformly in $g \in G$ and $x \in G/\Gamma$. Now $G/\Gamma$ is isomorphic to the product of a torus and a finite abelian group, and hence to some subgroup of a torus $(\R/\Z)^d$. By Lemma \ref{lip-extend}, we may assume that $F$ is defined on all of this torus. Let $0 < \eps < 1$ be arbitrary. By renormalising, we may also assume that $\Vert F \Vert_{\Lip} = 1$. By Lemma \ref{fourier-lip}, we may write
\[ F(x) = \sum_{j = 1}^J c_j e(m_j \cdot x) + O_{d}(\eps^{1/2})\]
(say), where $c_j = O(1)$ and $J = O_d(\eps^{-d})$. Writing $g = (\alpha_1,\dots,\alpha_d)$, we have
\[ F(g^n x) = \sum_{j = 1}^J c_j e(m_j \cdot x) e\big(n(\alpha_1 m_j^{(1)} + \dots + \alpha_d m_j^{(d)})\big) + O_d(\eps^{1/2}).\]
Multiplying by $\mu$ and taking the expectation over $n \leq N$, the contribution of each of the $J$ terms here is $O_A( \log^{-A} N)$ for any $A > 0$, thanks to \eqref{mu-alpha}. We therefore have
\[ \E_{n \in [N]} \mu(n) \overline{F(g^n x)} \ll_{A,d} \eps^{-d} \log^{-A} N + \eps^{1/2}.\]
Optimising this in $\eps$ and recalling that $A > 0$ was arbitrary, we obtain the claim.
\end{proof}

Our other goal in this section is to establish, in Proposition \ref{ortho-prog}, orthogonality of $\mu$ to phase functions which are almost linear on Bohr sets. 

\begin{definition}[Bohr sets]\label{bohr-set} Let $N \geq 1$.  Let $G/\Gamma$ be a $1$-step nilmanifold (i.e. a compact abelian Lie group).  Then $G/\Gamma$ can be embedded as a closed subgroup of a finite-dimensional torus $(\R/\Z)^d$, and we let $d_{G/\Gamma}(x,y) := \| x y^{-1} \|_{G/\Gamma}$ be the metric on $G/\Gamma$ induced from such an embedding (chosen arbitrarily), where we give the torus the metric induced by the $l^\infty$ norm \eqref{torus-infty}.
For any $g \in G$ and any $n \in \Z$, we define the ``norm'' $\|n\|_g = \|n\|_{g,N}$ for all $n \in \Z$ by the formula 
\[ \|n\|_g := \| g^n \|_{G/\Gamma} + |\frac{n}{N}|,\]
and then define the \emph{Bohr sets} $B_g(n_0,\rho) \subset \Z$ for any $n_0 \in \Z$ and $\rho > 0$ as
\[ B_g(n_0,\rho) := \{ n \in \Z: \|n-n_0\|_g < \rho \}.\]
Thus we have $B_g(n_0,\rho) = n_0 + B_g(0,\rho)$.
\end{definition}

\begin{remarks} These Bohr sets are closely related to the sets $B_N$ appearing in Theorem \ref{mainthm-tech}, and also to more ``traditional''
Bohr sets in the literature; see the proof of Lemma \ref{lem12.4} below.  We observe the sub-homogeneity property $\|nm\|_g \leq |n| \|m\|_g$ for all $n,m \in \Z$, with equality
$\|nm\|_g = |n| \|m\|_g$ holding whenever $|n| \|m\|_g < c$ for some constant $c_{G/\Gamma} > 0$.
We shall use these facts frequently in the sequel without further
comment.  
\end{remarks}

Some other easy properties of Bohr sets are contained in the following lemma.

\begin{lemma}[Bohr set estimates]\label{bohr-size}  
Let $N \geq 1$, let $G/\Gamma$ be a $1$-step nilmanifold, and let $g \in G$.  Let $0 < \rho < 1/2$.
\begin{itemize}
\item[(a)] \textup{(Lower bound)} We have $|B_g(0,\rho)| \gg_{G/\Gamma} \rho^{-O_{G/\Gamma}(1)} N$.
\item[(b)] \textup{(Doubling property)} We have $|B_g(0,2\rho)| \ll_{G/\Gamma} |B_g(0,\rho)|$.
\item[(c)] \textup{(Divisibility)} For any integer $d \geq 1$ we have
$$ |\{ n \in B_g(0,\rho): d|n \}| \gg_{G/\Gamma} \frac{1}{d} |B_g(0,\rho)|.$$
\end{itemize}
\end{lemma}

\begin{proof}  To obtain (a), we cover $G/\Gamma$ by $O_{G/\Gamma}(\rho^{-O_{G/\Gamma}(1)})$ balls $B$ of radius $\rho/4$, and also cover
$\{1,\ldots,N\}$ into intervals $I$ of length $\rho N/4$.  By the pigeonhole principle we can find an interval $I$ and a ball $B$
such that $S := \{n :  n \in I: g^n \in B \}$ has cardinality $\gg_{G/\Gamma} \rho^{-O_{G/\Gamma}(1)} N$.  The claim then follows from the triangle inequality. Indeed if $n, n_0 \in S$ then $|(n - n_0)/N| \leq \rho/2$ and $\Vert g^{n - n_0} \Vert_{G/\Gamma} \leq \rho/2$, and thus $S - n_0 \subseteq B_g(0,\rho)$. It follows that $|B_g(0,\rho)| \geq |S|$.

 The proof of (b) is very similar. We cover the ball with centre $0$ and radius $2\rho$ in $G/\Gamma$ by $O_{G/\Gamma}(1)$ balls $B$ of radius $\rho/4$, and the interval $\{1,\dots, \rho N\}$ by $O(1)$ intervals $I$ of length $\rho N/4$. By the pigeonhole principle, there is an interval $I$ and a ball $B$ such that the set $S :=  \{n \in B_g(0,2\rho) : n \in I : g^n \in B\}$ has cardinality $\gg_{G/\Gamma} |B_g(0,2\rho)|$. Note, however, that if $n,n_0 \in S$ then$|(n - n_0)/N| \leq \rho/2$ and $\Vert g^{n - n_0} \Vert_{G/\Gamma} \leq \rho/2$, and so $S - n_0 \subseteq B_g(0,\rho)$. It follows that $|B_g(0,\rho)| \geq |S|$.

 Finally, we establish (c). By the pigeonhole principle there is some residue class $X_b := \{x \in \Z : x \equiv b \md{d}\}$ for which $|B_g(0,\rho/2) \cap X_b| \geq d^{-1}|B_g(0,\rho/2)|$. Note, however, that if $n,n_0 \in B_g(0,\rho/2) \cap X_b$ then $d | (n - n_0)$ and $n - n_0 \in B_g(0,\rho)$. The result now follows from (b).\end{proof}

 As we have remarked, the next result will form the ``major arc'' part of our analysis of $2$-step nilsequences. It may appear a little technical at this point, but has been designed to cover everything we need in the later application.

\begin{proposition}[Orthogonality to almost linear phases on Bohr sets]\label{ortho-prog}
Let $N \in \N$ be large, let $G/\Gamma$ be a $1$-step nilmanifold, let $g \in G$, let $\rho \in (0,1)$ and let $B_g(n_0,\rho)$ be some Bohr set contained in $\{N+1,\ldots,2N\}$.  Let $\psi: \Z \to \R^+$ be a non-negative function supported on $B_g(n_0,\rho)$ which obeys
the Lipschitz estimate
\begin{equation}\label{lipmn}
|\psi(n) - \psi(m)| \ll \| n-m\|_g
\end{equation}
for all $n,m \in \Z$.  Let $q \in [1, N/100]$ be an integer, let $\eps \in (0,1)$, and let $\phi: \Z \to \R/\Z$ be a phase obeying the approximate linearity estimate
\begin{equation}\label{linearity}
\| \phi(x+h_1+h_2) - \phi(x+h_1) - \phi(x+h_2) + \phi(x) \|_{\R/\Z} \ll \eps 
\end{equation}
whenever
$x, x+h_1, x+h_2, x+h_1+h_2 \in B_g(n_0,10\rho)$ and $q | h_1,h_2$.  Then for any $\kappa \in (0,\rho]$ we have
\[ |\E_{N < n \leq 2N} \mu(n) \psi(n) e( -\phi(n) )| \ll_{A,G/\Gamma} \kappa^{-O_{G/\Gamma}(1)} q^3 \log^{-A} N + (\eps + \kappa) \E_{N < n \leq 2N} |\psi(n)|\]
for all $A > 0$ \textup{(}the constant is ineffective\textup{)}.
\end{proposition}

\begin{proof}  
We can divide the interval $\{ N+1,\ldots,2N\}$ into $q$ residue classes $X_1,\ldots,X_q$ modulo $q$.  By the triangle inequality it suffices to show
that
\begin{align*} |\E_{N < n \leq 2N} \mu(n) & 1_{X_s}(n) \psi(n) e( -\phi(n) )| \\ & \ll_{A,G/\Gamma} \kappa^{-C} q^2 \log^{-A} N + (\eps + \kappa) \E_{N < n \leq 2N} |\psi(n)| 1_{X_s}(n)\end{align*}
for all $s$, $1 \leq s \leq q$.  

 Fix $s$.  Without loss of generality we may assume that $X_s \cap B_g(n_0,\rho)$ is non-empty, thus we may choose $n_s \in X_s \cap B_g(n_0,\rho)$. We work in the group $\Z/p\Z$ where $p \in [10N,20N]$ is some prime, abusing notation by regarding functions on $[N,2N]$ as functions on $\Z/p\Z$ in an obvious way. Let $f: \Z/p\Z \to \C$ be the function
$f(x) := \psi(x) e(-\phi(x))$, and similarly let $\tilde \mu: \Z/p\Z \to \C$ be the function $\tilde \mu(x) := \mu(x)1_{N < x \leq 2N}$.  Then our task is to show
\begin{equation}\label{to-show-1} \E_{x \in \Z/p\Z} \tilde \mu(x) 1_{X_s}(n) f(x) \ll_{A,G/\Gamma} \kappa^{-C}q^2 \log^{-A} N + (\eps + \kappa)\E_{N < n \leq 2N} |\psi(n)| 1_{X_s}(n).\end{equation}
Now let $F: \Z/p\Z \to \C$ be the function defined by
\[ F(h) := 1_{q|h} 1_{B_g(0,\kappa)}(h) e(\phi(n_s+h)).\]  Observe that if $x \in X_s \cap B_g(n_0,\rho)$
and $h_1,h_2 \in B_g(0,\kappa)$ with $q | h_1,h_2$, then from three applications of \eqref{linearity} we have (since $\kappa \leq \rho$)
\begin{align*}
\phi(x+h_1) - \phi(x) - \phi(n_s+h_1) + \phi(n_s) &= O_{\R/\Z}(\eps) \\
\phi(x+h_2) - \phi(x) - \phi(n_s+h_2) + \phi(n_s) &= O_{\R/\Z}(\eps) \; \; \mbox{and} \\
\phi(x+h_1+h_2) - \phi(x+h_1) - \phi(x+h_2) + \phi(x) &= O_{\R/\Z}(\eps),
\end{align*}
where we use $O_{\R/\Z}(\eps)$ to denote a quantity whose $\| \cdot \|_{\R/\Z}$ norm is $O(\eps)$. Summing these three bounds yields
$$ \phi(x) = \phi(x+h_1+h_2) - \phi(n_s+h_1) - \phi(n_s+h_2) + 2 \phi(n_s) + O_{\R/\Z}(\eps),$$
which of course implies that
\[ e(-\phi(x)) = e(-\phi(x + h_1 + h_2))e(\phi(n_s + h_1))e(\phi(n_s + h_2))e(-2\phi(n_s)) + O(\eps).\]
From \eqref{lipmn}, the Lipschitz assumption on $\psi$, we know that $\psi(x + h_1 + h_2) = \psi(x) + O(\kappa)$ for $h_1,h_2 \in B_g(0,\kappa)$. Hence we conclude that
$$ f(x) = f(x+h_1+h_2) F(h_1) F(h_2) e(-2 \phi(n_s)) + O(\eps + \kappa)$$
for all $x \in \Z/p\Z$ and $h_1,h_2 \in B_g(0,\kappa)$ with $q | h_1,h_2$. Since $|f(x)| \leq \psi(x)$ pointwise, we may sum over $X_s$ and deduce that
\begin{align*}
\E_{x \in \Z/p\Z} \tilde \mu(x) 1_{X_s}(x) f(x) &= \E_{\substack{h_1,h_2 \in B_g(0,\kappa) \\ q | h_1,h_2}} \E_{x \in \Z/p\Z} \tilde \mu(x)
f(x+h_1+h_2) F(h_1) F(h_2) e(2 \phi(n_s)) \\
&\quad+ O((\eps + \kappa) \E_{N < n \leq 2N} |\psi(n)| 1_{X_s}(n)).
\end{align*}
To prove \eqref{to-show-1}, then, it suffices to show that
\[ \E_{h_1,h_2 \in B_g(0,\kappa); q | h_1,h_2} \E_{x \in \Z/p\Z} \tilde \mu(x) f(x+h_1+h_2) F(h_1) F(h_2) \ll_{A,G/\Gamma} q^2 \kappa^{-C} \log^{-A} N.\]
From Lemma \ref{bohr-size}(a) and (c) we have
$$ \#\{ h \in B_g(0,\kappa): q |h \} \gg \frac{1}{q} p \kappa^C,$$ and so it is enough to prove that
$$ \E_{h_1,h_2,x \in \Z/p\Z} \tilde \mu(x) f(x+h_1+h_2) F(h_1) F(h_2) \ll_{A} \log^{-A} N.$$
To demonstrate this we use the Fourier transform\footnote{If $g : \Z/p\Z \rightarrow \C$ is a function, and if $\xi \in \Z/p\Z$, we write $\widehat{g}(\xi) := \E_{x \in \Z/p\Z} g(x) e(-x\xi/p)$.} on $\Z/p\Z$, noting in particular the identity 
$$ \E_{x,h_1,h_2 \in \Z/p\Z} \widetilde \mu(x) f(x+h_1+h_2) F(h_1) F(h_2)
= \sum_{\xi \in \Z/p\Z} \widehat{\widetilde \mu}(\xi) \widehat{f}(-\xi) \widehat{F}(\xi)^2.$$
Since $f$ and $F$ are bounded, we see from Plancherel's formula that $|\hat f(-\xi)| = O(1)$ and $\sum_{\xi \in \Z/p\Z} |\widehat F(\xi)|^2 = O(1)$.  
Also, from \eqref{mu-alpha} we have $\widehat{\widetilde \mu}(\xi) \ll_A \log^{-A} N$ for any $\xi$.  The claim follows.
\end{proof}

\begin{remark} What we have in effect done here is approximate $\psi(n) e(-\phi(n))$ by something akin to a \textit{dual function} coming from the Gowers $U^2$-norm. By the general theory of this norm we know that any bounded function which is orthogonal to all linear exponentials (cf. \eqref{mu-alpha}) is orthogonal to all such dual functions. The Fourier argument at the end of the proof of Proposition \ref{ortho-prog} is basically the standard proof of this fact. See \cite{green-tao-inverseu3} for further discussion.
\end{remark}

\begin{remark}
 The results of this section may be used to show that $\mu$ is orthogonal to various other types of function, which need not be Lipschitz or even continuous, but which are still somehow ``approximately linear'' in $n$. Examples of such functions include the bracket-linear phases $e(\beta_1 \lfloor \alpha_1 n\rfloor + \dots + \beta_d \lfloor \alpha_d n\rfloor)$. We omit the details.
\end{remark}

\section{Orthogonality to quadratic phases}\label{quadsec}

 In this section our aim is to prove the estimate \eqref{mu-alphabeta}. Strictly speaking, this section is unnecessary, since \eqref{mu-alphabeta} does not represent the heart of the Main Theorem in the same way that \eqref{mu-alpha} forms the substance of \eqref{mu-1-step}. See the introduction for some remarks on this point.

 This section is included for two pedagogical reasons. First of all the argument does have some features in common with the (far more complicated) analysis of later sections, and thus introduces the main ideas of those sections in a simplified setting. Secondly, it represents a good opportunity to introduce some notation for inequalities which will be very helpful for the rest of the paper.

 The definition of asymptotic orthogonality involves establishing that $X \ll_A \log^{-A} N$, for various quantities $X$ and for all $A > 0$, and it is convenient to have a notation specific to this kind of situation. In each argument that follows, the value of $A$ will be arbitrary, but fixed throughout the argument. When we write $X \lessapprox Y$ or $Y \gtrapprox X$, we mean that \begin{equation}\label{longhand-1} |X| \leq C_AY \log^{C(A+1)} N\end{equation} for some constant $C$ which does not depend on $A$, and some constant $C_A$ which can depend (possibly in an ineffective manner) on $A$. The constants $C$ and $C_A$ can be different in different instances of this notation. In all our arguments the exponent $C$ can be chosen effectively, but it may not be possible to give an explicit value of $C_A$ due to the possibility of Siegel zeros.

 In some cases, statements of the form $X \lessapprox Y$ will appear as both hypotheses and conclusions of a proposition.  In such cases it is understood that the implied constants in the conclusions are dependent on the implied constants in the hypotheses.  Somewhat more subtly, in the course of an argument we may divide into several cases using this notation (e.g. we may divide into
two cases $X \lessapprox Y$ and $X \not \lessapprox Y$).  Once again, the implied constants in the conclusion of this argument will depend on the implied constants used to create the division of cases.  When necessary we shall draw attention to these dependence-of-constants issues\footnote{One can of course rewrite all the arguments in this paper replacing every appearance of $X \lessapprox Y$ or $Y \gtrapprox X$ by suitably explicit long-hand forms \eqref{longhand-1}, although some of the constants may be ineffective.  However we have found that this tended to clutter the estimates with distracting numerical constants, and so we have chosen instead to suppress all of these constants.}.

 Our argument here shall broadly follow that used to prove \eqref{mu-alpha} in \S \ref{linear-sec}.  We begin with the
analogue of Proposition \ref{corr-phase}.

\begin{proposition}[Correlation with quadratic phase implies major arc]\label{corr-phase-quad}  Let $\alpha,\beta,\gamma$ be real numbers, $A > 0$, and let $N$ be a large integer such that
\begin{equation}\label{mu-alpha-fail-quad}
 |\E_{N < n \leq 2N} \mu(n) e(-\alpha n^2 - \beta n - \gamma)| \geq \log^{-A} N.
\end{equation}
Then there exists $D$,  $1 \leq D \ll N^{2/3}$, an integer $q \lessapprox 1$ 
and a $\theta \in \R$ such that
\begin{equation}\label{d-spin-quad}
\#\{ d \in (D,2D]: \Vert q \alpha d^2-\theta\Vert_{\R/\Z} \lessapprox \frac{D^2}{N^2} \} \gtrapprox D.
\end{equation}
Furthermore if $D < N^{1/3}$ we can take $\theta = 0$.
\end{proposition}

\begin{proof} We can discard the constant phase $e(-\gamma)$.  As before, we 
apply Proposition \ref{inverse-prop} with $U=V=N^{1/3}$ and conclude one of the following statements:
\begin{itemize}
\item (Type I sum is large) There exists $D$, $1 \leq D \ll N^{2/3}$, such that
$$ \left|\E_{N/d < w \leq 2N/d} e(\alpha d^2w^2 + \beta dw)\right| \gtrapprox 1$$
for $\gtrapprox D$ values of $d \in (D,2D]$.\vspace{11pt}

\item (Type II sum is large) There exist integers $D,W$ with $N^{1/3} \ll D \ll N^{2/3}$ and $N/4 \leq DW \leq 4N$, such that
$$
 \left|\E_{D < d,d' \leq 2D} \E_{W < w,w' \leq 2W} e(\phi(dw) - \phi(d'w) - \phi(dw') + \phi(d'w'))\right| \gtrapprox 1
$$
where $\phi(n) := \alpha n^2 + \beta n$.
\end{itemize}

 Suppose first that the Type I sum is large.  Applying Lemma \ref{weyl-ineq}, we can find an integer $q \lessapprox 1$
such that $\|q d^2 \alpha\|_{\R/\Z} \lessapprox D^2/N^2$ for $\gtrapprox D$ values of $D < d \leq 2D$, which implies \eqref{d-spin-quad} (with $\theta=0$).

 Now suppose instead that the Type II sum is large.  By the pigeonhole principle, we can find $d', w'$ such that
$$
 |\E_{D < d \leq 2D} \E_{W < w \leq 2W} e(\phi(dw) - \phi(d'w) - \phi(dw') + \phi(d'w'))| \gtrapprox 1$$
and hence
$$ |\E_{W < w \leq 2W} e(\phi(dw) - \phi(d'w) - \phi(dw') + \phi(d'w'))| \gtrapprox 1$$
for $\gtrapprox D$ values of $d$. Now the phase $\phi(dw) - \phi(d'w) - \phi(dw') + \phi(d'w')$ is
quadratic in $w$ with a leading coefficient of $\alpha (d^2 - (d')^2)$.  We may thus apply
Lemma \ref{weyl-ineq} and conclude that there exists $q \lessapprox 1$ such that
\begin{equation}\label{qad}
\|q \alpha (d^2-(d')^2)\|_{\R/\Z} \lessapprox \frac{D^2}{N^2}.
\end{equation}
Pigeonholing in $q$, we conclude there exists a single value of $q$ such that \eqref{qad} follows
for $\gtrapprox D$ values of $d \in (D,2D]$.  Setting $\theta := q\alpha (d')^{2}$, the claim follows.
\end{proof}

 By using Lemma \ref{lem3.1}, we can now conclude the analogue of Corollary \ref{corr-phase-2}.

\begin{proposition}[Correlation with quadratic phase implies major arc, II]\label{corr-phase-quad-2}  Let $\alpha,\beta,\gamma$ be real numbers, $A > 0$, and let $N$ be a large integer such that \eqref{mu-alpha-fail-quad} holds.  Then we have
$$\|\alpha\|_{\R/\Z, Q} \lessapprox N^{-2}$$
for some $Q \lessapprox 1$.
\end{proposition}

\begin{proof}  We apply Proposition \ref{corr-phase-quad} to obtain $D$, $1 \leq D \leq N^{2/3}$, and
$q \lessapprox 1$ obeying \eqref{d-spin-quad}.  If\footnote{This is an instance of the subtlety of the $\lessapprox$ notation. By this we mean that $D \leq C_A \log^{C(A + 1)} N$, where $C$ is chosen so that if $D > C_A \log^{C(A + 1)}N$ then the later estimate \eqref{diagonal-negligible} holds true.} $D \lessapprox 1$
then certainly $D \ll N^{1/3}$, and so we may take $\theta = 0$. There then exists $d \in (D,2D]$ such that
$$ \|q \alpha d^2\|_{\R/\Z} \lessapprox \frac{D^2}{N^2} \lessapprox N^{-2},$$
and the claim follows on replacing $q$ by $qd^2$.

 Now let us suppose that $D \not\lessapprox 1$.  
We will not be able to apply Lemma \ref{lem4.6} as it is not sufficiently ``amplified'' for our use here.
Instead, we use the triangle inequality and \eqref{d-spin-quad} to obtain
$$ \#\{ d,d' \in (D,2D] : \|q \alpha (d^2-(d')^{2})\|_{\R/\Z} \lessapprox \frac{D^2}{N^2} \} \gtrapprox D^2.$$
The diagonal case $d=d'$ is negligible since $D \not\lessapprox 1$, i.e.
\begin{equation}\label{diagonal-negligible}  \#\{ d,d' \in (D,2D] : d \neq d', \|q \alpha (d^2-(d')^{2})\|_{\R/\Z} \lessapprox \frac{D^2}{N^2} \} \gtrapprox D^2.\end{equation}
  Writing $d^2 - (d')^{2} = d_1 d_2$, where $d_1 := d - d'$ and $d_2 := d+d'$,
we conclude
\[ \#\{ d_1,d_2 : 1 \leq |d_1|, |d_2| \leq 4D: \|q \alpha d_1 d_2\|_{\R/\Z} \lessapprox \frac{D^2}{N^2} \} \gtrapprox D^2.\]
By reflection symmetry we may take $d_1,d_2$ to be positive. In particular, for $\gtrapprox D$ values of $d_1$ in $[1,4D]$, 
we have
\[ \#\{ d_2 \in [1,4D]: \|q \alpha d_1 d_2\|_{\R/\Z} \lessapprox \frac{D^2}{N^2} \} \gtrapprox D.\]
Applying Lemma \ref{lem3.1} (ii) we thus conclude that for each such $d_1$, there exists $q_{d_1} \lessapprox 1$ such that
\[ \|q \alpha d_1 q_{d_1} \|_{\R/\Z} \lessapprox \frac{D}{N^2}.\]
Applying the pigeonhole principle, we can thus find $q' \lessapprox 1$ such that
\[ \#\{ 1 \leq d_1 \leq 4D: \|q \alpha d_1 q'\|_{\R/\Z} \lessapprox \frac{D}{N^2} \} \gtrapprox D.\]
Applying Lemma \ref{lem3.1} (ii) again, we conclude that there exists $q'' \lessapprox 1$ such that
\[ \|q \alpha q'' q'\|_{\R/\Z} \lessapprox \frac{1}{N^2}.\]
Since $qq'q'' \lessapprox 1$, the claim follows.
\end{proof}

 On the other hand, we have the quadratic analogue of Proposition \ref{major-orthog}:

\begin{proposition}[Major arc quadratic phases are orthogonal to M\"obius]\label{major-orthog-quad}  
Let $N$ be a large integer, let $\alpha, \beta, \gamma \in \R/\Z$, and let $Q, K \geq 1$
be such that $\|\alpha\|_{\R/\Z,Q} \leq K/N^2$.  Then we have
$$ \E_{N < n \leq 2N} \mu(n) e(-\alpha n^2 - \beta n - \gamma) \ll_{A'} Q^{1/3} K^{1/3} \log^{-A'} N
$$
for any $A' > 0$ \textup{(}the implied constant is ineffective\textup{)}.
\end{proposition}

\begin{proof}  
Let $1 \leq M < N$ be a parameter to be chosen later.  We can set $\gamma = 0$.
Arguing as in the proof of Proposition \ref{major-orthog}, we have
$$  \E_{N < n \leq 2N} \mu(n) e(-\alpha n^2 - \beta n)| \ll
\sup_{|I| = M; I \subset [N,2N]} \big|\frac{1}{M} \sum_{n \in I} \mu(n) e(-\alpha n^2 - \beta n)\big| + \frac{M}{N}.$$
By hypothesis, we have an integer $a$ and $1 \leq q \leq Q$ such that $|\alpha - \frac{a}{q}| \leq \frac{K}{N^2}$.  We thus have
\begin{align*}
e(\alpha n^2) &= e(a n^2/q) e( (\alpha - a/q) n^2) \\
&= e(a n^2/q) e( 2(\alpha - a/q) (n-n_I) ) e( (\alpha - a/q) n^2_I) + O( \frac{K M^2}{N^2} ) \\
&= e(a n^2/q) e( 2(\alpha - a/q) n) \b( \alpha, a/q, n_I) + O( \frac{K M^2}{N^2} )
\end{align*}
for any $n, n_I \in I$, where we use the $\b()$ notation from Appendix \ref{tools-sec}.  
Discarding the constant phase $\b( \alpha, a/q, n_I)$, and absorbing the linear phase
$e( 2(\alpha - a/q) n )$ into the $e(\beta n)$ factor we conclude
\begin{align*}
 \E_{N < n \leq 2N} \mu(n) e(-\alpha n^2-\beta n) &\ll \sup_{|I| = M; I \subset [N,2N]; \beta' \in \R} 
|\frac{1}{M} \sum_{n \in I} \mu(n) e(-a n^2 / q - \beta' n)| \\
&\quad + \frac{K M^2}{N^2}  + \frac{M}{N} .
\end{align*}
The function $e(a n^2/q)$ is periodic of period $q$, and can thus be decomposed as a Fourier series
$e(a n^2/q) = \sum_{b=0}^{q-1} c_b e(bn/q)$
where the coefficients $c_b$ are Gauss sums and can be computed explicitly.  From Plancherel's theorem
and the Cauchy-Schwarz inequality we have $\sum_{b=0}^{q-1} |c_b| = O( q^{1/2} )$ (cf. the proof of Proposition \ref{mob-period}).  
Applying \eqref{mu-alpha} (with $A$ replaced by $2A'$) we conclude that
$$ \sum_{n \in I} \mu(n) e(-a n^2 / q - \beta' n) \ll_{A'} N q^{1/2} \log^{-2A'} N,$$
and hence
$$  \E_{N < n \leq 2N} \mu(n) e(-\alpha n^2 - \beta n - \gamma) \ll_{A'}  \frac{N}{M} q^{1/2} \log^{-2A'} N +
\frac{K M^2}{N^2} + \frac{M}{N}.$$
If we set $M := K^{-1/3} q^{1/6} N \log^{-A'} N$ we obtain the claim.
\end{proof}

 Propositions \ref{corr-phase-quad-2} and \ref{major-orthog-quad} together imply \eqref{mu-alphabeta}, though the $\lessapprox$ notation does take some unravelling. Suppose for a contradiction that \eqref{mu-alpha-fail-quad} holds. Then Proposition \ref{corr-phase-quad-2} implies that $\Vert \alpha \Vert_{\R/\Z,Q} \leq K/N^2$, where we may take $K = Q = C_A \log^{C(A + 1)}N$ for some absolute $C$. Proposition \ref{major-orthog-quad} now implies, taking $A' = C(A + 1)$, that
\[ \E_{N < n \leq 2N} \mu(n) e(-\alpha n^2 - \beta n - \gamma) \ll_{A'} \log^{-A'/3} N.\]
We may clearly assume that $C > 3$, and so this does contradict our assumption that \eqref{mu-alpha-fail-quad} holds, at least if $N > N_0(A)$ is sufficiently large. To conclude the proof of \eqref{mu-alphabeta}, one simply applies Lemma \ref{telescope} with $\varphi \equiv 1$.
\endproof

\begin{remark} It is straightforward to iterate the above argument, as is done in the standard theory of Weyl exponential sums, to obtain
a generalisation of \eqref{mu-alphabeta} in which $\alpha n^2 + \beta n + \gamma$ is replaced by an arbitrary polynomial.  We will, however,
not pursue this generalisation here.  
\end{remark}

\section{Locally quadratic phase functions, I: a technical reduction}\label{sec10}

We now begin the (onerous) task of proving Theorem \ref{mainthm-tech}. Let us begin by recalling the statement: 

\begin{restate-tech-thm}[$\mu$ is strongly orthogonal to local quadratics] Let $G/\Gamma$ be a $1$-step nilmanifold, let
$F: G/\Gamma \to \C$ be a Lipschitz function, and let $g \in G$ and $x \in G/\Gamma$
be arbitrary.  Let $\phi : B_N \rightarrow \R/\Z$ be a phase which is locally quadratic on the \emph{Bohr set} $B_N := \{ n \in [N] : F(T_g^n x) \neq 0\}$. Then we have
\[ \E_{n \in [N]} \mu(n) \overline{F(T_g^n x)} e(-\phi(n)) \ll_{G/\Gamma, A} \Vert F \Vert_{\Lip} \log^{-A} N.\]
\end{restate-tech-thm}

Our objective in this (rather technical) section is to reduce this to a similar result which has certain important technical advantages. The most critical of these is that $\phi$ can be extended to a function which is quadratic somewhat beyond the domain $B_N = \Supp_n F(g^n x)$. This refined formulation reads as follows.

\begin{proposition}[$\mu$ is strongly orthogonal to extendible local quadratics]\label{ortho-prog-quad-2}
Let $g \in G$, $x \in G/\Gamma$, $n_0 \in \Z$, and let $\rho_0 \in (0,10^{-5})$ be a small radius.  Suppose that 
$B_g(n_0,100\rho_0)$ is contained in $\{ n \in \Z: N < n \leq 2N \}$, and suppose that $\phi: \Z \to \R/\Z$ is a function
which is locally quadratic when restricted to $B_g(n_0,100\rho_0)$.  Let $\psi: \Z \to \R^{+}$ be a function supported
on $B_g(n_0,\rho_0)$
which obeys the Lipschitz property
\begin{equation}\label{lip}
|\psi(n) - \psi(m)| \leq \| n - m \|_g \hbox{ for all } n, m \in \Z.
\end{equation}
Then we have
\begin{equation}\label{ortho-targ}
 |\E_{N < n \leq 2N} \mu(n) \psi(n) e(-\phi(n) )| \ll_{A,G/\Gamma}  \log^{-A} N.
 \end{equation}
\end{proposition}

\begin{proof}[Proof that Proposition \ref{ortho-prog-quad-2} implies Theorem \ref{mainthm-tech}] By renormalising, we may assume that $\Vert F \Vert_{\Lip} \leq 1$.

The essential idea is that a ``ball'' (say $B_N$) can be covered by balls $B_g(n_0, \epsilon)$ of a much smaller radius. Most of these will have the property that $B_g(n_0, 100\epsilon)$ is still contained in $B_N$, and hence that $\phi$ is still quadratic on $B_g(n_0, 100\epsilon)$. 

 We turn to the details. First of all, an application of Lemma \ref{telescope} implies that it suffices to establish the estimate
\begin{equation}\label{to-establish} \E_{N \leq n \leq 2N} \mu(n) \varphi(n/N) F(T_g^n x) e(-\phi(n)) \ll_{A, G/\Gamma} \log^{-A} N,\end{equation} where $\varphi : \R \rightarrow \R$ is the function
\[ \varphi(t) := \left\{ \begin{array}{ll} 6(t - \frac{7}{6}) & \mbox{if $\frac{7}{6} \leq t \leq \frac{3}{2}$};\\
6(\frac{11}{6} - t) & \mbox{if $\frac{3}{2} \leq t \leq \frac{11}{6}$}.
\end{array} \right.\]
(any similar function would work).
 The phase $e(\alpha n)$ which featured in that lemma has been absorbed into the quadratic phase $e(-\phi(n))$. 

We now replace $F$ by a ``smooth-thresholded'' function $\widetilde{F}$, as constructed in Lemma \ref{smooth-threshold}. Let $\rho_0 \in (0,10^{-5})$ be a parameter to be chosen later, and set $\delta := 10^4 \rho_0$ in Lemma \ref{smooth-threshold}. This provides a Lipschitz function $\widetilde{F} : G/\Gamma \rightarrow \R$ satisfying properties (i), (ii) and (iii) of that lemma. In particular from Lemma \ref{smooth-threshold} (iii) we see that
\begin{equation}\label{10.1} \E_{N \leq n \leq 2N} \mu(n) \varphi(n/N) F(T_g^n x) = \E_{N \leq n \leq 2N} \mu(n) \varphi(n/N) \widetilde{F}(T_g^n x) + O(\rho_0).\end{equation}
Now take a partition of unity $1 = \sum_{\alpha} \chi_{\alpha}$ on $G/\Gamma$, where
\begin{enumerate}
\item Each $\chi_{\alpha}$ is supported on a ball of diameter at most $\rho_0/2$;
\item Each $\chi_{\alpha}$ is bounded in magnitude by $1$ and satisfies $\Vert \chi_{\alpha} \Vert_{\Lip} \ll_{G/\Gamma} \rho_0^{-1}$;
\item The number of $\chi_{\alpha}$ is $O(\rho_0^{-O_{G/\Gamma}(1)})$.
\end{enumerate}
We leave the construction of such a partition to the reader: modelling $G/\Gamma$ by a torus, one may be quite explicit. This partition of unity induces a decomposition
\[ \widetilde{F} = \sum_{\alpha} F_{\alpha},\]
where $F_{\alpha} := F\chi_{\alpha}$. Note that since both $F$ and $\chi_{\alpha}$ are bounded we have, using Lemma \ref{smooth-threshold} (i), that
\begin{equation}\label{F-alpha-lip} \Vert F_{\alpha} \Vert_{\Lip} \ll \Vert \widetilde{F} \Vert_{\Lip} + \Vert \chi_{\alpha} \Vert_{\Lip} \ll_{G/\Gamma} \rho_0^{-1}.\end{equation}
We may also effect a Lipschitz decomposition
\[ \varphi = \sum_{\beta} \varphi_{\beta}\]
of $\varphi$ into $O(\rho_0^{-1})$ Lipschitz functions $\varphi_{\beta}$ with Lipschitz constant $O(\rho_0^{-1})$, each supported on an interval of diameter $\rho_0/2$.
Write $\psi_{\alpha,\beta}(n) := F_{\alpha}(T_g^n x) \varphi_{\beta}(n/N)$. Noting that 
\[ \varphi(n/N)\widetilde{F}(T_g^n x) = \sum_{\alpha, \beta} \psi_{\alpha,\beta}, \] it follows from \eqref{10.1} and the triangle inequality that 
\begin{align}
\E_{N \leq n \leq 2N} &\mu(n) \varphi(n/N) F(T_g^n x) \\ &\ll_{G/\Gamma} \rho_0^{-O_{G/\Gamma}(1)} \sup_{\alpha,\beta} |\E_{N < n \leq 2N} \mu(n) \psi_{\alpha,\beta}(n) e( -\phi(n) )|  + \rho_0 .\label{10.2}
\end{align}
Suppose that $n,n' \in \Supp(\psi_{\alpha,\beta})$. Then $\varphi_{\beta}(n/N), \varphi_{\beta}(n'/N) \neq 0$, which means that $|n - n'|/N \leq \rho_0/2$. Furthermore $F_{\alpha}(T_g^n x),F_{\alpha}(T_g^{n'}x) \neq 0$, meaning that $\Vert g^{n - n'} \Vert_{G/\Gamma} \leq \rho_0/2$. It follows that $\Vert n - n'\Vert_g \leq \rho_0$, and so the support of $\psi_{\alpha,\beta}$ is contained in some ball $B_g(n_0,\rho_0)$.

We are, of course, going to apply Proposition \ref{ortho-prog-quad-2}. It is therefore necessary to confirm that $\phi$ is defined on $B_g(n_0,100\rho_0)$, and also to say something concerning the Lipschitz constant of $\psi_{\alpha,\beta}$.

Starting with the first task, suppose that $\Supp (\psi_{\alpha,\beta}) \subseteq B_g(n_0,\rho_0)$ and that $\psi_{\alpha,\beta}(n_1) \neq 0$ for some $n_1 \in B_g(n_0,\rho_0)$ (we may clearly ignore those $\alpha,\beta$ for which $\psi_{\alpha,\beta} \equiv 0$). Then $\varphi_{\beta}(n/N) \neq 0$ and so, due to the choice of $\varphi$, we have $7/6 \leq n_1/N \leq 11/6$. It follows that if $n \in B_g(n_0,100\rho_0)$ then $|n - n_1|/N \leq 101\rho_0$ and thus, since $\rho_0$ is so small, that $N < n \leq 2N$. We also have that $F_{\alpha}(g^{n_1} x) \neq 0$, which implies that $\widetilde{F}(g^{n_1} x) \neq 0$. Now if $n \in B_g(n_0,100\rho_0)$ then $d_{G/\Gamma}(g^n x, g^{n_1} x) \leq 101\rho_0$. It follows from Lemma \ref{smooth-threshold} and our choice of $\delta$ that $F(g^n x) \neq 0$. We have shown that $B_g(n_0,100\rho_0) \subseteq B_N$, and hence $\phi$ is indeed defined on the desired set.

We now examine the Lipschitz constant of $\psi_{\alpha,\beta}$, with the $\Vert \cdot \Vert_g$ metric on $\Z$. We have, recalling \eqref{F-alpha-lip}, that
\[ |\varphi_{\beta}(n) - \varphi_{\beta}(n')| \ll \rho_0^{-1} \frac{|n - n'|}{N}\leq \rho_0^{-1} \Vert n - n' \Vert_g \] and
\[ |F_{\alpha}(g^n x) - F_{\alpha}(g^{n'} x)| \ll \rho_0^{-1} \Vert g^{n - n'} \Vert_{G/\Gamma} \leq \rho_0^{-1} \Vert n - n' \Vert_g.\] Since both $F_{\alpha}$ and $\varphi_{\beta}$ are bounded, the Lipschitz constant of $\psi_{\alpha,\beta}$ is $O_{G/\Gamma}(\rho_0^{-1})$.

We are now in a position to apply (a renormalised version of) Proposition \ref{ortho-prog-quad-2}. We deduce that
\[ \E_{N < n \leq 2N} \mu(n) \psi_{\alpha,\beta}(n) e( -\phi(n) ) \ll_{A} \rho_0^{-1}\log^{-A} N  \]
uniformly in $\alpha,\beta$. Thus, from \eqref{10.2}, we see that 
\[ \E_{N \leq n \leq 2N} \mu(n) \varphi(n/N) F(g^n x) \ll_{G/\Gamma} \rho_0^{-O_{G/\Gamma}(1)} \log^{-A} N + \rho_0.\]
Setting $\rho_0 := \log^{-A/2C} N$, and recalling that $A$ can be arbitrary, we do indeed conclude Theorem \ref{mainthm-tech}.
\end{proof}

 It will be convenient later on (in the proof of Lemma \ref{lem12.4}) to add some further technical assumptions to the hypotheses of Proposition \ref{ortho-prog-quad-2}. We may assume that $\psi$ is real. Next, recall that $G/\Gamma$ was embedded isometrically in a torus $(\R/\Z)^d$; we may in fact simply replace $G/\Gamma$ by that torus $(\R/\Z)^d$ (using Lemma \ref{lip-extend}) since this does not affect anything.  It will be convenient to work in $\Z/p\Z$ where $p$ is some prime between $10MN$ and $20MN$. We can approximate the group element $g$ by the nearest $p^\th$ root of unity $\widetilde g$ in $G/\Gamma$, thus $\|g^{-1} \widetilde g\|_{G/\Gamma} \ll 1/N$
and $\tilde g^{p} \in \Gamma$.  Observe that the $\| \cdot \|_g$ and $\| \cdot \|_{\widetilde g}$ norms are comparable, thanks to
the factor of $|\frac{n}{N}|$ in the definition of these norms.  Thus we may, after making some trivial adjustments to the constants such as 100 in the proof of Proposition \ref{ortho-prog-quad-2}, replace $g$ by $\widetilde g$, that is we may
assume that $g$ is a $p^\th$ root of unity.

\section{Locally quadratic phase functions, II: Explicit quadratic and quartic behaviour}\label{localquad-2-sec}

 We now begin the proof of Proposition \ref{ortho-prog-quad-2}.
We are going to show that if \eqref{ortho-targ} is false then the phase $\phi$ is somehow ``major arc''. Ultimately we will relate it to the type of phases in Proposition \ref{ortho-prog} which, in view of the main result of that proposition, will lead to a contradiction. We have already seen several instances where a hypothesis that the M\"obius function $\mu$ correlates with some phase implies that the phase is ``major arc'': Propositions \ref{corr-phase}, \ref{corr-phase-2}, \ref{corr-phase-quad} and \ref{corr-phase-quad-2} are examples of this. In those cases the phase involved, being either linear or quadratic, was of a simple algebraic kind, but the phases that interest us now are not so explicitly given. The two technical lemmas in this section show that these phases do, nevertheless, enjoy some algebraic structure.

 Suppose, for the remainder of the section, that $\phi : B_g(n_0, 100 \rho_0) \rightarrow \R/\Z$ is a locally quadratic phase. If $\|h_1\|_g, \|h_2\|_g \leq 30\rho_0$ then we define
\[   \phi''(h_1,h_2) := \phi(n_0+h_1+h_2) - \phi(n_0+h_1) - \phi(n_0+h_2) + \phi(n_0).\]
This expression is clearly symmetric in $h_1,h_2$.
Since $\phi$ is locally quadratic on the Bohr set $B_g(n_0,100\rho_0)$, we conclude the ``Taylor expansion''
\begin{equation}\label{taylor}
 \phi''(h_1,h_2) = \phi(n+h_1+h_2) - \phi(n+h_1) - \phi(n+h_2) + \phi(n)
\end{equation}
whenever $n \in B_g(n_0,40\rho_0)$.  By telescoping the right-hand side, we conclude the local bilinearity properties
\begin{equation}\label{bilinear}
 \phi''(h_1 + h'_1,h_2) = \phi''(h_1,h_2) + \phi''(h'_1,h_2); \quad \phi''(h_1,h_2+h'_2) = \phi''(h_1,h_2) + \phi''(h_1,h'_2)
\end{equation}
whenever $\|h_1\|_g, \|h_2\|_g, \|h'_1\|_g, \|h'_2\|_g \leq 15\rho_0$.

 As another corollary of Lemma \ref{taylor}, we see that $\phi$ behaves like a genuine quadratic function on certain short arithmetic progressions:

\begin{corollary}[Explicit quadratic structure]\label{quad-grow}  If $n \in B_g(n_0,20\rho_0)$, $L \in \Z$ and $h \in B_g(0, 20\rho_0/L)$,
then there exist $\alpha,\beta \in \R/\Z$ \textup{(}depending on $n$ and $h$\textup{)} such that
\[ \phi(n+hl) = \textstyle\frac{1}{2}\displaystyle l(l-1) \phi''(h,h) + \alpha l + \beta\]
for all $l$, $1 \leq l \leq L$.
\end{corollary}
\proof From \eqref{taylor} we obtain the recurrence
\[ \phi(n+h(l+2)) - 2\phi(n+h(l+1)) + \phi(n+hl) = \phi''(h,h)\]
for all $l$, $1 \leq l \leq L-2$.  The claim follows.\endproof

 This corollary is strong enough for us to understand the behaviour of the Type I sums which will appear when, in subsequent sections,  we analyse 
\[ \E_{n \in [N]} \mu(n) \psi(n) e(-\phi(n))\]
using
Proposition \ref{inverse-prop}. The corresponding Type II sums are more difficult.  The basic issue here is to understand
the algebraic structure of the expression $\psi(dw) e(\phi(dw))$, as a function of $d$ and $w$.  Since $\phi$ is already quadratic, the
phase $\phi(dw)$ here is \emph{quartic} (think of it as being like $d^2 w^2$).  We would like some analogue of Corollary \ref{quad-grow} that makes
this quartic structure manifest, for instance we would like $\phi( (d+sl) (w+tm) )$ to exhibit some explicitly quartic behaviour in
$l$ and $m$, under suitable hypotheses on $d,s,l,w,t,m$ of course.  This turns out to be a little tricky, because of the cross terms
$tdm$ and $slw$ present in the expression $(d+sl) (w+tm)$.  By introducing suitably many constraints (which will be available to us after later arguments) and taking enough differences of the phase, we can eliminate these cross terms and obtain the sought-after quartic structure.

\begin{lemma}[Explicit quartic structure]\label{quartic-grow}  Let $d,w; s, t$ and $L,M$ be integers such that 
\begin{equation}\label{limhk}
 LM \| st \|_g \leq \rho_0
\end{equation}
and let $P: \Z \times \Z \to \Z$
be the quadratic polynomial
\[ P(l,m) := (d+sl) (w+tm).\]
Suppose that the integers $l_0, l_1, l_2, m_0, m_1, m_2$ are such that $|l_i|, |m_i| \leq L$, and furthermore that all sixteen of the values
\begin{equation}\label{sixteen-constraints} P(l_0 + i_1 l_1 + i_2 l_2, j_1 m_1 + j_2 m_2), \; \; i_1, i_2, j_1, j_2 \in \{0,1\},\end{equation} lie in $B_g (n_0, \rho_0)$. Then we have
\begin{equation}\label{four-deriv}
\sum_{i_1,i_2,j_1,j_2 \in \{0,1\}} (-1)^{i_1+i_2+j_1+j_2} \phi( P( l_0 + i_1 l_1 + i_2 l_2, m_0 + j_1 m_1 + j_2 m_2 ) )
= 2 l_1 l_2 m_1 m_2 \phi''(st,st).
\end{equation}
\end{lemma}

\begin{remark} This lemma is a generalisation of the observation that if $\phi(n) = an^2 + bn + c$ is a quadratic, and
one differentiates $\phi(P(l,m))$ twice in the $l$ variable and 
twice in the $m$ variable, one gets $2 \times \phi'' \times st \times st$, where $\phi'' = 2a$ is the double derivative of $\phi$. It is key here that we have the sixteen constraints \eqref{sixteen-constraints}: this gives us sufficient instances where \eqref{taylor} and \eqref{bilinear} may be applied. Later arguments (involving many applications of the Cauchy-Schwarz inequality) will put us in a situation where we have such a multiplicity of constraints at our disposal.
\end{remark}

\proof By replacing $d, w$ by $d+l_0 s$ and $w+m_0 t$ we may assume that $l_0=m_0=0$.
Let $l_1,l_2,m_1,m_2$ be as in the hypothesis of the lemma, that is to say $|l_i|, |m_i| \leq L$ and the sixteen constraints \eqref{sixteen-constraints} are satisfied.  From the identities
\begin{align*}
dw &= P(0,0) \\
sw l_1 = P(l_1,0) - P(0,0), \; \;  & swl_2 = P(l_2, 0) - P(0,0) \\
td m_1 = P(0,m_1) - P(0,0), \; \; & tdm_2 = P(0,m_2) - P(0,0),
\end{align*} we see that
\begin{equation}\label{dweeb}
dw \in B_g(n_0,\rho_0) \quad \mbox{and} \quad sw l_1,swl_2, td m_1, tdm_2 \in B_g(0,2\rho_0).
\end{equation}
Now fix $i_1,i_2 \in \{0,1\}$ and consider the sum
\begin{equation}\label{twoderiv}
\phi( P( i_1 l_1 + i_2 l_2, m_1 + m_2 ) ) - \phi( P( i_1 l_1 + i_2 l_2, m_1 ) ) - \phi( P( i_1 l_1 + i_2 l_2, m_2 ) ) +
\phi( P( i_1 l_1 + i_2 l_2, 0 ) ).
\end{equation}
We can rewrite this as
\begin{equation}\label{use-in-a-few-seconds} \phi(n+h_1+h_2) - \phi(n+h_1) - \phi(n+h_2) + \phi(n)\end{equation}
where $n := w(d+i_1 s l_1 + i_2 s l_2)$, $h_1 := (d + i_1 sl_1 + i_2 sl_2) t m_1$ and $h_2 := (d + i_1 sl_1 + i_2 sl_2) t m_2$. From \eqref{limhk} and \eqref{dweeb} we see that $n \in B_g(n_0, 5\rho_0)$, and that $h_1, h_2 \in B_g(0, 4\rho_0)$. Thus all four of $n, n+h_1, n+ h_2, n+h_1 + h_2$ lie in $B_g(n_0, 13\rho_0)$ and \eqref{taylor} is applicable, which means we can rewrite \eqref{use-in-a-few-seconds} as
\[
\phi''( (d + i_1 sl_1 + i_2 sl_2) t m_1, (d + i_1 sl_1 + i_2 sl_2) t m_2 ).
\]
Applying \eqref{bilinear} and \eqref{dweeb}, \eqref{limhk}, we can expand this as
\[
\eqref{twoderiv} = X + i_1 Y + i_2 Z + 2 i_1 i_2 l_1 m_1 l_2 m_2 \phi''(st,st)\]
where $X,Y,Z$ are quantities which depend on $\phi$, $d,s,t, l_1,m_1,l_2,m_2$ but are independent of $i_1,i_2$.
If one then takes an alternating sum of this identity over the four possible choices of $i_1,i_2 \in \{0,1\}$ to eliminate the $X,Y,Z$ terms,
one obtains \eqref{four-deriv}.
\endproof

\section{Quadratic bias implies major arc}\label{quadmajor-sec}

 With the above preliminaries out of the way, we now begin the proof of Proposition \ref{ortho-prog-quad-2} in earnest.
In this section we shall establish the main step of this proof, namely that a quadratic bias necessarily implies a ``major arc''
condition on $\phi$.  We persist in our use of the notations $X \lessapprox Y$ and $X \gtrapprox Y$, which were introduced in \S \ref{quadsec}. Recall (cf. \eqref{longhand-1}) that $X \lessapprox Y$ means that 
\[ X \leq C_A Y \log^{C(A+1)} N\] for some constant $C$ which does not depend on $A$. That constant is, from now on, allowed to depend on the underlying $2$-step nilmanifold $G/\Gamma$ (in actuality, it will depend on the \emph{dimension} of that nilmanifold).  The constant $C_A$ is of course also allowed to depend on $G/\Gamma$.  Recall also from Appendix \ref{tools-sec} the notation 
\[ \| \alpha \|_{\R/\Z, Q} := \sup_{q \leq Q} \| q \alpha \|_{\R/\Z}.\]
The main result of this section is as follows.

\begin{proposition}\label{circle}  Let the notation and assumptions be as in the previous section.
Suppose that
\begin{equation}\label{fphi-2}
|\E_{N < n \leq 2N} \mu(n) \psi(n) e(-\phi(n))| \geq \log^{-A} N.
\end{equation}
Then there exist $X_0 \lessapprox 1$, $D \leq 4N^{2/3}$ and $Q \lessapprox 1$ with the following property: for any $X$ with $X_0 < X < N^{1/10}$, there exists a set $\mathcal{D} \subseteq [1,D]$, $|\mathcal{D}| \gtrapprox D/X^{1/2}$, such that if $d \in \mathcal{D}$ and $w \in \Z$ satisfies $\| dw \|_g \leq 1/X$ then 
\[ \|\phi''(dw, dw)\|_{\R/\Z, Q} \lessapprox X^{-2}. \]
\end{proposition}

\begin{remark} The conclusion here is an assertion that $\phi''(h,h)$ is major arc for many values of $h$.  We shall recast this
conclusion into a more tractable form in the next section (in particular it is necessary to show that as $d,w$ range over the values allowed in the conclusion of the proposition, $h = dw$ takes on many different values).
\end{remark}

\proof 
Since $\psi$ is Lipschitz and supported on $B(n_0,\rho_0)$, we have $\Vert \psi \Vert_{\infty} \ll \rho_0$, and so we conclude from \eqref{fphi-2} that
\begin{equation}\label{r-big}\rho_0 \gtrapprox 1
\end{equation}
In practice, this will make it fairly easy to verify hypotheses such as $\|h\|_g \leq \rho_0$ which occur in the lemmas of the previous section.

 We now apply Proposition \ref{inverse-prop} with $f(n) := \psi(n) e(\phi(n))$ and $U=V=N^{1/3}$ to conclude one of the following
statements must be true:
\begin{itemize}
\item (Type I sum is large) There exists an integer $1 \leq D \leq N^{2/3}$ such that
\begin{equation}\label{typeI-complex}
|\E_{N/d < w \leq 2N/d} \psi(dw) e(\phi(dw))| \gtrapprox 1
\end{equation}
for $\gtrapprox D$ integers $d$ such that $D < d \leq 2D$.\vspace{11pt}

\item (Type II sum is large) There exists integers $D,W$ with $\frac{1}{2}N^{1/3} \leq D \leq 4N^{2/3}$ and $N/4 \leq DW \leq 4N$, such that
\begin{align}
 |\E_{D < d,d' \leq 2D} \E_{W < w,w' \leq 2W}  & \nonumber\psi(dw) \psi(d'w) \psi(dw') \psi(d'w')\times \\ & \times  e(\phi(dw)-\phi(d'w)-\phi(dw')+\phi(d'w'))| \gtrapprox 1 \label{typeII-complex}
\end{align}
\end{itemize}

 We can thus assume that either \eqref{typeI-complex} or \eqref{typeII-complex} holds, and see what this implies about $\phi$.  We handle
the two cases separately.

\emph{Large Type I sums.} Let us consider the (substantially simpler) Type I case when \eqref{typeI-complex} holds for many values of $D$.  The bulk of the argument is contained inside the following lemma.

\begin{lemma}[Large Type I sum implies major arc]\label{logmajor}  Let $d$, $D \leq d < 2D$, be such that \eqref{typeI-complex} holds, that is to say
\[ 
|\E_{N/d < w \leq 2N/d} \psi(dw) e(\phi(dw))| \gtrapprox 1.\]
  Assume that $N$ is large depending on $A$.  Then there exist $Q \lessapprox 1$ and $\eps \gtrapprox 1$ such that
\[ \|\phi''(dt,dt)\|_{\R/\Z, Q} \lessapprox L^{-2}\]
whenever $L \geq 1$ and $t \in \Z$ is such that $\|dt\|_g \leq \eps/L$.
\end{lemma}

\proof 
The idea is to analyze the quantity in \eqref{typeI-complex} locally on short progressions of common difference $t$ and length $L$.
Since $\psi$ is supported on $(N,2N]$, we have
\[ |\sum_w \psi(dw) e(\phi(dw))| \gtrapprox \frac{N}{D}.\]
From the averaging identity
\[ \sum_w f(w) = \sum_w \E_{1 \leq l \leq L} f(w+tl),\]
valid for any compactly supported function $f: \Z \to \C$, we conclude
\[ |\sum_w \E_{1 \leq l \leq L} \psi(dw + dtl) e(\phi(dw+dtl))| \gtrapprox \frac{N}{D}.\]
Since $\psi$ is supported on $(N,2N]$ and
\[ |dtl| \leq LN |dtl/N| \leq LN \|dt\|_g \leq \eps N,\]
we see that in this sum we still have the constraint $|dw| = O(N)$, and whence $w = O(N/D)$.  Thus by the pigeonhole principle we can find $w$ such that
\[ |\E_{1 \leq l \leq L} \psi(dw + dtl) e(\phi(dw+dtl))| \gtrapprox 1.\]
By \eqref{lip} we have
\[ \psi(dw + dtl) = \psi(dw) + O( l \|dt\|_g ) = \psi(dw) + O( \eps )\]
and hence (if $\eps \gtrapprox 1$ is chosen suitably small)
\[ |\E_{1 \leq l \leq L} \psi(dw) e(\phi(dw+dtl))| \gtrapprox 1.\]
Since $\psi(dw)$ is bounded and independent of $l$, it can be discarded and this becomes
\[ |\E_{1 \leq l \leq L}  e(\phi(dw + dtl)| \gtrapprox  N.\] We apply Corollary \ref{quad-grow} with $n := dw$ and $h := dt$. We may assume, in view of \eqref{r-big}, that $\eps \leq 20 \rho_0$ which means that $h \in B_g(0,20\rho_0/L)$. Since $n \in \Supp(\psi) \subseteq B_g(n_0,\rho_0)$, Corollary \ref{quad-grow} does indeed apply and we may infer the existence of $\alpha, \beta \in \R/\Z$ 
such that 
\[ |\E_{1 \leq l \leq L} e(\textstyle\frac{1}{2}\displaystyle l(l-1) \phi''(dt,dt) + \alpha l + \beta)| \gtrapprox N.\]
Now if $L \geq \log^{C(A + 1)}N$, for sufficiently large $C$, then Lemma \ref{weyl-ineq} applies and we may indeed conclude that $\Vert \phi''(dt,dt)\Vert_{\R/\Z, Q} \lessapprox L^{-2}$. If $L$ is not this large then (because so much may be hidden inside the $\lessapprox$ symbol) the conclusion is trivial anyway.\endproof

 The deduction of Proposition \ref{circle} in the Type I case is almost immediate. Indeed from the preceding lemma we see that for $\gtrapprox D$ values of $d \in [D, 2D)$ we have \[ \|\phi''(dt,dt)\|_{\R/\Z, Q} \lessapprox L^{-2}\]
whenever $t \in \Z$ is such that $\|dt\|_g \leq \eps/L$. Now simply let $\mathcal{D}$ be the set of such $d$, set $L := X/\eps$, and require that $X_0 \lessapprox 1$ be large enough that $L \geq 1$ whenever $X > X_0$. 

\emph{Large Type II sums.} We move on now to the much more complicated Type II case where \eqref{typeII-complex} holds. That is to say, we work under the assumption that 

\begin{align}
 |\E_{D < d,d' \leq 2D} \E_{W < w,w' \leq 2W}  & \nonumber\psi(dw) \psi(d'w) \psi(dw') \psi(d'w')\times \\ & \times  e(\phi(dw)-\phi(d'w)-\phi(dw')+\phi(d'w'))| \gtrapprox 1 \nonumber
\end{align}
where $\frac{1}{2}N^{1/3} \leq D \leq 4N^{2/3}$ and $\frac{1}{4}N \leq DW \leq 4N$.

\begin{lemma}[Type II sum implies major arc]\label{logmajor-2}  Let $\frac{1}{2}N^{1/3} \leq D\leq 4N^{2/3}$ be such that $\frac{1}{4} N \leq DW \leq 4N$ and \eqref{typeII-complex} holds.  
Assume that $N$ is large depending on $A$.  
Then there exist $Q \lessapprox 1$ and $\eps \gtrapprox 1$ with the property that
\[ \| \phi''(st,st) \|_{\R/\Z,Q} \lessapprox 1/L^2 M^2\]
whenever $s,t \in \Z$ and $L,M \in \Z^+$ are such that $L|s| \leq \eps D$, $M|t| \leq \eps W$, $L,M \geq 1/\eps$ and $\Vert st \Vert_g \leq \eps^2/LM$.
\end{lemma}

\proof It will be convenient to use the $\b(x_1,\ldots,x_k)$ notation
introduced in Appendix \ref{tools-sec}.  Thus for instance we can
write \eqref{typeII-complex} as
\[
 |\E_{D < d,d' \leq 2D} \E_{W < w,w' \leq 2W} \psi(dw) e(\phi(dw)) \b(d,w') \b(d',w) \b(d',w')| \gtrapprox 1.
\]
By the pigeonhole principle, we can thus find $d',w'$ such that
\[
 |\E_{D < d \leq 2D} \E_{W < w \leq 2W} \psi(dw) e(\phi(dw)) \b(d,w') \b(d',w) \b(d',w')| \gtrapprox 1.
\]
which upon relabeling the bounded functions $\b$ becomes simply
\[
 |X| \gtrapprox 1
\]
where $X$ is the quantity
\[ X := \E_{D < d \leq 2D} \E_{W < w \leq 2W} \psi(dw) e(\phi(dw)) \b(d) \b(w).\]
Now we argue somewhat as in the proof of Lemma \ref{logmajor}, averaging $d$ \emph{and} $w$ over arithmetic progressions.  For any $1 \leq l \leq L$ and $1 \leq m \leq M$ we can make
the change of variables $d \to d + sl$, $w \to w + tm$ to obtain
\begin{align*} X = \E_{D-sl < d \leq 2D-sl} \E_{W-tm < w \leq 2W-tm} & \psi((d+sl)(w+tm)) \times \\ & \times e(\phi((d+sl)(w+tm))) \b(d+sl) \b(w+tm).\end{align*}
From our assumption that $L|s| \leq \eps D$ and $M|t| \leq \eps W$ we infer that
\begin{align*} X = \E_{D < d \leq 2D} \E_{W < w \leq 2W} \psi((d+sl)(w+tm)) & e(\phi((d+sl)(w+tm))) \times \\ & \times\b(d+sl) \b(w+tm) + O(\eps).\end{align*}
Averaging over $l$ and $m$ gives
\begin{align*} X = \E_{D < d \leq 2D} \E_{W < w \leq 2W} \E_{1 \leq l \leq L} \E_{1 \leq m \leq M}\psi((d+sl)(w+tm)) e(\phi((d+sl)(w+tm))) \times \\ \times \b(d+sl) \b(w+tm)  + O(\eps).
\end{align*}
If $\eps \gtrapprox 1$ is sufficiently small, the assumption that $X \gtrapprox 1$ implies that
\begin{align*} |\E_{D < d \leq 2D} \E_{W < w \leq 2W} \E_{1 \leq l \leq L} \E_{1 \leq m \leq M}
\psi((d+sl)(w+tm)) e(\phi((d+sl)(w+tm))) \times \\ \times \b(d+sl) \b(w+tm)| \gtrapprox 1\end{align*}
Hence by the pigeonhole principle there exist $d,w$ such that
\[ |\E_{1 \leq l \leq L} \E_{1 \leq m \leq M} \psi((d+sl)(w+tm)) e(\phi((d+sl)(w+tm))) \b(d+sl) \b(w+tm)| \gtrapprox 1.\]
Fix such $d,w$.  By relabeling the $\b$'s, we can write $\b(d+sl) \b(w+tm)$ simply as $\b(l) \b(m)$.  We also set
\[ P(l,m) := (d+sl)(w+tm).\]
We have, then, that
\[ |\sum_{l,m} f(l,m) \b(l) \b(m)| \gtrapprox LM\]
where
\begin{equation}\label{fdef}
 f(l,m) := \psi(P(l,m))) e(\phi(P(l,m))) 1_{1 \leq l \leq L} 1_{1 \leq m \leq M}.
 \end{equation}
Using Lemma \ref{cz} to eliminate the $\b(l) \b(m)$ factors, we conclude
$$ |\sum_{l,l',m,m'} f(l,m) \overline{f(l,m')} \overline{f(l',m)} f(l',m')| \gtrapprox L^2 M^2.$$
We write $l=l_0$, $l' = l_0+l_1$, $m=m_0$, $m' = m_0+m_1$ to obtain
\[ |\sum_{l_1,m_1} \sum_{l_0,m_0} F(l_0,m_0;l_1,m_1)| \gtrapprox L^2 M^2\]
where
\[ F(l_0,m_0;l_1,m_1) := f(l_0,m_0) \overline{f(l_0,m_0+m_1)} \overline{f(l_0+l_1,m_0)} f(l_0+l_1,m_0+m_1).\]
Applying Lemma \ref{cz} again, this time in the $l_0$ and $m_0$ variables, we see that
\begin{align} \nonumber |\sum_{l_1,m_1} \sum_{l_0,m_0,l'_0,m'_0} 
F(l_0,m_0;l_1,m_1) \overline{F(l_0,m'_0;l_1,m_1)}  \overline{F(l'_0,m_0;l_1,m_1)}& F(l'_0,m'_0;l_1,m_1)\label{dagger-1}
| \\ & \gtrapprox L^3 M^3.\end{align}
Writing $l'_0 = l_0 + l_2$, $m'_0 = m_0 + m_2$, this becomes
$$ |\sum_{l_0,l_1,l_2,m_0,m_1,m_2} G(l_0,l_1,l_2,m_0,m_1,m_2)
| \gtrapprox L^3 M^3$$
where 
\begin{align*}
G(l_0,l_1,l_2,m_0,m_1,m_2) &:=
F(l_0,m_0;l_1,m_1) \overline{F(l_0,m_0+m_2;l_1,m_1)} \\
&\quad\overline{F(l_0,m_0+m_2;l_1,m_1)} F(l_0+m_2,m_0+m_2;l_1,m_1) \\
&= \prod_{(i_1,i_2,j_1,j_2) \in \{0,1\}^4} \!\!\!{\mathcal C}^{i_1+i_2+j_1+j_2} f(l_0+i_1 l_1 + i_2 l_2, m_0 + j_1 m_1 + j_2 m_2)
\end{align*}
and ${\mathcal C}: z \mapsto \overline{z}$ is the conjugation operator. Observe that the support of the sum in \eqref{dagger-1} is still contained in the region $|l_i| \leq L$, $|m_i| \leq M$. By the pigeonhole principle, we can find $l_0$ and $m_0$ such that
\begin{equation}\label{dagger-2} |\sum_{l_1,l_2,m_1,m_2} \prod_{(i_1,i_2,j_1,j_2) \in \{0,1\}^4} {\mathcal C}^{i_1+i_2+j_1+j_2} f(l_0+i_1 l_1 + i_2 l_2, m_0 + j_1 m_1 + j_2 m_2)
| \gtrapprox L^2M^2.\end{equation}
Let us now expand the product using \eqref{fdef}; this creates a very long product involving sixteen phases (coming from the terms $e(\phi(P(l,m)))$ in the definition of $f$) and fourty-eight cutoffs (coming from the terms $\psi(P(l,m)) 1_{1\leq l \leq L}1_{1 \leq m \leq M}$).  
The sixteen phases $e(\phi(P(l,m)))$ combine to form a single phase
\[ e\big( \sum_{i_1,i_2,j_1,j_2 \in \{0,1\}} (-1)^{i_1+i_2+j_1+j_2} \phi( P( l_0+i_1 l_1 + i_2 l_2, m_0 + j_1 m_1 + j_2 m_2 ) ) \big).\]
The presence of the fourty-eight cutoffs is just what we need to apply Lemma \ref{quartic-grow}, which allows us write the phase in \eqref{dagger-2} as
\[ e\big( 2 l_1 l_2 m_1 m_2 \phi''(st, st)\big).\]
Note that the condition \eqref{limhk} required by that lemma is a consequence of the condition $\Vert st \Vert_g \leq \epsilon^2/LM$ we are working under here, provided that $\eps$ is chosen sufficiently small;  indeed recall from \eqref{r-big} that $\rho_0 \gtrapprox 1$.

 The fourty-eight cutoffs have now served their purpose of explicitly quartilinearising the phase, and we shall now set about obliterating them with further applications of the Cauchy-Schwarz inequality.
To do this, we observe by inspection that fourty-seven of these cutoffs depend on at most three of the variables $l_1,l_2,m_1,m_2$, with the lone
exception being $\psi(P( l_0+l_1+l_2, m_0+m_1+m_2 ))$.  Also, let us recall once more that the cutoffs restrict $l_1,l_2$ to have magnitude at most $L$, and $m_1,m_2$ to have magnitude at most $M$. We thus have
\begin{align}\nonumber 
 |\sum_{\substack{|l_1|,|l_2| \leq L \\|m_1|,|m_2| \leq M}} 
&e(2 l_1 l_2 m_1 m_2 \phi''(st,st)) \psi(P( l_0+l_1+l_2, m_0+m_1+m_2 ))\times \\
&\times \b( l_2,m_1,m_2 ) \b( l_1, m_1, m_2 ) \b( l_1,l_2,m_2) \b(l_1,l_2,m_1)| \gtrapprox L^2M^2.\label{llm}
\end{align}
We would like to eliminate all the $\b()$ factors using Lemma \ref{cz}, but we need to deal with the exceptional cutoff
$\psi(P( l_0+l_1+l_2, m_0+m_1+m_2 ))$ first. First observe that if $\psi$ were a multiplicative function
then the quadratic nature of $P$ would ensure that $\psi(P(l_0+l_1+l_2,m_0+m_1+m_2))$ would factor into the product of expressions, each of which only
depends on at most three (in fact, at most two) of the $l_1,l_2,m_1,m_2$.  Of course, $\psi$ is not multiplicative, but thanks to \eqref{lip}
we can write $\psi(n) = \Psi( g^n, n/N )$ for $N < n \leq 2N$, 
where $\Psi: G/\Gamma \times (\R/\Z) \to \R$ is Lipschitz on the orbit $\{ (g^n, n/N): N < n \leq 2N\}$ and hence, by Lemma \ref{lip-extend}, is the restriction of a Lipschitz function on all of $G/\Gamma \times (\R/\Z)$.  Let $\delta \gtrapprox 1$ be a parameter to be chosen later.
Using Lemma \ref{fourier-lip},
we can approximate $\Psi$ uniformly to accuracy $O(\delta)$ on $(N,2N]$ by a linear combination of at 
most $O(\delta^{-C})$ characters on $G/\Gamma \times (\R/\Z)$, each of which has the form $(x,\theta) \mapsto \chi(x)e(k\theta)$ where $\chi \in (G/\Gamma)^*$ and $k \in \Z$. The coefficients in this linear combination are all $O(1)$.  Thus we can estimate the left-hand side of
\eqref{llm} by
\begin{align*}
& O\big( \delta^{-C} \sup_{\substack{\chi \in (G/\Gamma)^* \\ k \in \Z}} |\sum_{\substack{|l_1|,|l_2| \leq L \\ |m_1|,|m_2| \leq M}} 
e(2 l_1 l_2 m_1 m_2 \phi''(st, st)) \chi(g^{P( l_0+l_1+l_2, m_0+m_1+m_2 )})\times \\
& \times e(kP ( l_0+l_1+l_2, m_0+m_1+m_2 ))\b( l_2,m_1,m_2 ) \b( l_1, m_1, m_2 ) \b( l_1,l_2,m_2) \b(l_1,l_2,m_1)| \big) \\ & \qquad\qquad\qquad+ O( \delta L^2 M^2 ).
\end{align*} Choosing $\delta \gtrapprox 1$ suitably small, we thus conclude that there
exist $\chi$ and $k$ such that the inner sum is $\gtrapprox \delta^C L^2 M^2 \gtrapprox L^2 M^2$. By the quadratic nature of $P$ we may absorb the terms $\chi(g^{P(l_0 + l_1 + l_2, m_0 + m_1 + m_2)})$ and $e(kP ( l_0+l_1+l_2, m_0+m_1+m_2 ))$ into the four unspecified bounded functions $\b()$, thereby obtaining
\begin{align*} 
|\sum_{\substack{|l_1|,|l_2| \leq L \\ |m_1|,|m_2| \leq M}} 
e(2 l_1 l_2 m_1 m_2 \phi''(st,st)) \b( l_2,m_1,m_2 ) \b( l_1, m_1, m_2 ) & \b( l_1,l_2,m_2) \b(l_1,l_2,m_1)| \\
&\gtrapprox L^2 M^2.
\end{align*}
Applying Lemma \ref{cz} to eliminate the $\b()$ factors, we deduce
\[ |\sum_{\substack{|l_1|,|l'_1|,|l_2|,|l'_2| \leq L \\ |m_1|,|m'_1|,|m_2|,|m'_2| \leq M}}
 e(2 (l_1-l'_1) (l_2-l'_2) (m_1-m'_1) (m_2-m'_2) \phi''(st,st)) | \gtrapprox L^4 M^4.\]
By the pigeonhole principle, we can thus find $l'_1, l'_2 = O(L)$ and $m'_1, m'_2 = O(M)$
such that
\[ |\sum_{\substack{|l_1|,|l_2| \leq L \\ |m_1|,|m_2| \leq M}}
 e(2 (l_1-l'_1) (l_2-l'_2) (m_1-m'_1) (m_2-m'_2) \phi''(st,st)) | \gtrapprox L^2 M^2.\]
Summing in $m_2$ using \eqref{exponential-sum}, we obtain
$$ \sum_{|l_1|,|l_2| \leq L} \sum_{|m_1| \leq M}
\min( 1, \frac{1}{M} \| 2 (l_1-l'_1) (l_2-l'_2) (m_1-m'_1) \phi''(st,st) \|_{\R/\Z}^{-1} ) \gtrapprox L^2 M.$$
Shifting $l_1$, $l_2$, $m_1$ by $l'_1$, $l'_2$, $m'_1$ respectively, and doubling $m_1$ to absorb the factor of two this creates, we thus have
$$ \sum_{|l_1|,|l_2| \leq 2L} \sum_{|m_1| \leq 4L}
\min( 1, \frac{1}{M \| l_1 l_2 m_1 \phi''(st,st) \|_{\R/\Z}} ) \gtrapprox L^2 M.$$
It follows that
$$ \min( 1, \frac{1}{M \| l_1 l_2 m_1 \phi''(st,st) \|_{\R/\Z}} ) \gtrapprox 1$$
for $\gtrapprox L^2 M$ triples $(l_1, l_2, m_1)$, which means that
$$ \| l_1 l_2 m_1 \phi''(st,st) \|_{\R/\Z} \lessapprox \frac{1}{M}$$ for those triples. 
In particular, we have $\gtrapprox L^2$ pairs $(l_1,l_2)$ for which this inequality holds for $\gtrapprox M$ values of $m_1 = O(M)$. If $M \geq \log^{C_1(A + 1)}N$ for some sufficiently large $C_1$ then we may apply Lemma \ref{lem3.1} (ii) with parameters $\delta_1 \lessapprox 1/M$, $\delta_2 \gtrapprox 1$ and $|I| \sim M$ to conclude that for each such pair $(l_1,l_2)$, there exists $q \lessapprox 1$ such that
\begin{equation}\label{12.5} \| l_1 l_2 q \phi''(st,st) \|_{\R/\Z} \lessapprox \frac{1}{M^2}.\end{equation}
This condition on $M$ may be met by choosing $\eps \gtrapprox 1$ sufficiently small, since one of the hypotheses of the lemma was that $M \geq 1/\eps$.
Applying the pigeonhole principle to \eqref{12.5}, we can now locate a single $q \lessapprox 1$ such that the above bound
holds for $\gtrapprox L^2$ pairs $(l_1, l_2)$.

 Taking $\eps$ sufficiently small we may assume that $L, M \geq \log^{C_2(A+1)}N$ for suitable $C_2$ and apply Lemma \ref{lem3.1} to $l_2$ instead of $m_1$. The parameters in that lemma are now $\delta_1 \approx 1/M^2$, $\delta_2 \approx 1$ and $|I| \sim L$, and we conclude the existence of $q' \lessapprox 1$ such that
$$ \| l_1 q' q \phi''(st,st) \|_{\R/\Z} \lessapprox \frac{1}{LM^2}$$
for $\gtrapprox L$ values of $l_1$.  Applying Lemma \ref{lem3.1} one last time, now with $\delta_1 \approx 1/LM^2$, $\delta_2 \approx 1$ and $|I| \sim L$,  we find a $q'' \lessapprox 1$
such that
$$ \| q'' q' q \phi''(st,st) \|_{\R/\Z} \lessapprox \frac{1}{L^2 M^2}.$$
Since $q''q'q' \lessapprox 1$, the proof of Lemma \ref{logmajor-2} is complete.
\endproof

 It remains to use this lemma to complete the proof of Proposition \ref{circle} in the Type II case. We take $\mathcal{D}$ to be simply the whole interval $[D/2X^{1/2}, D/X^{1/2}]$. There is a very important subtlety here: this set of integers can only be guaranteed to have size $\gtrapprox D/X^{1/2}$ if we assume that $D/X^{1/2} \gg 1$. Note, however, that in the Type II case this is so since we are working under that assumption that $D \gg N^{1/3}$ and $X \leq N^{1/10}$. This is not just a technical artefact of our approach -- it is simply not possible to bound a general bilinear form, such as the Type II sum
\[ T_{\mbox{\scriptsize II}} = \sum_{d \sim D} \sum_{w \sim W} a_d b_w f(dw)\] when one of the ranges $d \sim D$ or $w \sim W$ is too short, as the weights $a_d, b_w$ could conspire to give no cancellation.

 Suppose, then, that $d \in \mathcal{D}$ and that $w \in \Z$ satisfies the condition of Proposition \ref{circle}, namely that $\Vert dw \Vert_g \leq 1/X$. In Lemma \ref{logmajor-2} take $L = M := \eps X^{1/2}/10$ and $s := d$, $t := w$. If $X_0 \lessapprox 1$ is sufficiently large and $X > X_0$ then certainly the two conditions $L, M \geq 1/\eps$ are satisfied. Furthermore we have
\[ L|s| \leq \frac{\eps X^{1/2}}{10} \cdot \frac{D}{X^{1/2}} \leq \eps D\] and
\[ M|t| = \frac{\eps X^{1/2}}{10} \cdot \frac{|dw|}{|d|} \leq \frac{\eps X^{1/2}}{10} \cdot \frac{N \Vert dw \Vert_g}{|d|} \leq \frac{\eps N}{5D} \leq \eps W,\] and finally $\Vert st \Vert_g \leq \eps^2/LM$ by the definition of $L$ and $M$. All the conditions of Lemma \ref{logmajor-2} are thus satisfied, and we may infer that 
\[ \Vert \phi''(dw, dw) \Vert_{\R/\Z, Q} \lessapprox X^{-2}\] 
for some $Q \lessapprox 1$, as required.\endproof

 We may now forget about Type I and II sums, and work with the conclusion of Proposition \ref{circle} instead. In the next section we will use divisor moment estimates to cast this conclusion in a more tractible form.

\section{Massaging the major arc condition}\label{massage-sec}

 In the last two sections we showed that if $\psi(n) e(-\phi(n))$ correlates with M\"obius (specifically if \eqref{fphi-2} holds true) then $\phi$ must exhibit some kind of ``major arc'' behaviour. Indeed we proved Proposition \ref{circle}, which we urge the reader to recall now. Our first task in this section is to cast the conclusion of that proposition in a more useable form. Through this section, we assume that $\phi : B_g(n_0, 100\rho_0) \rightarrow \R/\Z$ is a phase for which \eqref{fphi-2}, and hence the conclusion of Proposition \ref{circle}, holds true.

\begin{proposition}\label{prop12.1}  Let $\phi$ be as above, and suppose that the parameter $\rho_1$ satisfies \begin{equation}\label{rho-1-assump} N^{-c} < \rho_1 < \rho_0\log^{-C_1 (A + 1)}N \end{equation} for some $c > 0$ and some $C_1$ which is sufficiently large depending on $G/\Gamma$ \textup{(}in reality $\rho_1$ will be much larger than $N^{-c}$, so the lower bound here is hardly relevant\textup{)}. 
Then 
\begin{equation}\label{phinn}
\|\phi''(n, n)\|_{\R/\Z, Q} \lessapprox \rho^2_1
\end{equation}
for $\gtrapprox \rho_1^{3/2}|B_g(0,\rho_1)|$ values of $n \in B_g(0,\rho_1)$, where $Q \lessapprox 1$.
\end{proposition}

\begin{remarks} Note that since $\rho_1^{3/2}$ is so much bigger than $\rho_1^2$, the conclusion is in the spirit of the hypotheses of Lemmas such as \ref{lem4.6}, where a quadratic whose fractional part was ``close to zero unexpectedly often'' was shown to be major arc. We will, in fact, apply exactly that lemma later in this section. The fact that we can arrange the exponents $3/2$ and $2$ in this way is ultimately due to the lower bound $|\mathcal{D}| \gtrapprox D/X^{1/2}$ in Proposition \ref{circle}; $|\mathcal{D}| \gtrapprox D/X$ would not suffice.
\end{remarks}

\begin{proof} Set $X := 1/\rho_1$ in Proposition \ref{circle}; we may certainly suppose that $C_1$ is so large that this is permissible. We find $D \ll N^{2/3}$ and a set ${\mathcal D} \subseteq \{ 1, \ldots, D\}$ of cardinality $|{\mathcal D}| \gtrapprox D/X^{1/2}$ such that
$$ \|\phi''(dw, dw)\|_{\R/\Z, Q} \lessapprox \rho_1^2$$
whenever $d \in {\mathcal D}$ and $w \in \Z$ are such that $dw \in B_g(0,\rho_1)$.  Thus, if we define the sets
$$ \Omega := B_g(0,\rho_1) \cap \Z^+; \quad \Omega_d := \{ n \in \Omega: d|n \}$$
for each integer $d > 1$, it will suffice (noting that $B_g(0,\rho_0)$ is symmetric about the origin) to prove the estimate
\begin{equation}\label{divpack-conclusion}
\big| \bigcup_{d \in {\mathcal D}} \Omega_d \big| \gtrapprox \rho_1^{3/2} |\Omega|.
\end{equation}Observing from Lemma \ref{bohr-size} that
$$ |\Omega| \gg \rho_1^C N \hbox{ and } |\Omega_d| \gg \frac{1}{D} |\Omega| \hbox{ for all } d \in {\mathcal D},$$ 
where $C$ depends only on $G/\Gamma$, 
it follows by taking $\kappa := 1/2C$ in Lemma \ref{div-pack} of Appendix \ref{div-section} that
\[ \big| \bigcup_{d \in {\mathcal D}} \Omega_d \big| \gg \frac{|\mathcal{D}|^2}{|D|^2} |\Omega| \rho_1^{1/2} \log^{-C_2} N.\]
Since $|\mathcal{D}| \gtrapprox D/X^{1/2}$, the result follows immediately.\end{proof}

 As we remarked, the conclusion of Proposition \ref{prop12.1} has the form ``$\phi''(n,n)$ is surprisingly close to an integer very often'' on a small Bohr set $B_g(0,\rho_1)$.  The next step is to amplify this to obtain $\phi''(h,h)$ major arc for a significantly larger set of $h$
 (working on $B_g(0,\rho_0)$ rather than $B_g(0,\rho_1)$).
More precisely, we now establish a more pleasant characterisation of major arc:

\begin{lemma}[Major arcs have small second derivative]\label{major-second}  Let $\phi$ be as above.  Then there
exists $Q_1 \lessapprox 1$ such that
$$ \|\phi''(h, h)\|_{\R/\Z, Q_1} \lessapprox \| h \|_g^2$$
for all $h \in B_g(0,\rho_0)$.
\end{lemma}

\begin{proof} 
The idea is to make the quadratic structure of $\phi''$ so explicit that we can apply Lemma \ref{lem4.6}.

We will choose $\rho_1 := \log^{-C_2 (A+1)} N$, where $C_2 \gg C_1$ is some constant to be specified later. In particular if $C_2$ is large enough then the conditions of Proposition \ref{prop12.1} are satisfied, and we can find as a result some set $\mathcal S \subseteq B_g(0,\rho_1)$ such that 
\begin{equation}\label{s-1} |\mathcal{S}| \gtrapprox \rho_1^{3/2} |B_g(0,\rho_1)|\end{equation} and 
\begin{equation}\label{s-2} \Vert \phi''(n,n) \Vert_{\R/\Z,Q} \lessapprox \rho^2_1\end{equation} for all $n \in \mathcal S$. Note that the implied constants in the $\gtrapprox$ and $\lessapprox$ notations here \emph{do not depend on $C_2$}. Note also that for reasons like this one must exercise extreme caution with these notations.

 Select some $C_3 \gg C_2$. If $\Vert h \Vert_g \geq \log^{-C_3(A + 1)} N$ then the lemma holds vacuously, and so we assume henceforth that $\Vert h \Vert_g \leq \log^{-C_3(A + 1)} N$. 
Now from  \eqref{s-1} and Lemma \ref{bohr-size}(b) we have 
$$ \E_{n \in B_g(0,2\rho_1)} 1_{\mathcal{S}}(n + m) \gtrapprox \rho_1^{3/2}$$
for all $m \in B_g(0,\rho_1)$.  Applying this to $m=hl$ for all $l \in \{1,\dots,L\}$, where $L := \lfloor \frac{\rho_1}{\|h\|_g} \rfloor$, and then averaging
in $L$, we conclude
$$ \E_{n \in B_g(0,2\rho_1)} \E_{1 \leq l \leq L} 1_{\mathcal{S}}(n + hl) \gtrapprox \rho_1^{3/2},$$
and thus by the pigeonhole principle we can find $n \in B_g(0,2\rho_1)$ such that
$$ \E_{1 \leq l \leq L} 1_{\mathcal{S}}(n + hl) \gtrapprox \rho_1^{3/2}.$$
In particular, we have
$$ \|\phi''(n+hl, n+hl)\|_{\R/\Z, Q} \lessapprox \rho_1^2$$
for $\gtrapprox \rho_1^{3/2}L$ values of $l \in \{1,\dots,L\}$.  Applying the pigeonhole principle again, we can thus find a single $q \lessapprox 1$ such that
$$ \|q \phi''(n+hl, n+hl)\|_{\R/\Z} \lessapprox \rho_1^2$$
for $\gtrapprox \rho_1^{3/2} L$ values of $l \in \{1,\dots,L\}$.  Now from Corollary \ref{quad-grow} (and \eqref{r-big}) we can write
$$ q \phi''(n+hl, n+hl) = q l^2 \phi''(h,h) + \alpha l + \beta $$
for some quantities $\alpha,\beta \in \R/\Z$ which depend on $q,\phi,n,h$ but are independent of $l$.  Thus
\[ \Vert ql^2 \phi''(h,h) + \alpha l + \beta \Vert_{\R/\Z} \lessapprox \rho_1^2\] for $\gtrapprox \rho_1^{3/2} L$ values of $l \in \{1,\dots,L\}$. Now Lemma \ref{lem4.6} applies to exactly this kind of situation. In that lemma we take $\delta_1 \approx \rho_1^{2}$ and $\delta_2 \approx \rho_1^{3/2}$, and note that the requisite conditions $\delta_1 \leq \frac{1}{4}\delta_2$ and $L \geq 2^{58} \delta_2^{-12}$ are handsomely satisfied if $C_2,C_3$ are chosen judiciously. The conclusion is that 
\[ \Vert q \phi''(h,h) \Vert_{\R/\Z, 2^{-43} \delta_2^{-9}} \leq 2^{141} \delta_2^{-28} L^{-2}.\] Setting $Q_1 := 2^{-43}\delta_2^{-9} \lessapprox 1$ and noting that $L \gtrapprox \Vert h \Vert_g^{-1}$, the conclusion follows.\end{proof}

 In the next lemma, we bootstrap Lemma \ref{major-second} to a depolarized version of itself. 

\begin{lemma}[Major arcs have small second derivative, II]\label{major-2}  Let $\phi$ be as above.  Then there exist $Q_2 \lessapprox 1$ and $\rho_2 \gtrapprox 1$
such that $\rho_2 \leq \rho_0$ and
$$ \|\phi''(h, h')\|_{\R/\Z, Q_2} \lessapprox \| h \|_g \| h' \|_g$$
for all $h, h' \in B_g(0,\rho_2)$.
\end{lemma}

\begin{proof}  Let $\rho_2 = \log^{-C_4 (A+1)} N$, for some large $C_4$ to be chosen later.
By symmetry we may assume $\|h'\|_g \leq \|h\|_g$.  Let $L > 1$ be the least integer such that $L \|h'\|_g > \|h\|_g$.
For any $l \in \{1,\dots,L\}$, we use \eqref{bilinear} and the hypotheses $h,h' \in B_g(0,\rho_2)$ to conclude
$$ 4 l \phi''(h, h') = \phi''(h + lh', h+lh') - \phi''(h-lh', h-lh').$$
Applying Lemma \ref{major-second} and the triangle inequality, we infer
$$ \| 4 l \phi''(h, h') \|_{\R/\Z, Q_1} \lessapprox \| h \|_g^2$$
and hence
$$ \| l \phi''(h, h') \|_{\R/\Z, 4Q_1} \lessapprox \| h \|_g^2$$
for all $l \in \{1,\dots,L\}$.
Let $C_5$ be a further constant to be specified later. If $L \leq \log^{-C_5(A + 1)} N$ then we can set $l=1$ and the argument is finished. Suppose, then, that $L \geq \log^{-C_5(A + 1)} N$.  By the pigeonhole principle, we can find $q \leq Q_1 \lessapprox 1$ such that
$$ \| q l \phi''(h, h') \|_{\R/\Z}  \lessapprox \| h \|_g^2$$
for $\gtrapprox L$ values of $l \in \{1,\dots,L\}$.  We are now in a position to apply Lemma \ref{lem3.1}(ii) with $\delta_1 \approx \Vert h \Vert_g^2$ and $\delta_2 \approx 1$. If $C_4$ is large enough then we certainly have $\delta_1 \leq \frac{1}{4}\delta_2$, whilst $C_5$ may be chosen so that $L > 2/\delta_2^2$. In those circumstances the lemma is applicable and we deduce that 
$$ \| q \phi''(h,h') \|_{\R/\Z} \lessapprox \|h\|_g^2 / L \leq \Vert h \Vert_g \Vert h' \Vert_g.$$
This concludes the proof.
\end{proof}

 The above lemma says that for any pair $h,h'$ each having small $\|\cdot\|_g$ norms, the second derivative $\phi''(h,h')$ is close to
a rational number $a/q$ for some small $q$.  However, this $q$ can currently depend on $h,h'$.  Fortunately, it is possible to ``clear denominators''
and make $q$ independent of $h,h'$, by taking advantage of a certain ``finite dimensionality'' of the Bohr set $B_g(0,\rho_2)$.  More precisely, we have

\begin{lemma}[Major arcs have small second derivative, III]\label{lem12.4}  Let $\phi$ be as above.  Then there exists $\rho_3 \gtrapprox 1$ and an integer $q \lessapprox 1$ such that
$$ \|q\phi''(h, h')\|_{\R/\Z} \lessapprox \| h \|_g \| h' \|_g$$
for all $h, h' \in B_g(0,\rho_3)$.
\end{lemma}
\begin{proof}  
We shall use some standard results from the geometry of numbers to obtain a ``basis'' for the Bohr set
$B_g(0,\rho_2)$. These result are discussed in several places: see, for example, \cite{bilu, green-tao-inverseu3} and \cite[Ch. 3]{tao-vu-book}. Recall at this point the discussion at the end of \S \ref{sec10}, where we remarked that $g \in (\R/\Z)^d$ can be taken to be an $p^\th$ root of unity, where $p \in [10N, 20N]$ is the prime we have associated to $N$ for those arguments where it is convenient to work in a cyclic group. This is such an argument.  We identify $B_g(0,\rho_2)$, which is certainly contained in $\{1,\dots,N\}$, with a subset of $\Z/p\Z$. Write 
\[ g = (\frac{\xi_1}{p},\dots, \frac{\xi_d}{p})\]
in $(\R/\Z)^d$, where $\xi_1,\dots, \xi_d \in \Z/p\Z$. Let $S \subseteq \Z/p\Z$ be the set of frequencies
\[ S := \{1,\xi_1,\dots, \xi_d\}.\]
In the notation of \cite{green-tao-inverseu3}, the Bohr set $B_g(0,\rho_2)$ is then comparable to a ``traditional''Bohr set 
$$B(S,\rho) := \{ x \in \Z/p\Z: \| \xi x / p \|_{\R/\Z} < \rho \}$$
in the sense that
\begin{equation}\label{bohr-compatibility} B(S,\textstyle\frac{1}{20}\displaystyle\rho_2) \subseteq B_g(0,\rho_2) \subseteq B(S,2\rho_2).\end{equation}
Applying \cite[Corollary 10.5]{green-tao-inverseu3}, and redefining $d := d+1$, we can then find a proper\footnote{By \emph{proper} we mean
that all the sums $l_1 v_1 + \ldots + l_d v_d$ are distinct.} generalized arithmetic progression
$$P = \{ l_1 v_1 + \ldots + l_d v_d: |l_j| \leq L_j \hbox{ for all } 1 \leq j \leq d \}$$
for some $L_1,\ldots,L_d \ge 1$ and $v_1,\ldots, v_d \in \Z/p\Z$, such that
$$ B_g(0, c \rho_2) \subseteq P \subseteq B_g(0,\rho_2)$$
for some $c = c(d) > 0$. In fact by applying that result to $B_g(0,\frac{1}{4}\rho_2)$ (and redefining $P$ and the $L_j$ slightly) we may insist on the slightly stronger inclusions
\begin{equation}\label{inclusions} B_g(0, \textstyle\frac{1}{4}\displaystyle c\rho_2) \subseteq P_{1/4} \subseteq P \subseteq B_g(0,\rho_2)\end{equation}
where $P_\theta$ is defined for any $\theta \in (0,1]$ by
$$P_\theta := \{ l_1 v_1 + \ldots + l_d v_d: |l_j| \leq \theta L_j \hbox{ for all } 1 \leq j \leq d \}.$$
We will prove the lemma with $\rho_3 := \frac{1}{4}c\rho_2$. Let us note from \eqref{inclusions} that 
$$ \Vert v_j \Vert_g \leq \frac{\rho_2}{L_j} \leq \frac{1}{L_j}$$
for each $j$, $1 \leq j \leq d$.  Thus by Lemma \ref{major-2} we may find for each $j,j'$,  $1 \leq j, j' \leq d$, a $q_{j,j'} \lessapprox 1$
such that
$$ \|q_{j,j'} \phi''(v_j, v_{j'})\|_{\R/\Z} \lessapprox \frac{1}{L_j L_{j'}}.$$
If we let $q$ be the least common multiple of all the $q_{j,j'}$, then we still have $q \lessapprox 1$ and
$$ \|q \phi''(v_j, v_{j'})\|_{\R/\Z} \lessapprox \frac{1}{L_j L_{j'}}$$
for all $j,j'$, $1 \leq j, j' \leq d$. Note that at this point the implied constants in the $\lessapprox$ notation have become heavily dependent on $d$.  By bilinearity \eqref{bilinear} it follows that
\begin{equation}\label{star-eq} \|q \phi''(h,h')\|_{\R/\Z} \lessapprox \| h \|_P \| h' \|_P\end{equation}
for all $h, h' \in P$, where the norm $\Vert \cdot \Vert_P$ on $P$ is defined by
$$ \| l_1 v_1 + \ldots + l_d v_d \|_P := \sup_{1 \leq j \leq d} \frac{|l_j|}{L_j}.$$
We claim that $\Vert h \Vert_P \lessapprox \Vert h \Vert_g$ for all $h \in B_g(0,\rho_3)$. In view of \eqref{star-eq}, this will suffice to prove the lemma.

 We may assume that $h \neq 0$ since the claim is trivial otherwise.  Observe that $h \in P_{1/2}$.
Let $M > 1$ be the smallest positive integer such that $Mh \not \in P_{1/4}$; since $Mh = (M-1)h + h$, we see that $Mh \in P_{1/2}$. Thus $\|Mh\|_P \leq 1/2$, which implies that $\|h\|_P \leq 1/2M$ (here we use the
hypothesis that $P$ is proper, which implies that the co-ordinates $l_1,\ldots,l_d$ of $Mh$ are $M$ times the co-ordinates of $h$).  
On the other hand, since $Mh \not\in P_{1/4}$, we have $Mh \not \in B_g(0,\rho_3)$, which implies $M \|h\|_g \geq \rho_3$ and hence that $\Vert h \Vert_g \geq \rho_3/M \gtrapprox 1/M$. Combining these estimates we obtain the claim, and hence the lemma.
\end{proof}

\section{Handling the major arcs}\label{major-sec}

 Let us summarise the current state of affairs. In our effort to prove Proposition \ref{ortho-prog-quad-2}, we assumed that its conclusion \eqref{ortho-targ} was false. After a long and complicated analysis, we deduced from this assumption that the phase $\phi$ is \emph{major arc}, in the sense that we have an estimate
\[ \Vert q \phi''(h,h') \Vert_{\R/\Z} \lessapprox \Vert h \Vert_g \Vert h' \Vert_g\] whenever $h,h' \in B_g(0,\rho_3)$, for some $q \lessapprox 1$ and some $\rho_3 \gtrapprox 1$. This was, of course, the content of Lemma \ref{lem12.4}. To close the argument, we relate major arc phases of this type to those appearing in Proposition \ref{ortho-prog}. This is not hard (though a little technical), and leads quickly to a contradiction (of the assumption that \eqref{ortho-targ} was false).

 Let $q$ be as above.  By bilinearity \eqref{bilinear} again, we see that
$$ \|\phi''(h, h')\|_{\R/\Z} \lessapprox \| h \|_g \| h' \|_g$$
for all $h, h' \in B_g(0,\rho_3)$ such that $q | h, h'$.  Let $\eps \ll \rho_3$, $\eps \gtrapprox 1$, be a small number to be chosen later.
Applying \eqref{taylor}, we conclude the approximate linearity relationship
\begin{equation}\label{star-1} \| \phi(n+h_1+h_2) - \phi(n+h_1) - \phi(n+h_2) + \phi(n)\|_{\R/\Z} \lessapprox \eps^2 \ll \epsilon \end{equation}
whenever $n \in B_g(n_0,2\rho_0)$, whenever $h_1,h_2 \in B_g(0,20\eps)$ are such that $q | h_1,h_2$, and provided that $\eps$ is small enough.

 Now due to the finite dimensionality of the space $(\R/\Z)^d \times \R$ from which the metric $\Vert n - m \Vert_g$ is naturally descended (cf. the remarks following Definition \ref{bohr-set}) we may cover $B_g(n_0,\rho_0)$ with $O(\eps^{-C})$ Bohr sets $B_g(n_\alpha,\eps)$ such that each point is contained in $O(1)$ of these Bohr sets.  
This induces a corresponding partition of $\psi$
into $O(\eps^{-C})$ functions $\psi_\alpha$, each of which is supported on a Bohr set $B_g(n_\alpha,\eps)$ and still obeys the Lipschitz bound \eqref{lip}.

 Now observe that if $n, n+h_1, n+h_2, n+h_1+h_2 \in B_g(n_{\alpha}, 10\eps)$, and if $q | h_1, h_2$ then \eqref{star-1} holds. Thus we may apply Proposition \ref{ortho-prog} (with $\rho = \eps$) to conclude that for any $\kappa \leq \epsilon$, and for some $A'$ to be chosen later, we have
$$ |\E_{N < n \leq 2N} \mu(n) \psi_\alpha(n) e( -\phi(n) )| \ll_{A'} \kappa^{-C} q^3 \log^{-A'} N + (\eps + \kappa) \E_{N < n \leq 2N} |\psi_\alpha|.$$
Summing in $\alpha$, using the bounded overlap of the Bohr sets and the fact that $\Vert \psi \Vert_{\infty} \ll \rho_0 \ll 1$, we conclude
\begin{equation}\label{dagger-eq} |\E_{N < n \leq 2N} \mu(n) \psi(n) e( -\phi(n) )| \ll_{A'} (\eps\kappa)^{-C} q^3 \log^{-A'} N + \eps + \kappa.\end{equation}
At this point we set\footnote{We kept the parameters $\eps$ and $\kappa$ separate in Proposition \ref{ortho-prog} for pedagogical reasons, to make the dependencies clear.} $\kappa = \eps = \log^{-C(A + 1)} N$ for some $C > 1$ which is so large that \eqref{star-1} holds. Recalling that $q \lessapprox 1$, we see that $A'$ may be chosen so that the right-hand side of \eqref{dagger-eq} is $\ll \log^{-A} N$.

 We have, at long last, contradicted the supposition that \eqref{ortho-targ} is false. This implies Proposition \ref{ortho-prog-quad-2}. By the analysis of \S \ref{sec10} Theorem \ref{mainthm-tech} is also true, and thus, by the deduction immediately after the statement of Theorem \ref{mainthm-tech}, so is the Main Theorem.\endproof

\appendix

\section{Some harmonic analysis tools}\label{tools-sec}

 In this appendix we collect some simple harmonic analysis tools which are used frequently in the paper. We begin by introducing some
norms on the unit circle $\R/\Z$, which can be lifted up to the real line $\R$.  

\begin{definition}[Circle norms]  If $\alpha$ is an element of the real line $\R$ or the circle $\R/\Z$, we use $\|\alpha\|_{\R/\Z}$ to denote
the distance from $\alpha$ to the nearest integer (if $\alpha$ is real) or to zero (if $\alpha$ is on the circle $\R/\Z$).  If $Q \geq 1$ is an integer, we use $\|\alpha\|_{\R/\Z,Q}$ to denote the quantity
$$ \|\alpha\|_{\R/\Z,Q} := \inf_{1 \leq q \leq Q} \|q\alpha\|_{\R/\Z}.$$
\end{definition}

 The quantity $\|\alpha\|_{\R/\Z}$ is subadditive, thus $\|\alpha + \beta\|_{\R/\Z} \leq \|\alpha\|_{\R/\Z} + \|\beta\|_{\R/\Z}$.  We caution however
that the quantity $\|\alpha\|_{\R/\Z,Q}$ (which is large when $\alpha$ lies in a ``minor arc'', and small when $\alpha$ lies in a ``major arc'') is \emph{not} subadditive.

 Define a \emph{discrete interval} to be any set of the form $\{ n \in \Z: a \leq n \leq b\}$ for some $a,b$.
By summing the geometric series, we observe the elementary exponential sum estimate
\begin{equation}\label{exponential-sum}
|\sum_{n \in I} e(\alpha n)| \leq 4\min( |I|, \frac{1}{ \|\alpha\|_{\R/\Z} } )
\end{equation}
for any discrete interval $I \subset \Z$ and any $\alpha \in \R/\Z$ (or any $\alpha \in \R$).  One consequence of this is the following
P\'olya-Vinogradov type completion of sums lemma, which allows one to estimate a partial sum by a completed sum at the cost of a logarithm
and an exponential phase.

\begin{lemma}[Completion of sums]\label{complete}  Let $I \subset \Z$ be a discrete interval, and $f: \Z \to \C$ be a function.  Then we have
$$ \sup_{J \subseteq I} |\sum_{n \in J} f(n)| \ll \log(1 + |I|) \sup_{\alpha \in \R/\Z} |\sum_{n \in I} f(n) e(\alpha n)|$$
where the supremum on the left ranges over discrete sub-intervals of $I$.
More generally, if $I' \subset \Z$ is another discrete interval, and $K: \Z \times \Z \to \C$ is a function, then we have
$$ \sum_{m \in I'} |\sum_{n \in I} 1_{J_m}(n) K(n,m)|^2 \ll \log^2(1 + |I|) \sup_{\alpha \in \R/\Z}  \sum_{m \in I'} |\sum_{n \in I} K(n,m) e(\alpha n)|^2$$
where for each $m \in I'$, $J_m \subset \Z$ is an arbitrary discrete interval.
\end{lemma}

\begin{proof}  We may assume $I$ is non-empty.  By translation we may take $I = \{1,\ldots,L\}$ for some $L \geq 1$, which we then identify with $\Z/L\Z$.  If $J$ is any interval in $\Z/L\Z$, we can use Fourier expansion in $\Z/L\Z$ to write
\begin{align*}
\sum_{n \in J} f(n) &= \sum_{n \in \Z/L\Z} 1_J(n) f(n) \\
&= \sum_{\xi \in \Z/L\Z} \widehat{1_J}(\xi) \sum_{n \in \Z/L\Z} e(n\xi/L) f(n)
\end{align*}
where 
$$ \widehat{1_J}(\xi) := \E_{n \in \Z/L\Z} 1_J(n) e(-n\xi/L).$$
Applying \eqref{exponential-sum}, we have
\[ |\widehat{1_J}(\xi)| \leq 4 \min( 1, \frac{1}{L \|\xi/L\|_{\R/\Z}} ).\]
Thus by the triangle inequality, we have
\begin{align*}
|\sum_{n \in J} f(n)| 
&\leq 4 \sum_{\xi \in \Z/L\Z} \min\big( 1, \frac{1}{L \|\xi/L\|_{\R/\Z}} \big) \bigg| \sum_{n \in \Z/L\Z} e(n\xi/L) f(n) \bigg| \\
&\ll \sum_{\xi \in \Z/L\Z} \min\big( 1, \frac{1}{L \|\xi/L\|_{\R/\Z}} \big) \sup_{\alpha \in \R/\Z} | \sum_{n \in I} e(n\alpha) f(n) | \\
&\ll \log(1+L) \sup_{\alpha \in \R/\Z} | \sum_{n \in I} e(n\alpha) f(n) |,
\end{align*}
which gives the first inequality.  Using similar arguments, as well as the triangle inequality in $l^2$, we have
\begin{align*} 
(\sum_{m \in I'} & |\sum_{n \in I} 1_{J_m}(n) K(n,m)|^2)^{1/2}  \\ &\ll
(\sum_{m \in I'} (\sum_{\xi \in \Z/N\Z} \min( 1, \frac{1}{L \|\xi/L\|_{\R/\Z}} ) |\sum_{n \in I} e(n\xi/L) K(n,m)| )^2)^{1/2} \\
&\ll \sum_{\xi \in \Z/N\Z} \min( 1, \frac{1}{L \|\xi/L\|_{\R/\Z}} ) (\sum_{m \in I'} |\sum_{n \in I} e(n\xi/L) K(n,m)|^2)^{1/2} \\
&\ll \log(1+L) \sup_{\alpha \in \R/\Z} (\sum_{m \in I'} |\sum_{n \in I} e(\alpha \xi) K(n,m)|^2)^{1/2}.
\end{align*}
\end{proof}

In a similar spirit, we now recall the well-known Erd\H{o}s-Tur\'an inequality:

\begin{proposition}[Erd\H{o}s-Tur\'an inequality]\label{et}
Let $(u_l)_{l=1}^L$ be a sequence in $\R/\Z$, and define the \emph{discrepancy} $\Delta(\alpha,\beta)$ for 
any $-\frac{1}{2} \leq \alpha < \beta < \frac{1}{2}$ by the formula
\[ \Delta(\alpha,\beta) := \# \{l \in \{1,\dots,L\} : u_l \in [\alpha,\beta] \} - (\beta - \alpha)L.\]
Then for any positive integer $Q$ we have
\[ |\Delta(\alpha,\beta)| \leq \frac{L}{Q} + 3 \sum_{q=1}^Q \frac{1}{q} \big| \sum_{l=1}^L e(qu_l)\big|.\]
\end{proposition}

\begin{proof}  See for instance \cite{montgomery}.  The constant $3$ is unimportant for us, and could be improved slightly. 
\end{proof}

 An important application of this inequality for us (which we will use extremely frequently) 
will be the following observation, which says that if a linear sequence $\alpha l$ stays close
to an integer for many $l$ in an interval $I$, then $\alpha$ must be ``major arc'', in the sense that $\|\alpha\|_{\R/\Z,Q}$ is small for some small $Q$.

\begin{lemma}[Recurrent linear functions are major arc]\label{lem3.1}  Let $I \subseteq \Z$ be a discrete interval, let $\alpha \in \R/\Z$,
and suppose that the set
$$
{\mathfrak L} := \{ l \in I: \| \alpha l \|_{\R/\Z} \leq \delta_1 \} 
$$
has cardinality at least $\delta_2 |I|$
for some $0 < \delta_1, \delta_2 < 1$ with $\delta_1 \leq \frac{1}{4} \delta_2$.
\begin{itemize}
\item[(i)] If $|I| > 1/\delta_2$, then 
$\displaystyle \|\alpha\|_{\R/\Z, 8/\delta_2} \leq \frac{2^8}{\delta_2^2 |I|}$.
\item[(ii)] If $|I| > 2/\delta_2^2$, then 
$\displaystyle \|\alpha\|_{\R/\Z, 16/\delta_2^2} \leq \frac{2^{15} \delta_1}{ \delta_2^6 |I|}$.
\end{itemize}
\end{lemma}

\begin{proof}  Write $I = \{ M+1,\ldots,M+L\}$, and let $(u_l)_{l=1}^L$ be the sequence $u_l := \alpha(M+l) \md 1$.  Then the lower bound
on ${\mathfrak L}$ implies the 
discrepancy estimate
$$ \Delta(-\delta_1, \delta_1) \geq (\delta_2 -2\delta_1)L \geq \frac{1}{2}\delta_2 L.$$

Let us now prove (i).  Applying Proposition \ref{et} we conclude
\[ \frac{1}{2}\delta_2 L \leq \frac{L}{Q} + 3 \sum_{q = 1}^Q \frac{1}{q}\big| \sum_{l=1}^L e(qu_l) \big|\]
for any $Q$.  Taking $Q =: \lceil 4/\delta_2\rceil$, this implies that there is $q \leq 8/\delta_2$ such that 
\[ \big| \sum_{l=1}^L e(qu_l) \big| \geq 2^{-6} \delta_2^2 L.\]
Applying \eqref{exponential-sum}, the result follows.

 We now use a standard ``amplification'' argument, exploiting the smallness of $\delta_1$ compared to $\delta_2$,
to bootstrap (i) to the stronger estimate (ii).  We may assume that $\delta_1 < \delta_2^2/16$ since the result follows immediately from (i)
otherwise.   Let $1 \leq m \leq L$ be an integer to be chosen later; then by the pigeonhole principle and the lower bound on
$|{\mathfrak L}|$, there exists some $b$ such that the set
$$ \mathfrak{L}_b := \{b+1,\dots, b+ m\} \cap \mathfrak{L}$$
has cardinality at least $\delta_2 m/2$.  We fix $b$, and note that if $x \in m \mathfrak{L} + \mathfrak{L}_b$, that is to say if $x = ml + l'$ with $l \in \mathfrak{L}$ and $l' \in \mathfrak{L}_b$, then $\|\alpha x\|_{\R/\Z} \leq 2m\delta_1$. Furthermore we have $|m \mathfrak{L} + \mathfrak{L}_b| \geq \delta_2^2 mL/2$, and also $m \mathfrak{L} + \mathfrak{L}_b$ is a subset of the interval \[ I' := \{m(M+1) + b + 1,\dots, m(M+L) + b + m\},\] which has cardinality at most $mL$.  We can apply (i) with $I$, $\delta_1$, $\delta_2$ replaced by $I'$, $2m\delta_1$, and $\delta_2^2/2$, provided that
$m \leq \delta_2^2/16\delta_1$ and $mL > 2/\delta_2^2$. It being sensible to take $m$ essentially as large as possible, set $m := \lfloor \delta_2^2/16\delta_1\rfloor$. The result follows quickly.
\end{proof}

 Next, we record a version of summation by parts.  Define the \emph{total variation} $\|\psi\|_{\TV}$ of a sequence $\psi: \Z \to \C$ to be the quantity
$$ \| \psi\|_{\TV} := \sup_{n \in \Z} |\psi(n)| + \sum_{n \in \Z} |\psi(n+1)-\psi(n)|,$$
and more generally define the \emph{total variation modulo $q$} for any $q \geq 1$ to be the quantity
$$ \| \psi\|_{\TV,q} := \sup_{n \in \Z} |\psi(n)| + \sum_{n \in \Z} |\psi(n+q)-\psi(n)|.$$

\begin{lemma}[Summation by parts]\label{sum-parts}  If $f, \psi: \Z \to \C$ and $I$ is an interval, then
$$ |\sum_{n\in I} f(n) \psi(n)| \leq \| \psi \|_{\TV} \sup_{J \subseteq I} |\sum_{n\in J} f(n)|.$$
More generally, for any $q \geq 1$ we have
$$ |\sum_{n\in I} f(n) \psi(n)| \leq q \|\psi \|_{\TV,q} 
\sup_{J \subseteq I, a \in \Z/q\Z} |\sum_{n\in J} f(n) 1_{n=a \md{q}}|.$$
\end{lemma}

\begin{proof} Write $I = \{u,\dots,v\}$, and denote by $S_n := \sum_{j = u}^n f(j)$ the partial sums of $f$. Recalling the summation by parts formula
\[ \sum_{n \in I} f(n) \psi(n) = S_v \psi(v) + \sum_{n = u}^{v-1} S_n (\psi(n) - \psi(n+1)),\]
the first inequality follows immediately. The second bound follows by splitting $I$ into $q$ residue classes modulo $q$
and applying a rescaled version of the first identity to each component.
\end{proof}

\begin{corollary}[Completion of sums, II]\label{complete-2}  Let $I \subset \Z$ be a discrete interval, and $f: \Z \to \C$ and $\psi: \Z \to \C$ be functions.  Then we have
$$ \sum_{n \in I} \psi(n) f(n) \ll \log(1 + |I|) \| \psi \|_{\TV} \sup_{\alpha \in \R/\Z} |\sum_{n \in I} f(n) e(\alpha n)|$$
and more generally for any $q \geq 1$
$$ \sum_{n \in I} \psi(n) f(n) \ll q\log(1 + |I|) \| \psi \|_{\TV,q} \sup_{\alpha \in \R/\Z} |\sum_{n \in I} f(n) e(\alpha n)|.$$
\end{corollary}
\begin{proof} The first part is immediate from Lemmas \ref{complete} and \ref{sum-parts}. To obtain the second bound, we begin with an invocation of the second bound in Lemma \ref{sum-parts}. It is now sufficient to prove that
\[ \sup_{J \subseteq I, a \in \Z/q\Z} |\sum_{n \in J} f(n) 1_{n \equiv a \md{q}} | \leq  \sup_{\alpha \in \R/\Z} |\sum_{n \in I} f(n) e(\alpha n)|.\]
To see this, expand $1_{n \equiv a \md{q}}$ as a Fourier series
\[ 1_{n \equiv a \md{q}} = \frac{1}{q} \sum_{\xi \in \Z/q\Z} e\bigg( \frac{(a - n)\xi}{q}\bigg),\] and apply Lemma \ref{complete} and the triangle inequality.
\end{proof}

 As a consequence of this Corollary, we can obtain the following convenient lemma, 
which allows us to replace the range $1 \leq n \leq N$ by a smooth cutoff
to the interval $N < n \leq 2N$, at the expense of adding an arbitrary linear phase to the function (which in our applications will be totally harmless). 

\begin{lemma}\label{telescope}
  Let $f: \N \to \C$ be a sequence bounded by $O(1)$.
  Let $\varphi: \R \to \R$ be a Lipschitz non-negative function of Lipschitz norm $O(1)$ which is at least 1 on $[4/3,5/3]$.  Suppose that we know that
$$ \E_{N < n \leq 2N} \varphi(\frac{n}{N}) f(n) e(\alpha n) \ll_A \log^{-A} N$$
for all $A > 0$, $N \geq 1$, and $\alpha \in \R/\Z$.  Then we have
$$ \E_{n \in [N]} f(n) \ll_{A,\varphi} \log^{-A} N$$
for all $A > 0$ and $N \geq 1$.
\end{lemma}

\begin{proof} For large $N$ we can write
$$ \E_{4N/3 < n\leq 5N/3} f(n) \ll_\varphi |\E_{N < n \leq 2N} \varphi(\frac{n}{N}) f(n) g(n)|$$
where 
$$ g(n) = 1_{4N/3 < n \leq 5N/3} \varphi^{-1}(\frac{n}{N}).$$
Since $\varphi^{-1}$ is Lipschitz on $[4/3,5/3]$, we have $\|g\|_{\TV} \ll_\varphi 1$, and hence by Corollary \ref{complete-2} and
hypothesis
\begin{equation}\label{eq700} \E_{4N/3 < n\leq 5N/3} f(n) \ll_{\varphi} \sup_{\alpha \in \R/\Z} |\E_{N < n \leq 2N} \varphi(\frac{n}{N}) f(n) e(\alpha n)|
\ll_A \log^{-A} N.\end{equation}
Now we may decompose the interval $\{1,\dots,N\}$ into $O(\log N)$ intervals of type $4M/3 < n \leq 5M/3$ together with $O(\log N)$ extra points. Combining \eqref{eq700} with the bound $f = O(1)$, we obtain the lemma.
\end{proof}

Another harmonic analysis tool we will need often is to approximate Lipschitz functions by exponentials.  We first recall a well-known extension lemma: 

\begin{lemma}[Lipschitz extension]\label{lip-extend} 
If $Y$ is a non-empty subset of a metric space $X = (X,d)$, and $f: Y \to \R$ is a Lipschitz function then there exists a Lipschitz extension $f_{\operatorname{ext}}: X \to \R$ of $f$ from $Y$ to $X$ with $\Vert f_{\operatorname{ext}} \Vert_{\Lip} = \Vert f \Vert_{\Lip}$.  Similarly, if $f: Y \to \C$ is Lipschitz then there exists an extension $f_{\operatorname{ext}}: X \to \C$ with $\Vert f_{\operatorname{ext}} \Vert_{\Lip} \leq 2 \Vert f \Vert_{\Lip}$.
\end{lemma}

\begin{proof} If $f$ is real-valued one can for instance define $f_{\operatorname{ext}}(x) := \min( \inf \{ f(y) + M d(x,y): y \in Y \}, \sup_{y \in Y} f(y))$, where $M := \Vert f \Vert_{\Lip}$. The complex case then follows by splitting $f$ into real and imaginary parts.
\end{proof}

\begin{lemma}[Fourier approximation of Lipschitz functions]\label{fourier-lip} 
Let $(\R/\Z)^d$ be the standard $d$-dimensional torus, with metric induced by the $l^\infty$ norm
\begin{equation}\label{torus-infty}
 \| (x_1,\ldots,x_d) \|_{(\R/\Z)^d} := \sup_{1 \leq j \leq d} \| x_j \|_{\R/\Z}.
\end{equation}
Let $Y$ be a subset of $(\R/\Z)^d$, and
let $f: Y \to \C$ be a Lipschitz function bounded in magnitude by $1$.  Then for any $N \geq 1$ there
exist $J = O_{d}(N^d)$, $c_1,\ldots,c_J = O(1)$, and $m_1,\ldots,m_J \in \Z^d$ such that
$$ f(x) = \sum_{j=1}^J c_j e(m_j \cdot x) + O_{d}\big(\frac{\Vert f \Vert_{\Lip} \log N}{N}\big)$$
for all $x \in Y$.  Furthermore, the values of $m_1,\ldots,m_J$ depend on $L$, $d$, $N$ but are otherwise independent of $f$ or $Y$.
\end{lemma}

\begin{proof} By Lemma \ref{lip-extend} we may take $Y = (\R/\Z)^d$.  Let $\sigma_N: (\R/\Z)^d \to \R^+$ be the Fej\'er kernel
$$ \sigma_N(x_1,\ldots,x_d) := \prod_{j=1}^d \frac{1}{N} \frac{\sin^2(\pi N x_j)}{\sin^2(\pi x_j)}.$$
Note that
\[\widehat{\sigma}_N(m) = \prod_{j = 1}^d \big( 1 - \frac{|m_j|}{N}\big) 1_{|m_j| \leq N}\]
for all $m \in \Z^d$. 
We have
\[ f \ast \sigma_N(x) = \sum_{m} \widehat{f \ast \sigma_N}(m) e(m \cdot x) = \sum_m \widehat{f}(m) \widehat{\sigma}_N(m) e(m \cdot x)\]
which, since $\Vert f \Vert_{\infty} = O(1)$, has the form $\sum_{j = 1}^J c_j e(m_j \cdot x)$ where $J = O_d(N^d)$ and $c_j = O(1)$. To conclude the proof of the lemma, then, it suffices to show that $\Vert f  - f \ast \sigma_N \Vert_{\infty} = O_{d}(\Vert f \Vert_{\Lip}\log N/N)$. To this end, note that
\[
|f(x) - f \ast \sigma_N(x)|  =  \big| \int_{(\R/\Z)^d} (f(x) - f(y)) \sigma_N(x - y)\, dy \big|,\]
and hence by the change of variables $z := x-y$ it will suffice to show that
\[
\int_{(\R/\Z)^d} \| z \|_{(\R/\Z)^d} \sigma_N(z)\, dz = O_d(\log N/N).\]
Since $\sigma_N$ has total mass one, the portion of the integral on the region $\| z \|_{(\R/\Z)^d} \leq N^{-1}$ is acceptable.
Now, for each integer $n \geq 0$, consider the portion of the integral on the annular region $2^n N^{-1} \leq \| z \|_{(\R/\Z)^d} \leq 2^{n+1} N^{-1}$.
We have
\begin{eqnarray*} \big| \int_{\Vert z \Vert_{(\R/\Z)^d} \sim 2^n N^{-1}} \|z\|_{(\R/\Z)^d} \sigma_N(z)\, dz \big| & \ll & 2^n N^{-1} \int_{\Vert t \Vert_{(\R/\Z)^d} \gg 2^n N^{-1}} |\sigma_N(t)|\, dt \\ & \ll_d & 2^n N^{-1} \int_{\Vert t_1 \Vert_{\R/\Z} \gg 2^n N^{-1}} \frac{1}{N} \frac{\sin^2 (\pi N t_1)}{\sin^2 (\pi t_1)}\, dt_1 \\ & \ll_d & 2^n N^{-1} 
\int_{\Vert t_1 \Vert_{\R/\Z} \gg 2^n N^{-1}} \frac{1}{N \| t_1 \|_{\R/\Z}^2}\, dt_1 \\
& \ll_d & \frac{1}{N}.
\end{eqnarray*}
Summing this over $n = 0,1,\dots,N$ we obtain the claim.
\end{proof}

We shall adopt the following convenient notation from \cite{green-tao-inverseu3}: we use $\b(x_1,\ldots,x_k)$ to denote any function of
the variables $x_1,\ldots,x_k$ which is bounded by $O(1)$; the exact value of $\b()$ may vary from line to line, just as with the $O()$ notation.
We use this notation to denote functions whose exact value is not of interest to us, invariably because they are destined to be annihilated in the course of a Cauchy-Schwarz argument such as the following.

\begin{lemma}[Cauchy-Schwarz inequality]\label{cz}  Let $X, Y$ be finite non-empty sets, and let $f: X \times Y \to \C$ be a function.  Then
$$ |\E_{x \in X} \E_{y \in Y} \b(x) f(x,y)| \ll |\E_{x \in X} \E_{y,y' \in Y} f(x,y) \overline{f(x,y')}|^{1/2}$$
and
$$ |\E_{x \in X} \E_{y \in Y} \b(x) \b(y) f(x,y)| \ll |\E_{x,x' \in X} \E_{y,y' \in Y} f(x,y) \overline{f(x,y')} \overline{f(x',y)} f(x',y')|^{1/4}.$$
Similarly, if $K: X^4 \to \C$ is a function, then
\begin{align*}
 &|\E_{x_1,x_2,x_3,x_4 \in X} \b(x_2,x_3,x_4) \b(x_1,x_3,x_4) \b(x_1,x_2,x_4) \b(x_1,x_2,x_3) K(x_1,x_2,x_3,x_4)| \\
 &\quad \ll \big|\E_{x_{1,0}, x_{1,1}, \ldots, x_{4,0}, x_{4,1} \in X} \prod_{i_1,i_2,i_3,i_4 \in \{0,1\}} {\mathcal C}^{i_1+\ldots+i_4}
 K(x_{1,i_1}, \ldots, x_{4,i_4})\big|^{1/16}
\end{align*}
where ${\mathcal C}: z \mapsto \overline{z}$ is the conjugation operator.
\end{lemma}

\begin{remark} These estimates are part of the theory of the Gowers uniformity norms $\|f\|_{U^d}$ and $\| K \|_{\Box^d}$; see for instance
\cite{green-tao-linearprimes,gowers-long-aps,green-fin-field,green-tao-primes,green-tao-inverseu3,tao:survey}. 
\end{remark}

\begin{proof}  From the triangle inequality and Cauchy-Schwarz we have
$$ |\E_{x \in X} \E_{y \in Y} \b(x) f(x,y)| \ll \E_{x \in X} |\E_{y \in Y} f(x,y)| \leq (\E_{x \in X} |\E_{y \in Y} f(x,y)|^2)^{1/2}$$
and the first claim follows.  The second claim follows by two iterations of the first, and the third follows from four iterations of the first.
\end{proof}

 Now, we develop some quadratic analogues to the linear phase estimates given above.
We begin with a quadratic counterpart to
\eqref{exponential-sum}. We do not pretend that the exponents here are even remotely optimal; we have opted for a statement which is conveniently derived from our earlier lemmas.

\begin{lemma}[Weyl's inequality]\label{weyl-ineq}
Let $\alpha,\beta,\gamma \in \R$ and let $\delta \in (0,1)$. Let $I \subset \Z$ be a discrete interval such that
$|I| \geq 2^{16}/\delta^6$ and
\[ \big| \E_{l \in I} e(\alpha l^2 + \beta l + \gamma) \big| \geq \delta.\]
Then we have
$$\|\alpha\|_{\R/\Z,2^{12}\delta^{-4}} \leq \frac{2^{43}}{\delta^{14}|I|^2}.$$
\end{lemma}

\begin{proof}  By translating $I$ we may take $I = \{1,\ldots,L\}$ for some $L$.
Squaring the expression gives a double sum over variables $l',l$; setting $l' = l+ h$, we find that
\[ |\sum_{h = -L}^{L} \sum_{l = \max(1 - h,1)}^{\min(L - h,L)} e(2\alpha h l + \alpha h^2 + \beta h)| \geq \delta^2 L^2.\] 
Summing the inner geometric series using \eqref{exponential-sum} we see that
\[ \sum_{h = -L}^L \min\big(L, \frac{1}{\|2\alpha h\|_{\R/\Z}}\big) \geq \delta^2 L^2/2\] and therefore that
\[ \sum_{h = 1}^L \min\big(L, \frac{1}{\|2\alpha h\|_{\R/\Z}}\big) \geq \delta^2 L^2/8.\]
It follows that there are at least $\delta^{2} L/16$ values of $h \in \{1,\dots,L\}$ such that $\|2\alpha h\|_{\R/\Z} \leq 16/\delta^2 L$. 
The claim then follows from Lemma \ref{lem3.1}(ii).
\end{proof}

 One can now repeat the proof of Lemma \ref{lem3.1}(i), using Lemma \ref{weyl-ineq} in place of \eqref{exponential-sum}, to conclude

\begin{lemma}[Recurrent quadratics are non-diophantine]\label{lem4.6}  Let $I \subseteq \Z$ be a discrete interval, let $\alpha, \beta, \gamma$ be real numbers, and suppose that the set
$$
\{ l \in I: \| \alpha l^2 + \beta l + \gamma \|_{\R/\Z} \leq \delta_1 \}
$$
has cardinality at least $\delta_2 |I|$ for some $0 < \delta_1, \delta_2 < 1$ with $\delta_1 \leq \frac{1}{4} \delta_2$.
If $|I| \geq 2^{58} \delta_2^{-12}$, then we have
$$\|\alpha\|_{\R/\Z,2^{43} \delta_2^{-9}} \leq 2^{141} \delta_2^{-28} |I|^{-2}.$$
\end{lemma}

The final tool we assemble in this appendix is a technical lemma used in \S \ref{sec10}. This allows us to approximate a Lipschitz function $F$ by a ``soft-thresholded'' function $\widetilde{F}$.

\begin{lemma}[Soft-thresholding a Lipschitz function]\label{smooth-threshold}
Let $F : X \rightarrow [-1,1]$ be any Lipschitz function on a metric space $(X,d)$, and let $\delta > 0$ be a parameter. Then there is a Lipschitz function $\widetilde{F} : X \rightarrow [-1,1]$ satisfying the following properties:
\begin{enumerate}
\item $\Vert \widetilde{F} \Vert_{\Lip} \leq \Vert F \Vert_{\Lip}$;
\item If $x \in \Supp(\widetilde{F})$ and $d(x,x') \leq \delta$ then $x' \in \Supp(F)$;
\item $\Vert F - \widetilde{F} \Vert_{\infty} \leq \delta \Vert F \Vert_{\Lip}$. \end{enumerate}
\end{lemma}
\proof We will set
\[ \widetilde{F}(x) := \max( |F(x)| - \lambda,0) \sgn(F(x))\] for an appropriate value of $\lambda \geq 0$ which we shall shortly specify. Let us first prove that \emph{any} such function satisfies (i). Since $|\widetilde{F}|$ is pointwise bounded by $|F|$, it suffices to show that if $x,x' \in X$ then
\[ | \widetilde{F}(x) - \widetilde{F}(x') | \leq |F(x) - F(x')|.\]
But this follows because the function $x \mapsto \max( |x| - \lambda, 0) \sgn(x)$ is easily seen to be a contraction. This proves (i).

Now set $\lambda := \delta \Vert F \Vert_{\Lip}$. Statement (iii) is then obvious. To prove (ii), note that if $x \in \Supp(\widetilde{F})$ then $|F(x)| > \lambda$. Thus if $d(x,x') \leq \delta$ then
\begin{equation}\boxeq |F(x')| \geq |F(x)| - |F(x) - F(x')| \geq |F(x)| - \delta \Vert F \Vert_{\Lip} > 0.\end{equation}

\section{Nilsequences and locally polynomial phases}\label{appendixA}

The purpose of this appendix if to give the proof of Proposition \ref{nil-to-local}, the statement of which we recall now.

\begin{nil-to-local-repeat}[$2$-step nilsequences are averages of twisted $1$-step nilsequences]  Let $G/\Gamma$ be a $2$-step nilmanifold and let $0 < \eps < 1/2$.
Let $F: G/\Gamma \to \C$ be a bounded Lipschitz function with $\Vert F \Vert_{\Lip} \leq 1$, and let $g \in G$ and $x \in G/\Gamma$ be arbitrary.  Then there exists a $1$-step nilmanifold $\widetilde G/\widetilde \Gamma$ depending only on $G/\Gamma$ and a decomposition \[
F(T_g^n x) = \E_{i \in I} w_i F_i(T_{g_i}^n x_i) e(-\phi_i(n)) + O(\eps)
\]
where
\begin{itemize}
 \item $I$ is a finite index set;
\item For each $i \in I$ the $w_i$ are complex numbers with $\E_{i \in I} |w_i| \ll \eps^{-O_{G/\Gamma}(1)}$;
\item $F_i: \tilde G/\tilde \Gamma \to \C$ is bounded $O_{G/\Gamma}(1)$-Lipschitz;
\item $g_i \in \tilde G$;
\item $x_i \in \tilde G/\tilde \Gamma$;
\item 
$\phi_i: B_i \to \R/\Z$ is a phase function which is locally quadratic on the \emph{generalized Bohr set} 
$B_i := \{ n \in [N] : F_i(T_{g_i}^n x_i) \neq 0\}$.
\end{itemize}
\end{nil-to-local-repeat}

As we remarked in \S \ref{technical-sec}, we are going to give a rather hands-on calculational approach to this theorem, using Mal'cev bases and the Heisenberg nilmanifold as an illustrative example. The reader interested in a comprehensive discussion of Mal'cev bases may consult the book \cite{corwin-greenleaf}.

Let $G$ be a connected, simply connected, $2$-step nilpotent Lie group. Thus $G$ is a Lie group, and the central series $G_0 = G_1 = G$, $G_2 := [G,G_1]$, $G_3 := [G,G_2]$ terminates at the third step, so that $G_3 = \{e\}$. Let $\Gamma$ be a discrete, cocompact subgroup of $G$. 

\textsc{The Heisenberg example.}  To motivate our arguments, let us first prove the above Proposition in the
model case of the Heisenberg nilmanifold
$G/\Gamma$, with
\[ G := \big\{ \left(\begin{smallmatrix}  1 & x_1 & x_3 \\  0 & 1 & x_2  \\ 0 & 0 & 1 \end{smallmatrix}\right) : x_1,x_2,x_3 \in \R\big\}\]
and
\[ \Gamma := \big\{ \left( \begin{smallmatrix}  1 & m_1 & m_3 \\  0 & 1 & m_2  \\ 0 & 0 & 1 \end{smallmatrix}\right) : m_1,m_2,m_3 \in \Z\big\}.\]
Clearly $G_1 = G$ and \[ G_2 := [G,G_1] = \big\{ \left( \begin{smallmatrix}  1 & 0 & t \\  0 & 1 & 0  \\ 0 & 0 & 1 \end{smallmatrix}\right) : t\in \R\big\}\] and $G_3 := [G,G_2] = \{I\}$.

Let us distinguish elements
\[ e_1 = \left( \begin{smallmatrix}  1 & 1 & 0 \\  0 & 1 & 0 \\ 0 & 0 & 1 \end{smallmatrix}\right), e_2 = \left( \begin{smallmatrix}  1 & 0 & 0 \\  0 & 1 & 1 \\ 0 & 0 & 1 \end{smallmatrix}\right), e_3 = \left( \begin{smallmatrix}  1 & 0 & 1 \\  0 & 1 & 0 \\ 0 & 0 & 1 \end{smallmatrix}\right).\] To these are associated the \emph{one-parameter subgroups} $(e_i^t)_{t \in \R}$:
\[ e_1^{t_1} = \left( \begin{smallmatrix}  1 & t_1 & 0 \\  0 & 1 & 0 \\ 0 & 0 & 1 \end{smallmatrix}\right), e_2^{t_2} = \left( \begin{smallmatrix}  1 & 0 & 0 \\  0 & 1 & t_2 \\ 0 & 0 & 1 \end{smallmatrix}\right), e_3^{t_3} = \left( \begin{smallmatrix}  1 & 0 & t_3 \\  0 & 1 & 0 \\ 0 & 0 & 1 \end{smallmatrix}\right).\]

 Note that 
\[ e_1^{t_1} e_2^{t_2} e_3^{t_3} = \left( \begin{smallmatrix}  1 & t_1 & t_3 + t_1 t_2 \\  0 & 1 & t_2 \\ 0 & 0 & 1 \end{smallmatrix}\right).\]
The collection $\{e_1,e_2,e_3\}$ is an example of a \emph{Mal'cev basis} for $G$ which respects $\Gamma$, the key feature to note being that $\Gamma$ is precisely the set $\{e_1^{m_1} e_2^{m_2} e_3^{m_3} : m_1,m_2,m_3 \in \Z\}$.

For Mal'cev coordinates to be of any use, we need to know how the group operation in $G$ interacts with them. It is easy to explore this for the Heisenberg nilmanifold. Every element $x = e_1^{t_1} e_2^{t_2} e_3^{t_3} \in G$ may be written in Mal'cev coordinates as $\langle t_1,t_2,t_3\rightii$. It is a simple matter to check that multiplication in $G$ is given by the rule
\begin{equation}\label{malcev-multiplication-1} \langle t_1,t_2,t_3 \rightii \ast \langle u_1, u_2, u_3 \rightii = \langle t_1 + u_1, t_2 + u_2, t_3 + u_3 - t_2 u_1\rightii.\end{equation}
A trivial induction confirms that if $g = \langle \alpha_1,\alpha_2,\alpha_3 \rightii$ then 
\begin{equation}\label{eq447} g^n = \langle n\alpha_1, n\alpha_2, n\alpha_3 - \half n(n-1) \alpha_1 \alpha_2\rightii ,\end{equation}
an expression which provides the first indication that 2-step nilmanifolds are somehow associated with ``quadratic'' types of behaviour.

 To coordinatize the nilmanifold $G/\Gamma$, we pick a fundamental domain for the action of $\Gamma$ on $G$. A very natural one is
\[ \mathcal{F} := \{ \langle x_1, x_2 , x_3 \rightii : -\half  < x_1,x_2,x_3 \leq \half \}.\] If $x = \langle x_1,x_2,x_3\rightii \in G$, then we write $\gamma_x$ for the unique element of $\Gamma$ such that $x \gamma_x \in \mathcal{F}$. We have
\[ \gamma_x = \langle -[x_1], -[x_2], -[x_3 - [x_1]x_2]\rightii,\] where $[u] = u - \{u\}$ denotes the nearest integer function (fractional parts are taken to have values in $(-\half ,\half )$). Defining
\[ \tau(x) = x\gamma_x,\]
we therefore have
\[ \tau(x) = \langle \{x_1\}, \{x_2\}, \{x_3 - [x_1]x_2\}\rightii.\]
For any element $x$ we have that $x$ and $\tau(x)$ are equivalent under the action of $\Gamma$ on $G$.

 We may now analyse the map $T_g : G/\Gamma \rightarrow G/\Gamma$. Recall that if $\psi : G \rightarrow G/\Gamma$ is the canonical projection then the transformation $T_g : G/\Gamma \rightarrow G/\Gamma$ is defined via the rule $T_g(\psi(x)) = \psi(g x)$. Persisting with the notation $g = \langle \alpha_1 ,\alpha_2,\alpha_3 \rightii$ and using coordinates on the fundamental domain $\mathcal{F}$ to represent $G/\Gamma$, we have
\begin{eqnarray}\nonumber T_g^n (0)&  = & \tau(g^n 0)\\ \nonumber & = & \langle \{n \alpha_1\}, \{ n\alpha_2 \}, \{n \alpha_3 - \half n(n-1) \alpha_1 \alpha_2 - [n\alpha_1]n\alpha_2\}\rightii \\ & \equiv & \langle n \alpha_1, n\alpha_2 , n \alpha_3 - \half n(n-1) \alpha_1 \alpha_2 - [n\alpha_1]n\alpha_2\rightii\;\; \md{1} .\label{heis-1}\end{eqnarray}
This provides the first indication that nilmanifolds encode behaviour somewhat more general than simply quadratic; here we have ``generalised'' quadratic behaviour typified by the appearance of the ``bracket quadratic'' $[n\alpha_1] n\alpha_2$. We have now assembled everything we need to prove Proposition \ref{nil-to-local} for the Heisenberg nilmanifold.

\begin{proof}[Proof of Proposition \ref{nil-to-local} for the Heisenberg nilmanifold] Let $F(T_g^n x)$ be a nilsequence on $G/\Gamma$. For the sake of exposition we take $x = 0$ so that \eqref{heis-1} applies. Let $\pi : G \rightarrow G/G_2$ be the canonical projection and, by abuse of notation, write $\pi : G/\Gamma \rightarrow G/\Gamma G_2$ for the induced projection. Now $G/\Gamma G_2$ is a $1$-step nilmanifold, being the quotient of $G/G_2$ by $\Gamma/\Gamma \cap G_2$, and we may identify it with $(\R/\Z)^2$ via the coordinatization
\[ \pi(\langle t_1, t_2, t_3\rightii) = (t_1, t_2).\]
Observe that $(\pi(T_g^n 0))_{n \in \N} = (T_{\pi(g)}^n 0)_{n \in \N}$ is an orbit on $G/\Gamma G_2$, generated by the rotation $T_{\pi(g)} : (t_1, t_2) \rightarrow (t_1 + \alpha_1 , t_2 + \alpha_2)$ on the torus.
 Let
\[ 1 = \sum_{l=1}^d \psi_l,\]$d = O(1)$, be a Lipschitz partition of unity on $(\R/\Z)^2$ with the property that for each $l$ there are $x_1,x_2$ such that 
\[ \mbox{Supp}(\psi_l) = [x_1, x_1 + \textstyle\frac{1}{10}] \times [x_2, x_2 + \frac{1}{10}].\]
 Then we have
\[ F(T^n_g 0) = \sum_{l=1}^d \psi_l (T^n_{\pi(g)} 0)) F(T^n_g 0).\]
 We will look at each constituent nilsequence $\psi_l(T^n_{\pi(g)} 0) F(T^n_g 0)$, and write it in terms of local quadratics on $1$-step Bohr sets defined on $G/\Gamma G_2$.

 Fix $l$, $1 \leq l \leq d$ together with the associated $x_1$ and $x_2$. Now the set $U := \{ x \in G/\Gamma : \pi(x) \in [x_1, x_1 + \textstyle\frac{1}{10}] \times [x_2, x_2 + \frac{1}{10}]\}$ is diffeomorphic to the direct product
\[ [x_1, x_1 + \textstyle\frac{1}{10}] \times [x_2, x_2 + \frac{1}{10}] \times \R/\Z,\] which itself is diffeomorphic to a subset of $(\R/\Z)^3$. 
Write $\pi_3 : U \rightarrow \R/\Z$ for projection onto the third coordinate. Write $S$ for the set of all $n \in \N$ such that $T_g^n 0 \in U$. Note that $S$ is a $1$-step Bohr set, since
\[ S = \{ n : \psi_l (T_{\pi(g)}^n 0) \neq 0\}.\]

\begin{lemma}[Local quadratic behaviour]\label{lqb}
Suppose that $n,h_1, h_2$ and $h_3$ are such that all eight of the points $n + \epsilon_1 h_1 + \epsilon_2 h_2 + \epsilon_3 h_3$, $\epsilon_1,\epsilon_2,\epsilon_3 \in \{0,1\}$, lie in $S$. Then the $\pi_3$-coordinates are subject to the quadratic constraint
\[ \sum_{\epsilon_1,\epsilon_2,\epsilon_3 \in \{0,1\}} (-1)^{\epsilon_1 + \epsilon_2 + \epsilon_3} \pi_3(T^{n + \epsilon_1 h_1 + \epsilon_2 h_2 + \epsilon_3 h_3}_g 0) = 0.\]
\end{lemma}

\begin{proof} Recall \eqref{heis-1}. Writing 
\[ f_1(n) := n \alpha_3 - \half n(n-1) \alpha_1 \alpha_2 - [n\alpha_1]n\alpha_2,\]
we are to show that 
\[ \sum_{\epsilon_1,\epsilon_2,\epsilon_3 \in \{0,1\}} (-1)^{\epsilon_1 + \epsilon_2 + \epsilon_3}f_1(n + \epsilon_1 h_1 + \epsilon_2 h_2 + \epsilon_3 h_3) = 0\]
whenever the $n + \epsilon_1 h_1 + \epsilon_2 h_2 + \epsilon_3 h_3$ are all in $S$. We may write $f_1$ as the sum of a quadratic polynomial and $f_2(n) := \{n \alpha_1\} n\alpha_2$. It suffices, then, to verify the result for this function $f_2$ instead. To do this, we note that the obvious relations
\[ \{ (\epsilon_1 h_1 + \epsilon_2 h_2 + \epsilon_3 h_3)\alpha_1\} \equiv \{ (n + \epsilon_1 h_1 + \epsilon_2 h_2 + \epsilon_3 h_3)\alpha_1\} - \{ n \alpha_1\} \md{1}\]
are actually equalities in $\R$, and not just in $\R/\Z$, by virtue of the constraint that all quantities $\{(n + \epsilon_1 h_1 + \epsilon_2 h_2 + \epsilon_3 h_3)\alpha_1 \}$ lie in the interval $[x_1, x_1 + \frac{1}{10}]$. Furthermore we have such relations as
\[ \{ h_1 \alpha_1 \} + \{ h_2 \alpha_1\} = \{(h_1 + h_2) \alpha_1\}.\]
By employing these together with a few simple manipulations, the lemma follows.
\end{proof}

 To introduce locally quadratic exponentials, we use Lemma \ref{fourier-lip} to approximate $F = F(u_1,u_2,u_3)$, considered as a function on $U \subseteq (\R/\Z)^3$, by a sum of exponentials. For any $\epsilon$ we may pick $J = O(\epsilon^{-3}\log^3(1/\epsilon))$ together with complex numbers $c_1,\dots,c_J = O(1)$ and frequencies $m_1,\dots, m_J \in \Z^3$ so that 
\[ F(u_1,u_2, u_3) = \sum_{j = 1}^J c_j e(m_j \cdot u) + O(\epsilon)\] for all $u = (u_1,u_2,u_3) \in U$. Using \eqref{heis-1} we obtain the formula
\[ F(T_g^n 0) = \sum_{j = 1}^J c_j e(m_j^{(1)} \{n \alpha_1 \} + m_j^{(2)} \{n \alpha_2\} + m_j^{(3)} \pi_3 (T_g^n 0)) + O(\epsilon).\] Each function $e(m_j^{(1)} \{n \alpha_1 \} + m_j^{(2)} \{n \alpha_2\})$ is a Lipschitz nilsequence on $G/\Gamma G_2$, that is to say it can be written in the form $f_{k}(T_{\pi(g)}^n 0)$. Thus we can write
\[ \psi_l(\pi(T^n_g 0)) F(T^n_g 0) = \sum_{j = 1}^J \widetilde f_{j}(T_{\pi(g)}^n 0) e(m_j^{(3)} \pi_3 (T_n^g 0)) + O(\epsilon).\]
By Lemma \ref{lqb}, each of the constituents here is a local quadratic on a 1-step Bohr set. This concludes the proof of Proposition \ref{nil-to-local} in the special case of the Heisenberg nilmanifold.
\end{proof}

\textsc{The general case.} The above arguments can be can be extended to more general nilpotent groups.  To do so, we need to involve the Lie algebra $\mathfrak{g}$ associated to $G$ together with the exponential map
\[ \exp : \mathfrak{g} \rightarrow G.\]
For the Heisenberg nilmanifold $\mathfrak{g}$ may be identified with the Lie algebra of strictly upper triangular $3 \times 3$ matrices over $\R$ with $0$'s on the diagonal, that is to say
\[ \mathfrak{g} = \left\{ \left( \begin{smallmatrix}  0 & u_1 & u_3 \\  0 & 0 & u_2  \\ 0 & 0 & 0 \end{smallmatrix}\right) : u_1,u_2,u_3 \in \R\right\}.\]
The exponential map is given by matrix exponentiation, so $\exp(X) = e^X$, which in practice means that if
\[ X =  \left( \begin{smallmatrix}  0 & u_1 & u_3 \\  0 & 0 & u_2  \\ 0 & 0 & 0 \end{smallmatrix}\right)\]
then
\[ \exp(X) =  \left( \begin{smallmatrix}  1 & u_1 & u_3 + \frac{1}{2}u_1 u_2 \\  0 & 1 & u_2  \\ 0 & 0 & 1 \end{smallmatrix}\right).\]
Wth the notation of Lie algebras and the exponential map it is possible to define, for a connected, simply-connected, nilpotent Lie group $G$, the $1$-parameter subgroup $(g^t)_{t \in \R}$ associated to an element $g \in G$. Thus we set
\[ \exp(X)^t := \exp (t X),\]
for all $X \in \mathfrak{g}$ and $t \in \R$.

We can now obtain Mal'cev coordinates for any nilmanifold arising from a connected and simply connected Lie group:

\begin{proposition}[Mal'cev coordinates of the second kind]\label{malcev-coordinates-2}
Let $G$ be a connected and simply connected $s$-step nilpotent Lie group with central series 
\[ G = G_0 = G_1 \supseteq G_2 \supseteq G_3 \supseteq \dots \supseteq G_{s+1} = \{e\}.\]
Let $\Gamma$ be a discrete, cocompact subgroup of $G$. Then there is a collection 
\[ \{ e_1,\dots, e_{i_1}, e_{i_1 + 1}, \dots, e_{i_2}, e_{i_2 + 1}, \dots, e_{i_k}\}\] such that 
\begin{itemize}
\item[(i)] Suppose that $j \in \{1,\dots,s+1\}$, and define $i_0 : = 1$. Then every element of $G_j$ can be written uniquely as $e_{i_j + 1}^{t_{i_j + 1}} \dots e_{i_{s+1}}^{t_{s+1}}$, for real numbers $t_{i_j + 1}, \dots, t_{s+1}$. 
\item[(ii)] We have
\[ \Gamma = \{e_1^{m_1} \dots e_{s+1}^{m_{s+1}} : m_1,\dots,m_{s+1} \in \Z\}.\]
\end{itemize}
\end{proposition}

It turns out to be more natural to deal with \emph{coordinates of the first kind}, which are defined on the Lie algebra $\mathfrak{g}$. Before defining these, we assemble some slightly disparate facts about how the exponential map provides a link between $\mathfrak{g}$ and $G$ in the nilpotent case. It is not particularly easy to find proofs of all of these statements in one place: our main resources were \cite{bourbaki} and \cite{corwin-greenleaf}. 

\begin{proposition}[Nilpotent Lie algebras and groups]\label{nilpotent-properties} Let $G$ be a connected, simply connected, $s$-step nilpotent Lie group. Let $\mathfrak{g}$ be the corresponding Lie algebra, and let $\exp: \mathfrak{g} \rightarrow G$ be the exponential map. We have the following statements.
\begin{itemize}
\item[(i)] $\exp$ is a diffeomorphism between $\mathfrak{g}$ and $G$, both of which are diffeomorphic to some Euclidean space $\R^d$. 
\item[(ii)] Define the central series of $\mathfrak{g}$ by $\mathfrak{g}_0 = \mathfrak{g}_1 := \mathfrak{g}$ and $\mathfrak{g}_{i+1} = [\mathfrak{g}, \mathfrak{g}_i]$ for $i \geq 1$. Then $\exp(\mathfrak{g}_i) = G_i$. In particular, the Lie algebra $\mathfrak{g}$ is $s$-step nilpotent. We have the relations $[\mathfrak{g}_i , \mathfrak{g}_j] \subseteq \mathfrak{g}_{i + j}$ and $[G_i, G_j] \subseteq G_{i + j}$. 
\item[(iii)] \textup{(Baker-Campbell-Hausdorff Formula)} We have
\[ \exp(X) \exp(Y) = \exp(Z),\]
where
\[ Z = X + Y + \frac{1}{2}[X,Y] + \frac{1}{12}[X,[X,Y]] + \frac{1}{12}[Y, [Y, X]] + \dots\]
\end{itemize}
\end{proposition}

\begin{remarks} The dots in (iii) are supposed to indicate that the Baker-Campbell-Hausdorff formula has terms involving commutators of fourth and higher order. Note, however, that since $\mathfrak{g}$ is nilpotent, the series does terminate. It is possible to give a description of the whole series, though it does not have a particularly simple closed form. See \cite{bourbaki}.
\end{remarks}

We describe now the Mal'cev coordinates of the first kind:

\begin{theorem}[Mal'cev coordinates of the first kind] Let $G$ be a connected, simply-connected, nilpotent Lie group with Mal'cev basis $\{e_1,\dots,e_k\}$. Thus any element $g \in G$ may be written uniquely as $e_1^{t_1}\dots e_k^{t_k}$, giving rise to the Mal'cev coordinates of the second kind $\langle t_1,\dots, t_k\rangle_{\mbox{\emph{\scriptsize II}}}$. Write $e_i = \exp(X_i)$, where $X_i \in \mathfrak{g}$. Then for any $g \in G$ there are unique $\xi_1,\dots,\xi_k \in \R$ such that $g = \exp(\xi_1 X_1 + \dots + \xi_k X_k)$. We refer to the elements of the $k$-tuple $\langle\xi_1,\dots,\xi_k\rangle_{\mbox{\emph{\scriptsize I}}}$ as the \emph{Mal'cev coordinates of the first kind}.\endproof
\end{theorem}

\begin{remark} In view of Proposition \ref{malcev-coordinates-2} (i) and Proposition \ref{nilpotent-properties} (ii), we have
\[ \mathfrak{g}_j = \mbox{Span}_{\R} (X_{i_j + 1}, \dots, X_{i_k}).\]
For the Heisenberg nilmanifold, note that $\langle t_1, t_2 , t_3 \rightii = \langle t_1, t_2, t_3 + \frac{1}{2}t_1 t_2\righti$. 
\end{remark}

Writing $\tau : \R^3 \rightarrow G$ for the map which identifies coordinates of the first kind with the element in $G$ they represent, we see that $\tau^{-1}(\Gamma)$ is not a lattice. Fortunately, something nearly as good is true.

\begin{proposition}[Fundamental domain description of a nilmanifold] \cite[Ch IV.6]{auslander-green-hahn}. \label{fund-domain} Let $G/\Gamma$ be a nilmanifold, and suppose that $X_1,\dots,X_k$ is a Mal'cev basis of the first kind in $\mathfrak{g}$. Let $\tau : \langle\xi_1,\dots,\xi_k\rangle_{\mbox{\emph{\scriptsize I}}} \mapsto \exp(\xi_1 X_1 + \dots + \xi_k X_k)$ be the coordinate map, and let $\mathcal{F}$ be any region of the form
\[ \{\langle\xi_1,\dots, \xi_k\rangle_{\mbox{\emph{\scriptsize I}}} : a_i \leq \xi_i < a_i + 1 \; \; \mbox{for all $i$}\}. \] Then each point of $G$ is equivalent, under the right action of $\Gamma$, to precisely one point in $\exp(\mathcal{F})$. Furthermore the natural projection map $\pi : G \rightarrow G/\Gamma$ is continuous on $\exp(\mathcal{F})$ and is a homeomorphism when restricted to the interior $\exp(\mathcal{F})^{\circ}$.
\end{proposition}

Our aim now is to describe the action of some $g = \langle \beta_1,\dots,\beta_k\righti$ on $G/\Gamma$ by finding formul{\ae} analogous to \eqref{malcev-multiplication-1}, \eqref{eq447} and \eqref{heis-1}. The key tool is the Baker-Campbell-Hausdorff formula. For notational simplicity we restrict to the $2$-step case from now on, and write $m := i_2$ and $n := i_3$. Thus the Mal'cev basis of the first kind for $G$ is $\{X_1,\dots,X_m,X_{m+1}, \dots, X_n\}$, where
\[ \mbox{Span}_{\R}(X_{m+1}, \dots, X_n) = \mathfrak{g}_2 = [\mathfrak{g}, \mathfrak{g}].\]
The Lie algebra $\mathfrak{g}$ is completely specified by its \emph{structure constants}, a collection of real numbers $(a_{ijk})_{1 \leq i,j \leq m, m+1 \leq k \leq n}$ such that 
\begin{equation}\label{structure-constants} [X_i, X_j] = \sum_{k = m+1}^n a_{ijk} X_k.\end{equation}
These constants can be arbitrary so long as $(a_{ijk})_{i,j \leq m}$ is antisymmetric for each $k$, though if we want $G$ to possess a cocompact subgroup $\Gamma$ then certain rationality conditions must hold \cite{malcev}.

\begin{lemma}[Multiplication in coordinates of the first kind]\label{mult-lem}
Suppose that $G$ is a connected and simply-connected $2$-step nilpotent Lie group with group operation $*$, and abuse notation by identifying elements of $G$ with their coordinates of the first kind. Then we have
\begin{align} \nonumber & \langle \xi_1, \dots, \xi_n\righti \ast \langle \nu_1, \dots, \nu_n \righti \\ = & \langle \xi_1 + \nu_1 , \dots, \xi_m + \nu_m, \xi_{m+1} + \nu_{m+1} + \phi_{m+1}(\xi_{\leq m}, \nu_{\leq m}), \dots, \xi_n + \nu_n + \phi_n(\xi_{\leq m}, \nu_{\leq m})\righti, \nonumber \end{align}
where the $\xi_{\leq m} := (\xi_1,\dots, \xi_m)$, $\nu_{\leq m} := (\nu_1, \dots, \nu_m)$ and the $\phi_j$ are antisymmetric bilinear forms.
\end{lemma}
\proof This is a simple matter of combining the Baker-Campbell-Hausdorff formula with the existence of structure constants \eqref{structure-constants}. We remark that the presentation of a $2$-step nilmanifold in this form is essentially the same as an example discussed by Furstenberg in \cite{furstenberg-vonneumann}.

Observe in particular that 
\begin{equation}\label{eq447-b} g^n = \langle n\beta_1, \dots, n\beta_n\righti,\end{equation}
and thus
\begin{equation}\label{eq447-c}
T^n_g x = \langle n \beta_1 + x_1, \dots, n \beta_m + x_m, n\beta'_1 + x_{m+1}, \dots, n \beta'_n + x_n\righti
\end{equation}
for certain constants $\beta'_j$ depending on $g, x$ and the bilinear forms $\phi_j$.

 To coordinatize $G/\Gamma$ we pick, in view of Proposition \ref{fund-domain}, the very natural fundamental domain
\[ \mathcal{F} := \{ \langle x_1, \dots, x_n \righti : -\half  < x_1,\dots, x_n \leq \half \}.\] If $x = \langle x_1,\dots,x_n \righti \in G$, then we write $\gamma_x$ for the unique element of $\Gamma$ such that $x \gamma_x \in \mathcal{F}$. Write $\tau(x) = x\gamma_x$. We need a formula for $\gamma_x$ in terms of coordinates of the first kind, and to obtain such a result we need a description of the lattice $\Gamma$ in terms of these coordinates. Since $\Gamma$ may be identified with $\Z^n$ in coordinates of the second kind, such a description can be obtained by finding the relation between the two types of coordinate. Such a relation is easy to obtain. Indeed by definition we have
\[ \langle t_1, \dots, t_n \rightii = \langle t_1,0,\dots,0\righti \ast \dots \ast \langle 0,\dots, 0, t_n\righti.\]
By inductive use of Lemma \ref{mult-lem} this quickly implies that 
\begin{equation}\label{i-ii-relation}
 \langle t_1, \dots, t_n \rightii = \langle t_1,\dots, t_m, q_{m+1}(t_{\leq m}), \dots, q_{n}(t_{\leq m})\righti\end{equation} for certain quadratic forms $q_j$. In fact these forms are rather related to the alternating forms $\psi_j$; if $\psi(x,y) = \sum_{k,l \leq m} a_{kl} x_l y_k$ then $q(x) = \sum_{k < l} a_{kl}x_k x_l$. 

 In terms of coordinates of the first kind, then, we see that 
\[ \Gamma = \{ \langle r_1,\dots, r_m, r_{m+1} + q_{m+1}(r_{\leq m}), \dots,r_n + q_n(r_{\leq m}) : r_1,\dots, r_n \in \Z\}.\]
It follows that 
\begin{align} \nonumber
\gamma_x = \langle -[x_1], \dots, -[x_m],&  -[x_{m+1} - \phi_{m+1}(x_{\leq m}, [x]_{\leq m}) + q_{m+1}([x]_{\leq m})], \\ & \dots, -[x_n - \phi_n(x_{\leq m}, [x]_{\leq m}) + q_n([x]_{\leq m}) \righti \nonumber 
\end{align}
and that 
\begin{align}\nonumber \tau(x) = \langle \{x_1\}, \dots, \{x_m\}, & \{x_{m+1} - \phi_{m+1}(x_{\leq m}, [x]_{\leq m}) + q_{m+1}([x]_{\leq m})\} ,\\ & \dots, \{x_n - \phi_n(x_{\leq m}, [x]_{\leq m}) + q_n([x]_{\leq m})\} \righti.\nonumber\end{align}
We remark that we have essentially provided an independent confirmation of Proposition \ref{fund-domain} for 2-step nilmanifolds. The proof in the $s$-step case merely involves more notation.

 Combining this with \eqref{eq447-c} leads to the analogue of \eqref{heis-1}:
\[ T_g^n x \equiv \langle  n\beta_1, \dots, n \beta_m, \psi_{m+1}(n), \dots, \psi_n(n)\righti \; \; \md{1},\] 
where each $\psi_j$ has the form
\[ \psi(n) = an + b + \sum_{i = 1}^m c_i n [n\beta_i] + \sum_{l < k \leq m} c_{lk} \{ n \beta_l \} \{ n \beta_k\}.\] 
The remainder of the proof of Proposition \ref{nil-to-local} is, from this point, almost identical to the special case of the Heisenberg nilmanifold. We leave the details to the reader.\endproof

\section{Divisor moment estimates}\label{div-section}

We collect some standard moment estimates for the divisor function $\tau(n) := \sum_{d|n} 1$.  These are used to prove Proposition \ref{div-pack}, which is used in \S \ref{massage-sec} to show that there are not too many ``collisions'' occuring in sets such as $\{ dw: D < d \leq 2D; W < w \leq 2W \}$.

The basic estimate we need is

\begin{lemma}\label{lem5.1} Let $m, N \geq 1$ be integers. Then
we have the moment estimate \[ \E_{n \in [N]} \tau(n)^m \ll_m (\log N)^{2^m - 1}.\]
\end{lemma}

\begin{proof} 
This is very standard: see, for example, \cite{bourgain2} or \cite{ruzsa}.  For our application, the precise value of exponent $2^m-1$ does not need to
be attained; any bound of the form $\log^{C_m} N$ would suffice.
\end{proof}

 In particular, we have the second moment estimate
$$ \E_{n \in [N]} \tau(n)^2 \ll \log^3 N$$
which by dyadic decomposition then implies
\begin{equation}\label{taud}
 \sum_{n \in [N]} \frac{\tau(n)^2}{n} \ll \log^4 N.
\end{equation}

 Now if $A \subseteq \{1,\ldots,N\}$ is a nonempty set of size $\alpha N$ and $m \geq 2$ is an integer, then from H\"older's inequality we have
\begin{align*}
\E_{n \in A} \tau(n)^2 &\leq (\E_{n \in A} \tau(n)^{m})^{2/m} \\
&\leq \alpha^{-2/m} (\E_{n \in [N]} \tau(n)^m)^{2/m}) \\
&\ll_m \alpha^{-2/m} (\log N)^{2(2^m-1)/m}.
\end{align*}
In particular, for any $\kappa < 1/2$ we have the moment estimate
\begin{equation}\label{tau-rider}
\E_{n \in A} \tau(n)^2 \ll_\kappa \alpha^{-\kappa} \log^{2^{2/\kappa}} N.
\end{equation}
 This estimate has the following consequence.

\begin{lemma}[Divisor packing lemma]\label{div-pack}  Let $A \subseteq \{1,\ldots,N\}$ be a non-empty set of size $\alpha N$, and for each $d \geq 1$ let $A_d := \{ n \in A: d|n\}$
denote those elements of $A$ which are multiples of $d$.  Suppose ${\mathfrak D} \subset \Z^+$ is a finite set of positive integers such that
$$ |A_d| \geq \delta |A|$$
for all $d \in {\mathfrak D}$ and some $\delta > 0$.  Then for any positive $\kappa < 1/2$ we have
$$ |\bigcup_{d \in {\mathfrak D}} A_d| \gg_\kappa \delta^2 |{\mathfrak D}|^2 |A| \alpha^{\kappa} \log^{-2^{2/\kappa}} N.$$
\end{lemma}

\begin{proof}  From hypothesis we have
$$ \E_{n \in A} \sum_{d \in {\mathfrak D}} 1_{A_d}(n) = \sum_{d \in {\mathfrak D}} \frac{|A_d|}{|A|} \geq \delta |{\mathfrak D}|.$$
By Cauchy-Schwarz we conclude that
$$ \frac{|\bigcup_{d \in {\mathfrak D}} A_d|}{|A|} \E_{n \in A} (\sum_{d \in {\mathfrak D}} 1_{A_d}(n))^2 \geq \delta^2 |{\mathfrak D}|^2.$$
From the trivial bound
$$ \sum_{d \in {\mathfrak D}} 1_{A_d}(n) \leq \sum_{d|n} 1 = \tau(n)$$
and \eqref{tau-rider} we thus have
$$ \frac{|\bigcup_{d \in {\mathfrak D}} A_d|}{|A|} \alpha^{-\kappa} \log^{2^{2/\kappa}} N \gg_\kappa \delta^2 |{\mathfrak D}|^2$$
and the claim follows.
\end{proof}


\begin{thebibliography}{99}



\bibitem{auslander-green-hahn} L.~Auslander, L.~Green and F.~Hahn, \emph{Flows on Homogeneous spaces,} Annals of Math. Studies \textbf{53} (1963).

\bibitem{baker-harman} 
R.~C.~Baker and G.~ Harman, \emph{Exponential sums formed with the M\"obius function,}
J. London Math. Soc. (2) \textbf{43} (1991), no. 2, 193--198.

\bibitem{bilu} Y.~Bilu, \emph{Structure of sets with small sumset,} Structure theory of set addition.  Ast\'erisque  \textbf{258} (1999), xi, 77--108.

\bibitem{bourbaki} N.~Bourbaki, \emph{Lie groups and Lie algebras,} Chapters 1--3. Translated from the French. Reprint of the 1989 English translation. Elements of Mathematics (Berlin). Springer-Verlag, Berlin, 1998. xviii+450 pp.


\bibitem{bourgain2} J.~Bourgain, \emph{On $\Lambda(p)$-subsets of squares,} Israel J. Math. \textbf{67} (1989), no. 3, 291--311.

\bibitem{corwin-greenleaf} L.~Corwin and F.~P.~Greenleaf, \emph{Representations of nilpotent Lie groups and their applications, Part I: Basic theory and examples,} Cambridge Advanced Studies in Math. \textbf{18}, CUP 1990.

\bibitem{davenport-early} H. Davenport, \emph{On some infinite series involving arithmetical functions. II}, Quart. J. Math. Oxf. \textbf{8} (1937), 313--320 
 
\bibitem{davenport} H.~Davenport, \emph{Multiplicative number theory,} Third edition. Graduate Texts in Mathematics, \textbf{74}. Springer-Verlag, New York, 2000. xiv+177 pp

\bibitem{furstenberg-vonneumann} H.~Furstenberg, \emph{Nonconventional ergodic averages,} in The legacy of John von Neumann (Hempstead, NY, 1988),  43--56, Proc. Sympos. Pure Math., 50, Amer. Math. Soc., Providence, RI, 1990. 

\bibitem{gowers-long-aps} W.~T.~Gowers, \emph{A new proof of Szemer\'edi's Theorem,} Geom. Funct. Anal. \textbf{11} (2001), no. 3, 465--588.

\bibitem{green-fin-field} B.~J.~Green, \emph{Finite field models in additive combinatorics}, Surveys in Combinatorics 2005, London Math. Soc. Lecture Notes \textbf{327}, 1--27.

\bibitem{green-tao-primes} B.~J.~Green and T.~C.~Tao, \emph{The primes contain arbitrarily long arithmetic progressions,} to appear, Annals of Math.

\bibitem{green-tao-inverseu3} B.~J.~Green and T.~C.~Tao, \emph{An inverse theorem for the Gowers $U^3$-norm,} to appear, Proc. Edinburgh Math. Soc.

\bibitem{green-tao-linearprimes} B.~J.~Green and T.~C.~Tao, \emph{Linear equations in primes,} to appear, Annals of Math.

\bibitem{hua-1} L.~K.~Hua, \emph{Some results in the additive prime number theory,} Quart. J. Math. Oxford \textbf{9} (1938), 68--80.

\bibitem{iwaniec-kowalski} H.~Iwaniec and E.~Kowalski, 
\emph{Analytic number theory,} 
American Mathematical Society Colloquium Publications, \textbf{53}. 
American Mathematical Society, Providence, RI, 2004. xii+615 pp

\bibitem{landau}
E.~Landau, \emph{Handbuch der Lehre von der Verteilung der Primzahlen.}, Leipzig, Germany: Teubner, 1909. 

\bibitem{malcev}
A. Mal'cev, \emph{On a class of homogeneous spaces}, Izvestiya Akad. Nauk SSSR, Ser Mat. \textbf{13} (1949), 9--32.

\bibitem{montgomery} H.~L.~Montgomery, \emph{Ten lectures on the interface between analytic number theory and harmonic analysis,}
CBMS Regional Conference Series in Mathematics, \textbf{84}. 
Published for the Conference Board of the Mathematical Sciences, Washington, DC by the American Mathematical Society, Providence, RI, 1994. xiv+220 pp.

\bibitem{ruzsa} I.~Z.~Ruzsa, \emph{On an additive property of squares and primes,} Acta Arith. \textbf{49} (1988), no. 3, 281--289.

\bibitem{tao:survey}
T.~C.~Tao, \emph{Arithmetic progressions in the primes}, Collectanea Mathematica (2006), Vol. Extra., 37--88. [Proceedings, 7th International Conference on Harmonic Analysis and Partial Differential Equations.]

\bibitem{tao-vu-book} T.~C.~Tao and V.~H.~Vu, \emph{Additive combinatorics,} CUP 2006.

\bibitem{vaughan} R.~C.~Vaughan, \emph{Sommes trigonom\'etriques sur les nombres premiers,} 
C. R. Acad. Sci. Paris S\'er. A-B \textbf{285} (1977), no. 16, A981--A983.

\bibitem{vaughan-hlm} R.~C.~ Vaughan, \emph{The Hardy-Littlewood method,} Second edition. Cambridge Tracts in Mathematics, \textbf{125}. Cambridge University Press, Cambridge, 1997. xiv+232 pp

\bibitem{van-der-corput}
J.G. van der Corput, \emph{\"Uber Summen von Primzahlen und Primzahlquadraten,} Math. Ann. \textbf{116} (1939), 1--50.

\bibitem{vinogradov-1} I.~M.~Vinogradov, \emph{Some theorems concerning the primes,} Mat. Sbornik. N.S. \textbf{2} (1937), 179--195.



\end{thebibliography}
\end{document}